\documentclass[a4paper,12pt]{article}
\usepackage{amsmath,latexsym,amssymb,amsfonts,natbib}
\usepackage[a4paper,hmargin=1.87 cm]{geometry}
\usepackage[dvips]{graphicx}
\usepackage{pifont}
\usepackage{pstcol,pst-fill,pst-grad}
\usepackage{rotating}
\usepackage{subfigure}
\usepackage{graphics}
\usepackage{hyperref}
\usepackage{multirow}
\usepackage{multicol}
\usepackage{courier}
\usepackage{pstricks}
\usepackage[latin1]{inputenc}
\usepackage{indentfirst}
\usepackage{ccaption}
\usepackage{verbatim}
\usepackage{makeidx}
\usepackage{listings}
\usepackage{pgf,tikz}
\usetikzlibrary{decorations}
\usetikzlibrary{decorations.pathmorphing}
\usetikzlibrary{decorations.pathreplacing}
\usetikzlibrary{decorations.shapes}
\usetikzlibrary{decorations.text}
\usetikzlibrary{decorations.markings}
\usetikzlibrary{decorations.fractals}
\usetikzlibrary{decorations.footprints}

\usepackage{graphicx}
\usepackage{pgf, tikz}
\usepackage{algorithm}
\usepackage{algorithmic}
\usepackage{array}
\usepackage{dsfont}

\newcommand{\proof}     {\paragraph{Proof}}
\newcommand{\carre}     {\hfill$\Box$}
\numberwithin{equation}{section}
\newtheorem{defi}{Definition}
\newtheorem{lem}{Lemma}
\newtheorem{theo}{Theorem}

\newtheorem{prop}{Proposition}
\newtheorem{rem}{Remark}
\newtheorem{hyp}{Hypothesis}

\newcommand{\N}{\mathbb N}

\title{Asymptotic adaptive threshold for connectivity in a  random geometric social network}
\author{Ahmed Sid-Ali and Khader Khadraoui \\
{\small Laval University, Department of Mathematics and Statistics}\\
{\small Qu\'ebec city G1V 0A6, Canada}\\
}
\date{}
\begin{document}
\maketitle
\begin{abstract}
Consider a dynamic random geometric social network identified by $s_t$ independent points $x_t^1,\ldots,x_t^{s_t}$ in the unit square $[0,1]^2$ that 
interact in continuous time $t\geq 0$. The generative model of the random points is a Poisson point measures. Each point $x_t^i$ can be active or not in the 
network with a Bernoulli probability $p$. Each pair being connected by affinity thanks to a step connection function if the interpoint distance 
$\|x_t^i-x_t^j\|\leq a_\mathsf{f}^\star$ for any $i\neq j$. We prove that when $a_\mathsf{f}^\star=\sqrt{\frac{(s_t)^{l-1}}{p\pi}}$ for $l\in(0,1)$, the 
number of isolated points is governed by a Poisson approximation as $s_t\to\infty$. This offers a natural threshold for the construction of a  
$a_\mathsf{f}^\star$-neighborhood  procedure tailored to the dynamic clustering of the network adaptively from the data.
\end{abstract}

\section{Introduction}
\label{intro}

A problem common to many disciplines is that of adequately studying complex networks. Research on this problem occurs in applied mathematics (dynamic systems 
governed by ODEs and PDEs), probability and statistics (random graph models) and in computer science and engineering (statistical learning networks). A less familiar 
studying in this area is stochastic individual based models (IBM). The theory of IBM is a quickly growing interdisciplinary area with a very broad spectrum of 
motivations and applications. In consideration of complex networks, we characterized recently in \cite{A+K} social network systems by  individuality of components and 
localized interaction mechanisms among components realized via  densities dependent spatial state which leads to a globally regular behavior of the system. 
In modern stochastic (with diffusion or point processes) and infinite dimensional modeling, interacting particle models for the study of complex 
systems forms a rich and powerful direction. For instance, interacting particle systems are widely used to model population biology \citep{Kh15}, 
ecology \citep{Fink+al09}, condensed matter physics, chemical kinetics, sociology and economics (agent based models). For an account on the subject of interacting 
particle systems, we refer the reader to the book \cite{Ligg85}.

We study the following random geometric social network $(r_t)_{t\geq 0}$ presumed to be described by the punctual measure 
\begin{align*}
 r_t(dx)=\sum_{i=1}^{s_t}\delta_{x_{t}^i}(dx),
\end{align*}
where $\delta_x$ denotes the Dirac measure centered on $x$ with $x$ in some measurable space (defined in a rigorous way afterwards), $s_t$ denotes the 
size of the network at time $t$ and the independent points $x_t^1,\ldots,x_t^{s_t}$ denote the spatial positions of the network members at time $t$. 
The generative dynamic model of the random points is a Poisson point measures that we describe here its  infinitesimal construction for completeness. As 
known, the main issue in this kind of study is the connection between points. In particular, we choose some deterministic rule where two members interact 
if and only if their distance does not exceed some threshold  which is precisely what has been done in random geometric graphs \citep{Gilb61}.  We refer the 
reader to \cite{Pern03} for a thorough presentation of the many  properties of Gilbert's graph. However, the two main differences in our network compared to 
Gilbert's graph are: $(i)$ the state of our network is dynamic in continuous time and not static; $(ii)$ the density of the random points is not uniform 
in $[0,1]^d$ with $d\geq 2$ as usually assumed by authors (see for instance \cite{Brout+al16,Brout+al14} and references therein).

It is known that the main obstacle to connectivity inside certain networks is the existence of isolated members. In the present paper we prove that 
the number of isolated points has asymptotic Poisson distribution by employing a concentration bound on Poisson random variables together with a tailored 
discretization method of the unit square $[0,1]^2$. We establish the optimal threshold associated with this approximation and show that this result holds over 
connection functions that are zero beyond the threshold. This was previously known for the random geometric graph and its variants where points are 
independently and uniformly distributed \citep{Det+Hen89}. Hence, our result may be seen as a generalization of this Poisson approximation for the number of 
isolated points to random points generated from Poisson point processes which encompass the uniform case.  To our best knowledge, this is the first result 
in this direction and it open a door to other complex problems such as the extension to higher dimensions  or other connection functions (Rayleigh fading 
functions, functions that decay exponentially in some power of distance, etc). Moreover, we take advantage of the asymptotic distribution of isolated points to built a data adaptive 
dynamic statistical cluster method tailored to our network. This conceivable strategy enables us to detect groups and isolated 
members at each time from the dynamic of the network.  In particular, an important question investigated recently  which is the 
community detection inside the network \citep[and the references therein]{Zhao+al12,Cast+Ver14,Jin15,Bick+al15}. In the Erd\"os-Rényi graph context, 
the stochastic block model (SBM) \citep{Hol+Pac83} is usually used to model communities where the probability of a connection occurring between two members
depends solely on their community membership. There are many extensions of the SBM for various applications, including the biological, communication and 
social networks, for instance in \cite{Bick+al09}, \cite{Sni+al97} and \cite{Park+al12}.

The outline of this paper is as follows: In Section \ref{sec:1}, we describe in detail the generative model of the network data by giving an  
explicit representation and the exact Monte Carlo scheme for computation.  In Section  \ref{sec:2}, we establish the asymptotic adaptive 
threshold needed for the Poisson approximation law of the number of isolated members inside the network. Section \ref{sec:3} is devoted to the tailored 
dynamic clustering procedure using the optimal threshold established in this paper in order to detect communities and isolated members in the network. 
Moreover, we present some numerical simulations. We then discuss our results and the outlook in Section \ref{sec:4}. Section \ref{proofs} contains the 
proofs of the main results of the paper.

\section{Preliminaries: Generation of the point sets}
\label{sec:1}

We introduce the notation and the generative model for the random social network as an interacting particle system. The spatio-temporal paradigm retained here is 
represented by a random dynamic for the network in terms of its instantaneous size together with the spatial (geometric) patterns of the members. In a rigorous 
sense, the  system of particles considered is a Markov process with values in a space of punctual 
measures and where each member of the network is tracked through time. Basically, we construct a geometric virtual space that is a closure $\bar{\mathcal{D}}$ of 
an open connected subset $\mathcal{D}$ of $\mathbb{R}^d$, for some $d\geq 1$. We represent each member of the network located at a virtual state $x$ as a Dirac 
measure $\delta_x$. The idea behind the use of  virtual space is for simple managing  of  interactions between members which is nothing beyond some distance. 
Roughly speaking, closer the members are in the virtual space bigger is  their affinity (friendship) mechanism. We 
assume that each member may invite another individual at a given rate. When a new individual is arrived, it immediately disperses from the member who invited and 
becomes a member of the network. We also assume that members are subject to departures. That  is, each member quits at a rate that depends on the local 
network state. Thus, in this social network new members could arrive at continuous time and become members as well as the members can leave the network at any time.

Formally, let denote by $\mathcal{S}_F(\bar{\mathcal{D}})$  the set of finite nonnegative measures  on 
$\bar{\mathcal{D}}$ and $\mathcal{S}\subset\mathcal{S}_F(\bar{\mathcal{D}})$ that consists of all finite point measures on $\bar{\mathcal{D}}$:
\begin{align*}
 \mathcal{S}=\Big\{\sum_{i=1}^s \delta_{x^i},\;s\geq0,\;x^i\in\bar{\mathcal{D}}    \Big\},
\end{align*} 
where the states $x^1,\ldots,x^s$  of members represent their spatial locations in the space $\bar{\mathcal{D}}$ and $s\in\mathbb{N}$ stands for the size of the 
network.  According to our previous description the network that we denote by $r_t$ at time $t\geq 0$ is characterized by the distribution 
of the members at any given time  inside the virtual space and is given by
\begin{align}
 r_t(dx)=\sum_{i=1}^{s_t}\delta_{x_{t}^i}(dx),
\label{net}
\end{align}
where the stochastic process $r_t$ take its values in  $\mathcal{S}$ and the indexes $i$ in the positions $\{x_{t}^i\}_{i=1,\ldots,s_{t}}$ are ordered here 
from an arbitrary order point of view. Before giving an explicit description of the process $(r_t )_{t\geq0}$ we introduce now the heuristics of 
the network.

\subsection{Heuristics of the network}

The  dynamic of the random social network considered here can be roughly summarized by  endowing the system with three events as follows:
\begin{enumerate}
\item a recruitment by invitation event: A member located at the virtual position $x$ in the network could send an invitation to another individual in the 
outside of the network in order to join the network community. Then, the individual can accept the invitation and joins the network at the virtual position 
$y=x+z$ where $z$ is chosen randomly following a given dispersion kernel;
\item a departure from the network event: Each member can leave the network at any moment and its position inside the virtual space becomes empty immediately; 
\item  a recruitment by affinity with the network event: An individual outside the network can be interesting to join the network thanks to a certain 
affinity with the network. Then, this individual choose a position in the network following a given dispersion kernel.
\end{enumerate}

We shall describe the system by the evolution in time of the measure $r_t$. For this, let define the parameters of the previous events: 
\begin{itemize}
\item[(i)] $v_r\in[0,\infty)$ denotes the invitation rate for each member in the network;
\item[(ii)] $K(x,dz)$ denotes the dispersion law for the new individual invited by a member located at $x$ and it is assumed to satisfy, for each $x\in\bar{\mathcal{D}}$,
\begin{align*}
 \int_{z\in\mathbb{R}^d, x+z\in \bar{\mathcal{D}}} K(x,dz)= 1 \quad \textrm{ and } \quad \int_{z\in\mathbb{R}, x+z\notin \bar{\mathcal{D}}} K(x,dz)= 0;
\end{align*}
 \item[(iii)] $d_r\in[0,\infty)$ denotes a departure rate for each individual at some $x\in\bar{\mathcal{D}}$; 
 \item[(iv)] $K^{\mathsf{af}}(dy)$ denotes the affinity dispersion law for the new arrived individual at some $y\in\bar{\mathcal{D}}$ and it is assumed to 
 satisfy,
 \begin{align*}
 \int_{y\in \bar{\mathcal{D}}} K^{\mathsf{af}}(dy)=1 \quad \textrm{ and } \quad \int_{y\notin \bar{\mathcal{D}}} K^{\mathsf{af}}(dy)= 0;
\end{align*}
\item[(v)] for all $x,y\in\bar{\mathcal{D}}$, $\mathsf{aff}(x,y)=\mathsf{aff}(y,x)\in[0,\infty)$ is the affinity kernel which describes the strength of 
 affinity between members located at $x$ and $y$.
\end{itemize}
To explain in more details the concept of affinity and its spatial dependence, we introduce  another function $w^{\mathsf{af}}$ which describes the affinity 
that may have an individual in the outside with the network; This individual may choose randomly a position $y$ for its recruitment by the network in 
function of the current member localizations  in the neighborhood of  this position $y$. To this end and for the sake of simplicity, we propose to sum the 
local affinities between the position $y$ and all the members around the future position $y$ such that: 
\begin{align*}
w^{\mathsf{af}}(y,r)=\sum_{x\in r}\mathsf{aff}(x,y)=\int_{\bar{\mathcal{D}}}\mathsf{aff}(x,y) r(dx), 
\end{align*}
for all $ x,y\in\bar{\mathcal{D}}$ which leads to $w^{\mathsf{af}}(y,r)\in[0,\infty)$. Concerning the local affinity we consider the following indicator function:  
 \begin{align}
  \mathsf{aff}(x,y)= A_\mathsf{f} \mathds{1}_{\{\|x-y\|\leq a_\mathsf{f}\}}(x,y),
  \label{localaff}
 \end{align}
where $\mathds{1}_{A}$ denotes the indicator function of the set $A$,  $A_{\mathsf{f}}\in[0,\infty)$ denotes  the maximum  that  interaction by affinity 
between two members can reach  and $a_{\mathsf{f}}\in[0,\infty)$ denotes the radius (affinity threshold) of the zone of interaction by affinity. It is easily 
seen that the affinity between $x$ and $y$ is equal to zero when the Euclidean distance $\| x-y\|> a_\mathsf{f}$. The affinity presented here generalizes 
without difficulty to more complex connection scenarios  provided that the affinity rate  investigated  is bounded.

We assume that all of these basic mechanisms are independent. Furthermore, to avoid explosion phenomenon, we suppose that the spatial dependence of the 
introduced kernels and rates is bounded. Indeed, we assume that the dispersion kernels induce  densities with respect to  the Lebesgue measure such that 
$K(x,dz)=k(x,z)dz$ and $K^{\mathsf{af}}(dy)=k^{\mathsf{af}}(y)dy$. We assume that there exists some positive  reals $\gamma_1>0$ and $\gamma_2>0$ and two 
probability densities $\tilde{k}$ on $\mathbb{R}^d$ and $\tilde{k}^{\mathsf{af}}$ on $\bar{\mathcal{D}}$ such that, for all $x\in\bar{\mathcal{D}}$,
\begin{align*}
 k(x,z)\leq \gamma_1 \tilde{k}(z) \qquad \textrm{ and } \qquad  k^{\mathsf{af}}(y)\leq \gamma_2 \tilde{k}^{\mathsf{af}}(y).  
\end{align*}
We assume that there exists a constant 
$A_\mathsf{f}$ such that, for all $x,y\in\bar{\mathcal{D}}$ and for all $r\in\mathcal{S}$,
\begin{align*}
\mathsf{aff}(x,y)\leq A_\mathsf{f} \qquad \textrm{ then } \qquad  w^{\mathsf{af}}(y,r)\leq A_\mathsf{f} s.  
\end{align*}

Now, we focus on the pathwise description of the  $\mathcal{S}$-valued stochastic process $(r_t)_{t\geq 0}$. After assuming that the spatial dependence of all the 
parameters is bounded in some sense, we shall introduce now the Poisson point measures associated with the process $r_t$.

\subsection{Explicit representation of the network}

First we derive an explicit representation of a function of the random network  for a large class of functions. Second we deduce the explicit expression of the infinitesimal generator 
on that class. To obtain this representation for the random social network $(r_t)_{t\geq 0}$ we introduce some Poisson random measures which manage the 
recruitment of new members (by invitation and affinity) and the departure of members. To this end, let $(\Omega,\mathcal{A},\mathbb{P})$ be a 
sufficiently large probability space and let consider three  punctual Poisson random measures $N_{v_r}(d\tau,di,dz,d\alpha)$ defined on 
$[0,\infty)\times\mathbb{N}^*\times\mathbb{R}^d\times[0,1]$, $N_{d_r}(d\tau,di,d\alpha)$ defined on $[0,\infty)\times\mathbb{N}^*\times[0,1]$  and 
$N_{\mathsf{af}}(d\tau,di,dy,d\alpha)$ defined on $[0,\infty)\times\mathbb{N}^*\times\bar{\mathcal{D}}\times[0,1]$ with respective intensity measures: 
\begin{align*}
 n_{v_r}(d\tau,di,dz,d\alpha)&=v_r \gamma_1\tilde{k}(z)d\tau didzd\alpha,\\
 n_{d_r}(d\tau,di,d\alpha)&=d_r d\tau did\alpha,\\
 n_{\mathsf{af}}(d\tau,di,dy,d\alpha)&= A_{\mathsf{f}}\gamma_2 \tilde{k}^{\mathsf{af}}(y) d\tau didyd\alpha,
\end{align*}
where $d\tau,dz,dy$ and $d\alpha$ are the Lebesgue measures on $[0,\infty),\mathbb{R}^d,\bar{\mathcal{D}}$ and $[0,1]$ and $di$ 
is the counting measure on $\mathbb{N}^*$. 

Let  denote by $r_0$ the initial condition of the process, it is a random variable with values in $\mathcal{S}$. Suppose that $N_{v_r}$, $N_{d_r}$, $N_{d_r}$ and $r_0$ 
are mutually independent. We also consider the canonical filtration $(\mathcal{F}_t)_{t \geq0}$ generated by  the random objects $N_{v_r}$, $N_{d_r}$, $N_{d_r}$ and 
$r_0$. In the following, we denote by $r_{t-}$ the network process at time $t$  before any possible event. The stochastic process $(r_t )_{t\geq0}$ features a jump 
dynamics (recruitments and departures). We can therefore recall a well-known formula for the pure jump processes:
\begin{align}
 \Phi(r_t)=\Phi(r_0)+\sum_{\tau\leq t}[\Phi(r_{\tau-}+\{r_\tau-r_{\tau-}\})-\Phi(r_{\tau-})], \qquad \textrm{ a.s. for }t\geq 0, 
 \label{jump}
\end{align}
for any function $\Phi$ defined on $\mathcal{S}_F(\bar{\mathcal{D}})$ for all $r\in\mathcal{S}$. We remark that in equation (\ref{jump}) the sum $\sum_{\tau\leq t}$ contains 
only a finite number of terms as the network process $(r_t )_{t\in[0,T]}$ admits only a finite number of jumps for any $T<\infty$. Note that the 
number of jumps in the network $r_t$ is bounded by a linear recruitment and departure process with arrival rate $v_r+ A_\mathsf{f}$ and departure rate $d_r$ 
\citep{Allen03}. According to the formula (\ref{jump}), for any function $\Phi$ on 
$\mathcal{S}_F(\bar{\mathcal{D}})$, we can write
\begin{align}
\begin{split}
 \Phi(r_t)=\Phi(r_0)&+\int_{0}^t\int_{\N^*}\int_{\mathbb{R}^d}\int_0^1  \mathds{1}_{\{i\leq s_{\tau-}\}}
 \mathds{1}_{\{\alpha\leq ( k(x_{\tau-}^i,z))/(\gamma_1\tilde{k}(z)) \}} \\
&\qquad\qquad \times[\Phi(r_{\tau-}+\{\delta_{(x_{\tau-}^i+z)}\})-\Phi(r_{\tau-})]  N_{v_r}(d\tau,di,dz,d\alpha) \\
&+\int_{0}^t\int_{\N^{\ast}}\int_0^1 \mathds{1}_{\{i\leq s_{\tau-}\}}
 \mathds{1}_{\{\alpha\leq(d_r/d_r)\}}  \\
&\qquad\qquad \times[\Phi(r_{\tau-}-\{\delta_{(x_{\tau-}^i)}\})-\Phi(r_{\tau-})]  N_{d_r}(d\tau,di,d\alpha) \\
&+\int_{0}^t \int_{\N^{\ast}}\int_{\bar{\mathcal{D}}}\int_0^1 \mathds{1}_{\{i\leq s_{\tau-}\}}
 \mathds{1}_{\{\alpha\leq (\mathsf{aff}(x_{\tau-}^i,y)k^{\mathsf{af}}(y))/(A_\mathsf{f}\gamma_2 \tilde{k}^{\mathsf{af}}(y) )\}} \\
&\qquad\qquad \times[\Phi(r_{\tau-}+\{\delta_{(y)}\})-\Phi(r_{\tau-})] N_{\mathsf{af}} (d\tau,di,dy,d\alpha),
\end{split}
\label{Phifst}
\end{align}
where the three integral terms are associated with the three basic independent mechanisms and $\alpha\in[0,1]$ (for Monte Carlo acceptance-rejection). 
In this context, let us do a small remark on the Monte Carlo method. Parameter $\alpha$ corresponds to a decisional value for the type of the event that 
will occur by acceptance-rejection, which is usually defined as a random uniform realization in $[0,1]$. Let us explain the foundation of the Monte Carlo 
simulation scheme in detail. In particular, for all $t\geq0$, the explicit representation of the process $r_t$ is given a.s. by
\begin{align}
\begin{split}
 r_t=&r_0+\int_{0}^t \int_{\N^*}\int_{\mathbb{R}^d}\int_0^1 \mathds{1}_{\{i\leq s_{\tau-}\}}
\mathds{1}_{\{\alpha\leq ( k(x_{\tau-}^i,z))/(\gamma_1\tilde{k}(z)) \}}
  \delta_{(x_{\tau-}^i+z)} N_{v_r}(d\tau,di,dz,d\alpha)\\
&-\int_{0}^t \int_{\N^{\ast}}\int_0^1 \mathds{1}_{\{i\leq s_{\tau-}\}}\mathds{1}_{\{\alpha\leq(d_r/d_r)\}}
\delta_{(x_{\tau-}^i)}  N_{d_r}(d\tau,di,d\alpha)\\
&+\int_{0}^t\int_{\N^{\ast}}\int_{\bar{\mathcal{D}}}\int_0^1 \mathds{1}_{\{i\leq s_{\tau-}\}}
\mathds{1}_{\{\alpha\leq (\mathsf{aff}(x_{\tau-}^i,y)k^{\mathsf{af}}(y))/(A_\mathsf{f}\gamma_2 \tilde{k}^{\mathsf{af}}(y) )\}}\delta_{(y)} N_{\mathsf{af}}(d\tau,di,dy,d\alpha).
\end{split}
\label{netdef}
\end{align}
Thus, the network is expressed by (\ref{netdef}) in terms of the stochastic objects introduced above. Even the explicit expression looks somewhat complicated, the 
interpretation  is easy to discuss. The indicator functions that involve $\alpha$  are related to the rates and the 
dispersion kernels.  Indeed, the first term describes the recruitment by invitation event, the second term describes the departure from the 
network event and the last term describes the mechanism of recruitment by affinity. 

The $(\mathcal{F}_t)_{t\geq0}$-adapted stochastic process $(r_t)_{t\geq0}$ given by (\ref{netdef}) is Markovian with values in $\mathcal{S}$. For all bounded 
and measurable maps $\Phi:\mathcal{S}_F(\bar{\mathcal{D}})\mapsto \mathbb{R}$   and for all $r\in\mathcal{S}$, the  infinitesimal 
generator $G$ of the process $r_t$  is defined  by,
\begin{align}
\begin{split}
G\Phi(r)&=v_r\int_{\bar{\mathcal{D}}}r(dx) \int_{\mathbb{R}^d}\Big\{\Phi(r+\delta_{x+z})-\Phi(r)\Big\}k(x,z)dz 
+d_r \int_{\bar{\mathcal{D}}} \Big\{\Phi(r-\delta_{x})-\Phi(r)\Big\} r(dx)   \\
&\qquad+\int_{\bar{\mathcal{D}}}\Big\{ \int_{\bar{\mathcal{D}}}\big( \Phi(r+\delta_{y})-
\Phi(r)\big) \mathsf{aff}(x,y) {k}^{\mathsf{af}}(y)dy \Big\}r(dx).
\end{split}
\label{generat}
\end{align}
It is not hard to see that $r_t$ is a Markov process by classic arguments and to establish the expression (\ref{generat}) by differentiating the expectation 
of (\ref{Phifst}) at $t=0$. Note that if we define the extinction time as the stopping time:
\begin{align*}
 t_0 = \inf \{ t\geq 0, s_t=0\}
\end{align*}
with the convention $\inf \varnothing = \infty $ then before $t_0$ the infinitesimal generator is given by (\ref{generat}). After the extinction time $r_t$ is 
the null measure,  i.e. the network does not contain any individual and the infinitesimal generator is simply reduced to null measure. Furthermore, note that the 
representation  (\ref{netdef})  allows obtaining an easy simulation scheme for the numerical computation of the network which  will be explained further on below.

\subsection{Monte Carlo algorithm}

The distribution (law) of the network process is characterized by its infinitesimal generator (\ref{generat}). This characterization
is relatively abstract, so we subsequently propose now an exact Monte Carlo algorithm that simulates the network and provides an empirical
representation of its law. The method is exact as, up to the pseudo-random numbers generator approximation, it generates a network which has the same 
distribution as the considered Markov process $(r_t)_{t\geq 0}$. To give a computational representation of the network (\ref{netdef}), we shall detail   the 
associated algorithm. We suppose that the members in the network are independent and, as we have seen previously, the network has three 
possible events. To describe these events, we endow each member with three independent exponential clocks:
\begin{itemize}
\item a recruitment exponential affinity clock with rate $v_r$;
\item an invitation exponential clock with rate $A_{\mathsf{f}}$;
\item a departure exponential clock with rate $d_r$. 
\end{itemize}

So, endowing each member with 3 exponential clocks leads to a very high  total number of clocks ($3\times s_t$ clocks at time $t$) which is 
computationally intensive. Instead, a more efficient strategy consist in reducing  the number of exponential realizations by considering only  3
"fast" clocks: A global clock for  recruitment by affinity, a global clock for recruitment by invitation and a global clock for departure. Another more 
efficient strategy consist in considering one clock that control all the events thanks to the properties of the exponential law. We define one global clock 
that control all punctual mechanisms: 
\begin{align*}
\mathsf{H}_t= \mathsf{h}_t^{v_r} +\mathsf{h}_t^{\mathsf{af}} +\mathsf{h}_t^{d_r},\qquad 
\end{align*}
where
\begin{align*}
    \mathsf{h}_t^{v}= v_r s_t,  \quad  \textrm{ and } \quad
  \mathsf{h}_t^{\mathsf{af}}= A_{\mathsf{f}}s_t,  \quad  \textrm{ and } \quad \mathsf{h}_t^{d} = d_r s_t.
\end{align*}
To choose which event occurs, at each time step, we calculate the probabilities of each event. Let the 
time of the last event $T_{k-1}$ and the corresponding random network $r_{T_{k-1}}$, we simulate $t_k$ and $r_{t_k}$ as follows: We set 
$$T_k=T_{k-1}+S_k,$$ with $S_k\sim\textrm{Exp}(\mathsf{H}_{T_{k-1}})$ and 
\begin{align*}
 r_t=r_{T_{k-1}},  \qquad \textrm{ for  }\;  t \in[T_{k-1},T_k).
\end{align*} 
	
\noindent We calculate the probabilities of each event using the following objects:
\begin{align*} 
  \alpha_k^{v_r}
  &=  \frac{\mathsf{h}_{T_{k-1}}^{v_r} }  {   \mathsf{H}_{T_{k-1}} } \,,
  &
  \alpha_k^{\mathsf{af}}
  &=   \frac{\mathsf{h}_{T_{k-1}}^{\mathsf{af}} }{ \mathsf{H}_{T_{k-1}} } \,, \qquad\qquad \textrm{ and }
  &
  \alpha_k^{d_r}
  &=  \frac{ \mathsf{h}_{T_{k-1}}^{d_r}}  { \mathsf{H}_{T_{k-1}}}.
\end{align*}	
We draw the events as follows:
\begin{enumerate}
\item[(i)] With probability $\alpha_k^{v_r}$ an invitation event occurs. We draw a member $x_{T_{k-1}}^i$ where the index   
$i\sim U\{1,\dots, s_{T_{k-1}}\}$. We draw $z\in\mathbb{R}^d$ with the dispersal kernel  $K(x_{T_{k-1}}^i,dz)$. We add, with probability 
$\frac{k(x_{\tau-}^i,z)}{\gamma_1\tilde{k}(z)}$, a new 
individual to the virtual position $y=x_{T_{k-1}}^i+z$ and we set $$r_{T_k}=r_{T_{k-1}}+\delta_{x_{T_{k-1}}^i+z}.$$ 
\item[(ii)] With probability $\alpha_k^{d_r}$ a departure event occurs. We draw a member $x_{T_{k-1}}^i$ where the index
$i\sim U\{1,\dots, s_{T_{k-1}}\}$ and we set $$r_{T_k}=r_{T_{k-1}}-\delta_{x_{T_{k-1}}^i}.$$ 
\item[(iii)] With probability $\alpha_k^{\mathsf{af}}$ an affinity recruitment event occurs. We draw a state $y$ with the affinity kernel 
$K^{\mathsf{af}}(dy)$ and a member $x_{T_{k-1}}^i$ where the index  
$i\sim U\{1,\dots,s_{T_{k-1}}\}$. We reject the event of affinity recruitment with probability  
$$1-\frac{\mathsf{aff}(x_{T_{k-1}}^i,y)k^{\mathsf{af}}(y)}{{A}_{\mathsf{f}}\gamma_2\tilde{k}^{\mathsf{af}}(y)},$$ otherwise we add a new member to the virtual 
position $y$ and 
we set $$r_{T_k}=r_{T_{k-1}}+\delta_{y}.$$ 
\end{enumerate}     
\begin{algorithm}
\caption{Exact Monte Carlo algorithm of the random social network}
\begin{algorithmic} 
\STATE draw $r_0$ and set $T_0=0$
   \FOR{$k=1,2,...,N$}
    \STATE  \% Setting clocks 
    \STATE $\mathsf{h}_{T_{k-1}}^{v}= v_r s_{T_{k-1}},
  \mathsf{h}_{T_{k-1}}^{\mathsf{af}}= A_{\mathsf{f}}s_{T_{k-1}},  \mathsf{h}_{T_{k-1}}^{d_r} = d_r s_{T_{k-1}}$
    \STATE $\mathsf{H}_{T_{k-1}}= \mathsf{h}_{T_{k-1}}^{v_r} +\mathsf{h}_{T_{k-1}}^{\mathsf{af}} +\mathsf{h}_{T_{k-1}}^{d_r},\qquad$
    \STATE $S_k\sim \exp(H_{T_{k-1}})$ 
    \STATE $T_k=T_{k-1}+S_k$
    \item \noindent\% We calculate the probabilities:
    
    \STATE $\alpha_k^{v_r}=\mathsf{h}_{T_{k-1}}^{v_r}\big/\mathsf{H}_{T_{k-1}}; \; 
    \alpha_k^{\mathsf{af}}=\mathsf{h}_{T_{k-1}}^{\mathsf{af}} \big/ \mathsf{H}_{T_{k-1}} ; \;
    \alpha_k^{d_r}= \mathsf{h}_{T_{k-1}}^{d_r} \big/ \mathsf{H}_{T_{k-1}}$.   
    \STATE $u\sim\mathcal{U}[0,1]$
    \IF{$u\in[0,\alpha_k^{v_r}]$}
      \STATE  $i\sim\mathcal{U}\{1,...,s_{T_{k-1}}\}$
      \STATE $z\sim K(x_{T_{k-1}}^i,dz)$
       \STATE $u'\sim\mathcal{U}[0,1]$
      \IF {$u'\leq \frac{k(x_{\tau-}^i,z)}{\gamma_1\tilde{k}(z)}$}
      \STATE $y=x_{T_{k-1}}^i+z$
      \STATE $r_{T_k}=r_{T_{k-1}}+\delta_y$ \% Add the new invited member 
      \ENDIF
    \ELSE
      \IF {$u\in[\alpha_k^{v_r},\alpha_k^{v_r}+\mathsf{h}_{T_{k-1}}^{d_r}]$}
        \STATE $i\sim\mathcal{U}\{1,...,s_{T_{k-1}}\}$
        \STATE $r_{T_k}=r_{T_{k-1}}-\delta_{x_{T_{k-1}}^i}$ \% Departure of the member
      \ELSE 
        \STATE $y\sim k^{\mathsf{af}}(y)$
        \STATE $i\sim\mathcal{U}\{1,...,s_{T_{k-1}}\}$\\
        \STATE $u'\sim\mathcal{U}[0,1]$
          \IF {$u'\leq\frac{\mathsf{aff}(x_{T_{k-1}}^i,y)k^{\mathsf{af}}(y)}{A_\mathsf{f}\gamma_2\tilde{k}^{\mathsf{af}}(y)}$}
           \STATE $r_{T_k}=r_{T_{k-1}}+\delta_y$ 
          \ENDIF
       \ENDIF
     \ENDIF  

\ENDFOR
\end{algorithmic}
\label{algo}
\end{algorithm}

\begin{figure}
\center
\subfigure[Initial network state]
    {\includegraphics[width=5.5cm,height=5cm]
                 {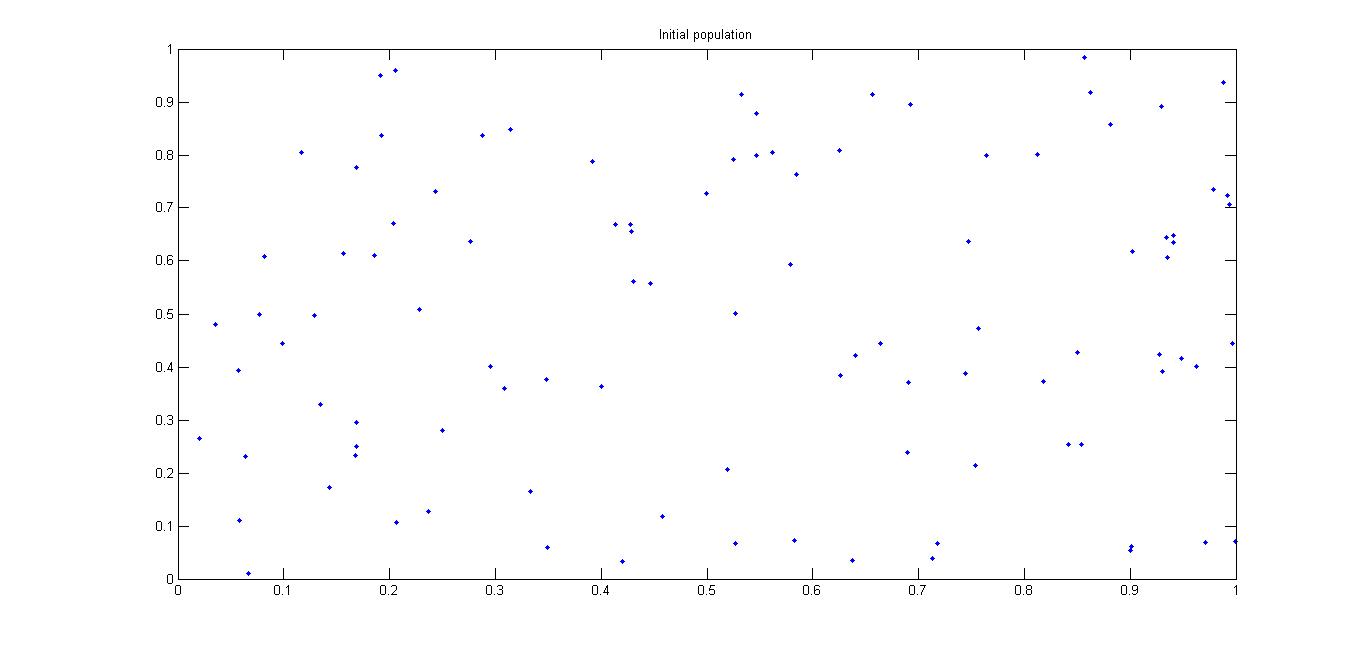}}
\subfigure[Final state of the network with a dispersion $\sigma=0.01$]
    {\includegraphics[width=5.5cm,height=5cm]
                 {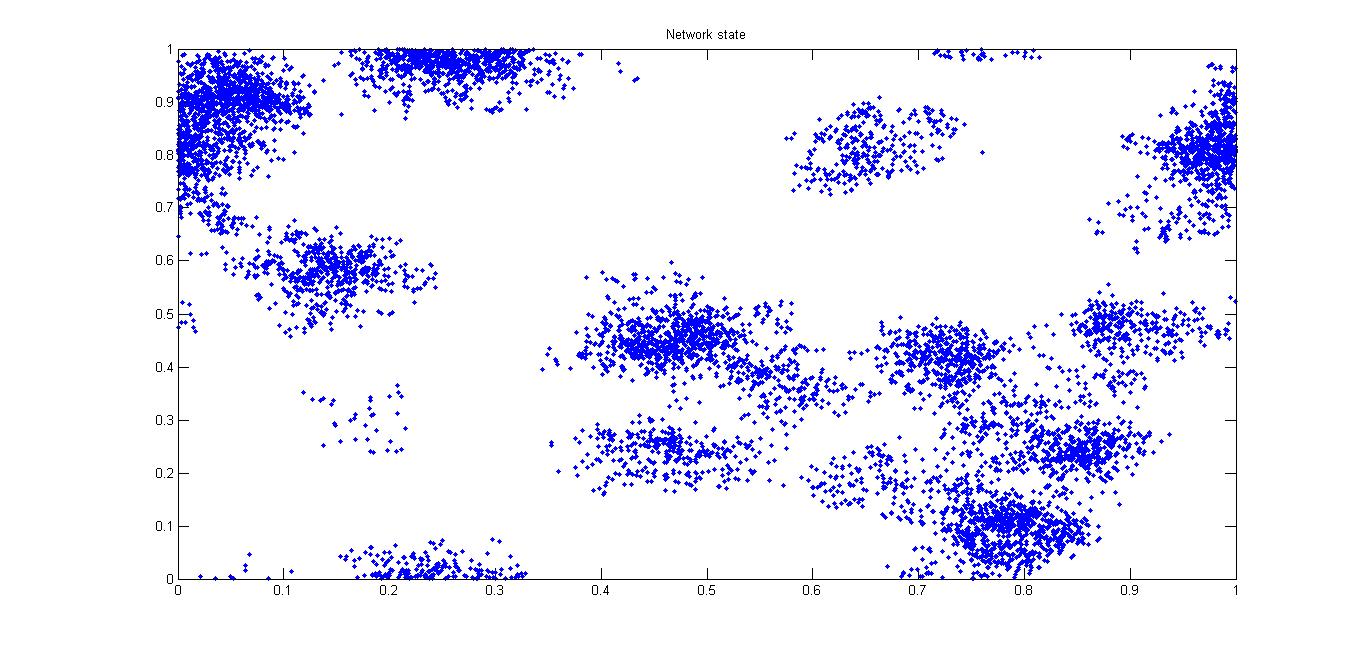}}
\subfigure[Final state of the network with a dispersion $\sigma=0.005$]
    {\includegraphics[width=5.5cm,height=5cm]
                 {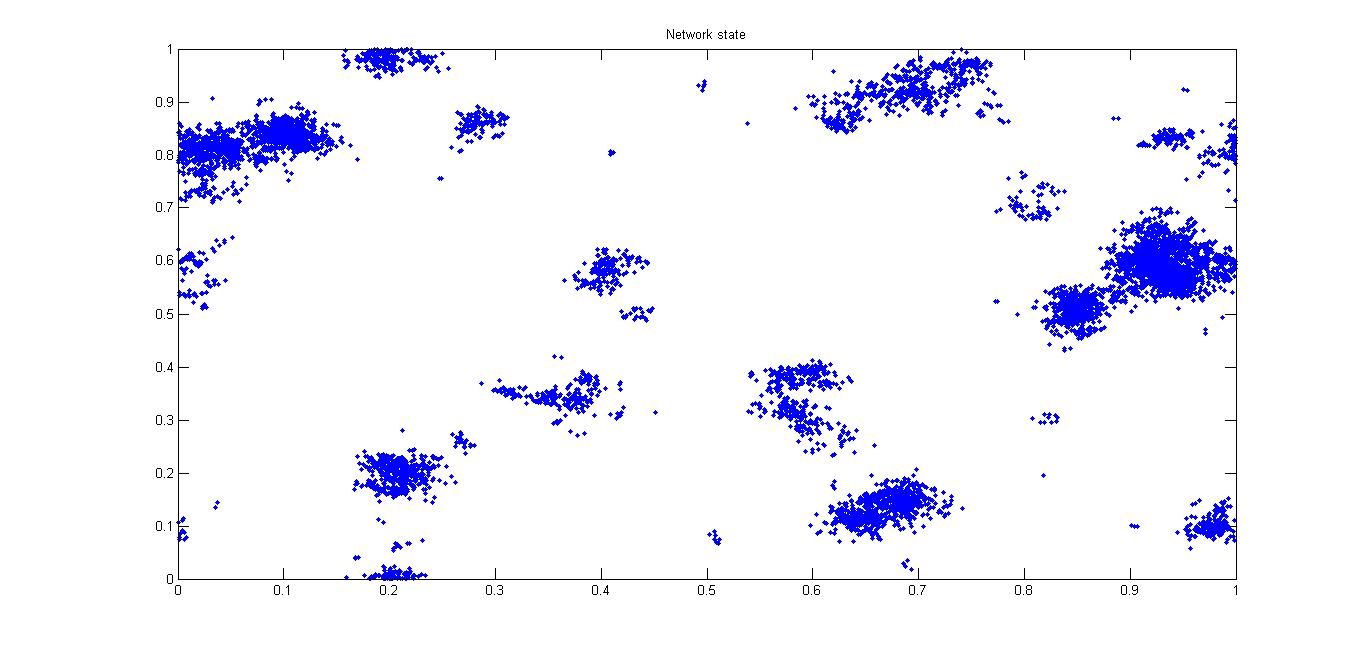}} 
\subfigure[Final state of the network with a dispersion $\sigma=0.001$]
    {\includegraphics[width=5.5cm,height=5cm]
                 {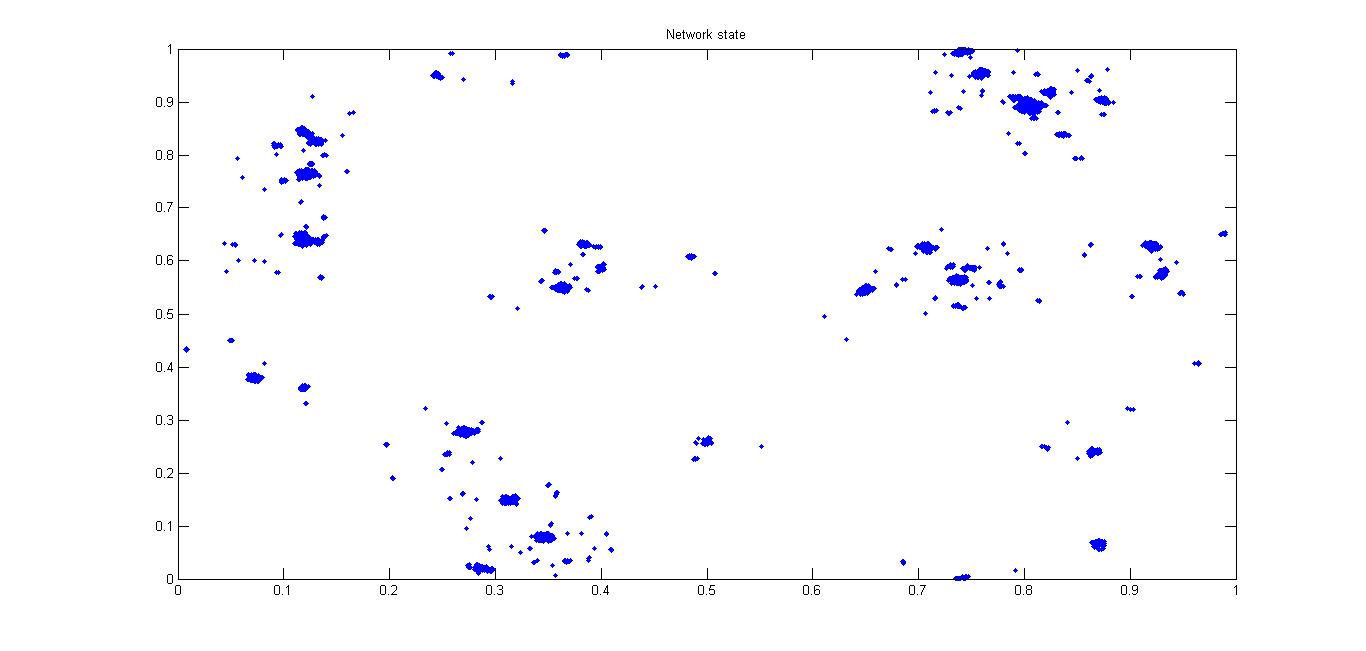}}             
\caption{{\small Figures shown the initial ($s_0=100$) and the final states of the random network $r_t$ (after $10^5$ iterations) using different 
dispersion levels  $\sigma=0.01$, $\sigma=0.005$ and $\sigma=0.001$.}}
\label{fig-sim}
\end{figure}

The exact Monte Carlo  scheme is detailed in Algorithm \ref{algo}. We are now in position to compute the network using a Monte Carlo strategy. Let us give a 
brief exposition of the use of the exact scheme to simulate some example of dynamics since an initial network (first members) through continuous time. Let 
the virtual space $\bar{\mathcal{D}}$ be the unit square $[0,1]^2$. Then, each member of the network is characterized by a point with two coordinates inside 
the square. We let the algorithm  run three times for the same number of iterations ($10^5$ updates) starting from the state $r_0\sim U_{[0,1]^2}$ with 
$s_0=100$. Following the description of the function $\mathsf{aff}(x,y)$, the friendship between two members inside the network is function  
to their Euclidean distance in the virtual space. Furthermore, the dynamic of the network depends heavily of the chosen parameters and dispersion kernels. We 
run the algorithm with $v_r=3$, $d_r=1.6$, $A_{\mathsf{f}}=2 $, $a_\mathsf{f}= 0.1$ and  we consider a normal kernel for the recruitments (by invitation and by affinity) 
with a certain dispersion parameter $\sigma$. We plot in Figure \ref{fig-sim} the final state of the three simulated dynamics 
(with $\sigma=0.01$, $\sigma=0.005$ and $\sigma=0.001$) 
from the initial state $r_0$. Our aim is to understand the behaviour of the network dynamics in function for instance of the dispersion. Thus,  we are 
interested in understanding the influence of spatial dispersion on the formation of patterns and other aspects of social network dynamics. The analysis of 
 Figures  \ref{fig-sim}(b-d)  shows that with high dispersion new members tend to occupy almost all the space and their dispersion is uniformly in the 
vicinity of already present members. With low dispersion new  members tend to occupy only the space located in the vicinity of their friend members. This 
favours the formation of clusters. Unsurprisingly with very small dispersion ($\sigma=0.001$), the network dynamic evolves again through a process of clustering. 
On the basis of the numerical study evidence, it appears that the impact of the dispersion level is strongest.

\section{Asymptotic adaptive threshold}
\label{sec:2}

As we have seen in the numerical tests, the model proposed in this framework contains several parameters that have significant effect on the dynamic and on the 
spatial patterns. Thus, a thorough analysis of the clusters (communities) seems crucial to ensure a well understanding of the different connections inside the network. 
Concerning the connectivity, it 
is well known that  a fundamental question about any network is whether or not it is connected. We study from now on the clustering problem of the set 
of random points $\{x_{t}^1,\ldots,x_{t}^{s_t}\}$ in $[0,1]^2$. To avoid technicalities arising from irregularities around the borders of $[0,1]^2$, we 
consider the unit square as a torus.  As in the random geometric graph of Gilbert, given the radius of affinity 
$a_\mathsf{f}>0$, we may consider our network $r_t$ in which members $i$ and $j$ are connected if and only if the distance of $x_{t}^i$ and 
$x_{t}^j$ does not exceed the threshold $a_\mathsf{f}$. We shall consider an extension to this threshold by introducing a soft assumption that all members 
are active independently with a  probability $p$ (be active $\sim$ Bernoulli$(p)$, $0<p\leq 1$). Such soft assumption is motivated by ad-hoc social networks. From a practical viewpoint, a 
member may be inactive in the network for different reasons and will not take part in a community membership. Hence, taking into account the activity of the 
members, our social network is said to be connected if each inactive member is adjacent to at least one active member together with all active members form 
a connected network. 

In the exact Algorithm \ref{algo}, the generative model for the data $\mathbf{r}_{T_k}=\{x_{T_k}^1,\ldots,x_{T_k}^{s_{T_k}}\}$, for $1\leq k\leq N$, is a inhomogeneous 
 Poisson point process with a spatial dependence intensity $\lambda_k^\star$ (such that  $\lambda_k^\star(x)=v_rk(y,x) + d_r + \sum_{y\in r_{T_k}}\mathsf{aff}(x,y)k^{\mathsf{af}}(x)$, 
 for $x\in[0,1]^2$). We note that the intensity $\lambda_k^\star$ is a locally integrable positive function. Now consider the sequences 
$(\mathbf{r}_{T_k})_{k\geq 1}$ and $(s_{T_k})_{k\geq 1}$ of random variables respectively in $[0,1]^2$ and $\mathbb{N}$, both adapted to a filtration 
$(\mathcal{F}_k)_{k\geq 1}$. We study the problem of the clustering (known as community detection in network literature) of the random points  
$\mathbf{r}_{T_k}$ at each time $T_k$ based on the observation of $(\mathbf{r}_{T_1},\ldots,\mathbf{r}_{T_N})$ and $(s_{T_1},\ldots,s_{T_N})$, where 
$N\geq 1$ is a finite $(\mathcal{F}_k)$-stopping time.

\begin{rem}
 The results given from now on are stated in a setting where one observes $(\mathbf{r}_{T_k},s_{T_k})_{k= 1}^N$ with $N$ a stopping time. It is 
 worth pointing out that this contains the usual case $N\equiv 1$ and $s_{T_1}\equiv n$, where $n$ is a fixed sample size. This strategy includes situations 
 where the statistician decides to stop the recording process of network data according to some design of experiment rule. 
\end{rem}

The analysis in this section  is conducted under the following assumption.

\begin{hyp}[$H_1$]
 There is a $(\mathcal{F}_k)$-adapted sequence of functions $\big(\mathfrak{f}_k({s_{T_k}})\big)_{k\geq 1}$ of positive random variables. This sequence has 
a limit assumed to be known  such that for $\mathfrak{f}_k:\mathbb{N}\mapsto\mathbb{R}_+$
 \begin{enumerate} 
\item  For all $k\geq 1$,
\begin{align}
 \psi_k({s_{T_k}}):=\log\Big(\frac{s_{T_k} \pi}{4(s_{T_k}-1)!}\Big)+\mathfrak{f}_k({s_{T_k}})+ (s_{T_k}-1)\log(\alpha_k^\star(s_{T_k}))-
 \alpha_k^\star(s_{T_k})  >0 ,
\end{align}
where the sequence of functions $(\alpha_k^\star)_{k\geq 1}$ will be stated later and depends on, among other, the constants 
$v_r, d_r, A_{\mathsf{f}},\sigma,a_\mathsf{f}$.
\item For all $k\geq 1$,
  \begin{align}
   \lim_{s_{T_k}\to \infty} \frac{ \psi_k({s_{T_k}})}{s_{T_k}} = 0, \qquad \textrm{ and } \qquad  \lim_{s_{T_k}\to \infty}  \psi_k({s_{T_k}}) = +\infty.
  \end{align}
 \end{enumerate}
\label{hyp2}
\end{hyp}

 The first condition of  assumption \ref{hyp2} allows us to overcome the threshold $a_\mathsf{f}$ problem in the clustering  whereas the second condition 
  means  that $\psi_k $ (a function that depends on the current size of the network) does not grow to $\infty$ more faster than  $s_{T_k}$. 
 
 We now come to our main results. Define $\rho(x)=|B(x,a_\mathsf{f}) \cap [0,1]^2|$ to be the number of points (active and inactive) in affinity with $x$, 
 where $B(x,a_\mathsf{f})=\{y\in[0,1]^2: \|x-y\|\leq a_\mathsf{f}  \}$. Furthermore, we consider that a member $x$ is isolated in the network when 
 $\rho(x)$ contains only the number of inactive members adjacent to $x$ (which means that $x$ is considered isolated if there is no active members in the 
 ball $B(x,a_\mathsf{f})$ except at least $x$). Let $B(x)$ denotes the intersection of the ball $B(x,a_\mathsf{f})$ with $[0,1]^2$. Next, with 
 $\textrm{vol}(x)$ the volume of $B(x)$ (in other words the Lebesgue measure of the measurable set $B(x)$), for all $x\in [0,1]^2$ we have 
 $\textrm{vol}(x)\leq \pi a_\mathsf{f}^2$.  For a purely notational reason, let 
 $E_{ik}$ be the event that $x_{T_k}^i$ is isolated and $F_{ik}$ be the event that $x_{T_k}^i$ is an active isolated member for 
 $i\in\{1,\ldots,s_{T_k}\}$ and $k\in\{1,\ldots,N\}$.
 
 The following lemma gathers some standard deviation estimate (a concentration bound) on Poisson random variables; see \cite{Bouch+al13} for more 
 details.
 
 \begin{lem}
  Let $\mathcal{P}(\lambda_0)$ be a Poisson random distribution with mean $\lambda_0$. Then, there exists a $\delta_0>0$ such that for all 
  $\delta\in[0,\delta_0]$ we have
\begin{align}
 \mathbb{P}(|\mathcal{P}(\lambda_0)-\lambda_0|\geq \delta \lambda_0) \leq  2 \textrm{e}^{-\lambda_0\delta^2/3}.
\end{align}
\label{concent}
 \end{lem}
\proof For $\delta\geq \lambda_0$ set $\mu=\delta/\lambda_0$ in the inequality $\mu^\delta \mathbb{P}(\mathcal{P}(\lambda_0) \geq \delta)
\leq \mathbb{E}[\mu^X]=\textrm{e}^{\lambda_0(\mu-1)}$ where $X\sim\mathcal{P}(\lambda_0)$, to obtain
\begin{align*}
 \mathbb{P}(\mathcal{P}(\lambda_0) \geq \delta)\leq \textrm{e}^{-\delta\log(\mu)+\lambda_0\mu-\lambda_0}=\textrm{e}^{-\lambda_0\Big(\frac{\delta}{\lambda_0}
 \log(\frac{\delta}{\lambda_0})-\frac{\delta}{\lambda_0}+1\Big)}.
\end{align*}
Hence,  we get using a second order Taylor approximation that
\begin{align*}
 \mathbb{P}\Big(\mathcal{P}(\lambda_0)-\lambda_0 \geq \delta \lambda_0\Big)\leq  \textrm{e}^{-\lambda_0\Big((1+\delta)\log(1+\delta)-\delta \Big)} \leq  
 \textrm{e}^{-\lambda_0\delta^2/3}.
\end{align*}
Now, if $0<\delta\leq \lambda_0$, one can set  $\mu=\delta/\lambda_0$ in the inequality $\mathbb{P}(\mathcal{P}(\lambda_0) 
\leq \delta)\leq \mu^{-\delta} \mathbb{E}[\mu^X]$ to establish
\begin{align*}
  \mathbb{P}(\mathcal{P}(\lambda_0) \leq \delta)\leq \textrm{e}^{-\lambda_0\Big(\frac{\delta}{\lambda_0}
 \log(\frac{\delta}{\lambda_0})-\frac{\delta}{\lambda_0}+1\Big)}.
\end{align*}
A similar bound for $\mathbb{P}\Big(\mathcal{P}(\lambda_0)-\lambda_0 \leq \delta \lambda_0\Big)$ follows similarly. \carre
 
 For the random geometric social network, the number of isolated members (denoted from now on $N_0^k$) enjoys a Poisson approximation when the size of the 
 network tends to $\infty$. So, for all $k\geq 1$ and for $m\in\mathbb{N}$, we have 
 \begin{align}
  \mathbb{P}\Big(N_0^k=m\Big) {\stackrel{{{\scriptscriptstyle s_{T_k}\to+\infty}}}{\longrightarrow}} \textrm{e}^{-\mathbb{E}[N_0^k]} \frac{\big(\mathbb{E}[N_0^k]\big)^m}{m!}.
 \label{approxpois}
 \end{align}
In the present section we prove result of this kind for the class of random network model described in Section \ref{sec:1} when we connect each pair of 
members with an indicator function $ \mathsf{aff}(\cdot,\cdot)$ of the distance between them. We show that the approximation $(\ref{approxpois})$ holds for 
the network $r_{T_k}$ (for all $k\geq 1$) for large network size, uniformly over affinity functions that are zero beyond a given distance. The proof  relies 
heavily on some levels of discretization of the unit square into smaller subsquares (open boxes). Assume that $\bar{\mathcal{D}}=[0,1]^2$ is  compact in $\mathbb{R}^2$ 
and consider the partitioning of $[0,1]^2$ into a family $\{\boldsymbol{C}_\ell\}_{\ell=1,\ldots, L_k}$  of disjoint boxes of $\mathbb{R}^2$ with side 
$a_\mathsf{f}'/\sqrt{2}$ that we need to cover the state space  $\bar{\mathcal{D}}$ where $a_\mathsf{f}'\leq a_\mathsf{f} $. In the course of the proofs, for all $k\geq1$, we condition on the 
locations of the points $ x_{T_k}^1,\ldots,x_{T_k}^{s_{T_k}}$ and assume that they are sufficiently regularly distributed. The probability that this holds is proved 
in the Lemma \ref{regul} that relies on concentration bound of large deviations for Poisson random variables presented in Lemma \ref{concent}. For a 
box $\boldsymbol{C}_\ell$, we have
\begin{align*}
 \mathbb{E}\Big[|\boldsymbol{C}_\ell \cap [0,1]^2|\Big] = \frac{s_{T_k}}{L_k}= \frac{s_{T_k}\times (a_\mathsf{f}')^2}{2}.
\end{align*}
Here, we need a definition for the regularity of boxes $\{\boldsymbol{C}_\ell\}_{\ell=1,\ldots, L_k}$.
\begin{defi}
 Fix $\nu\in(0,1)$. A box $\boldsymbol{C}_\ell$ is called $\nu$-\emph{regular} if one has
 \begin{align*}
  \frac{(1-\nu)s_{T_k} (a_\mathsf{f}')^2}{2}  \leq  |\boldsymbol{C}_\ell \cap [0,1]^2|  \leq \frac{(1+\nu)s_{T_k} (a_\mathsf{f}')^2}{2}.
 \end{align*}
\end{defi}

 In the following simple lemma we estimate the probability that all boxes $\boldsymbol{C}_\ell$ are $\nu$-regular for $s_{T_k}$ large enough.
 \begin{lem}
  There exists  $\delta_k'>0$ and $\gamma_k>0$ such that for all $\delta_k\in[0,\delta_k']$ and if $a_\mathsf{f}'(s_{T_k})=\gamma_k\sqrt{\log(s_{T_k})/s_{T_k}}$ for 
  all $k\geq1$, then for all $s_{T_k}$ large enough,
  \begin{align*}
   \inf_{\ell=1,\ldots,L_k} \mathbb{P}(\boldsymbol{C}_\ell \textrm{ is } \nu-\emph{regular}) \geq 1-2 L_k s_{T_k}^{-\gamma_k^2 \delta_k^2/6}.
  \end{align*}
In particular, if $\gamma_k^2>6/\delta_k^2$ then
\begin{align*}
 \mathbb{P}(\textrm{every box }  \boldsymbol{C}_\ell \textrm{ is } \nu-\emph{regular}) 
 {\stackrel{{{\scriptscriptstyle s_{T_k}\to+\infty}}}{\longrightarrow}} 1.
\end{align*}
  \label{regul}
 \end{lem}

 \proof 
 For any box $\boldsymbol{C}_\ell$, the number of points $ |\boldsymbol{C}_\ell \cap [0,1]^2|$ is distributed like a Poisson random variable with mean 
 $s_{T_k}\times (a_\mathsf{f}')^2/2$. By Lemma \ref{concent}, we have for $\delta_k'>0$ and $\delta_k\in[0,\delta_k']$,
 \begin{align}
  \mathbb{P}\Big(\Big|\mathcal{P}(\frac{s_{T_k} (a_\mathsf{f}')^2}{2})-\frac{s_{T_k} (a_\mathsf{f}')^2}{2}\Big|\geq \delta_k \frac{s_{T_k} (a_\mathsf{f}')^2}{2}\Big) \leq
  2 \textrm{e}^{-s_{T_k} (a_\mathsf{f}')^2\delta_k^2/2.3}.
 \end{align}
Now, for every box and for all $s_{T_k}$ large enough, 
\begin{align*}
 \mathbb{P}( \boldsymbol{C}_\ell \textrm{ is not } \nu-\textrm{regular}) & \leq  2 \textrm{e}^{-s_{T_k} (a_\mathsf{f}')^2\delta_k^2/6}\\
 & \leq 2  s_{T_k}^{-\gamma_k^2\delta_k^2/6},
\end{align*}
since $a_\mathsf{f}'(s_{T_k})=\gamma_k\sqrt{\log(s_{T_k})/s_{T_k}}$ for any $k\geq 1$. In addition, if there exists one box that is not $\nu$-regular, 
then one of the $L_k$ boxes has a number of points that is out of range, so that as $s_{T_k}\to\infty$,
\begin{align*}
 \mathbb{P}(\exists\boldsymbol{C}_\ell:  \boldsymbol{C}_\ell \textrm{ is not } \nu-\textrm{regular}) &\leq  2 \times L_k \times s_{T_k}^{-\gamma_k^2\delta_k^2/6}
 = 2 \times \frac{2s_{T_k}}{\gamma_k^2\log(s_{T_k})} \times s_{T_k}^{-\gamma_k^2\delta_k^2/6}\\
 &\leq  s_{T_k}^{1-\gamma_k^2\delta_k^2/6+o(1)},
\end{align*}
which tends to zero provided that $\gamma_k^2>6/\delta_k^2$. \carre 
 
 For simplicity and from now on we make the following technical assumption.      
 \begin{hyp}[$H_2$]
 \begin{enumerate}
  \item  Assume that $a_\mathsf{f}(s_{T_k})\geq a_\mathsf{f}'(s_{T_k})$ and $a_\mathsf{f}'(s_{T_k})\sim \gamma_k\sqrt{\log(s_{T_k})/s_{T_k}}$ for all $k\geq 1$ 
  with $\gamma_k>\sqrt{6}/\delta_k$ and $\delta_k\in[0,\delta_k']$ and $\delta_k'>0$.  Assume also that there exists  a
  $(\mathcal{F}_k)$-adapted sequence of functions $\big(g_{k}(x,s_{T_k}) \big)_{k\geq 1}$, with $g_{k}:[0,1]^2\times \mathbb{N}\mapsto\mathbb{R}_+$, such that, for all 
  $x,x_{T_k}^1,\ldots,x_{T_k}^{s_{T_k}} \in[0,1]^2$ and all $k\geq 1$, we have 
\begin{align*}
 \sum_{i=1}^{s_{T_k}} \mathds{1}_{\{|x-x_{T_k}^i|\leq a_\mathsf{f}\}} = g_{s_{T_k}}(x) \sum_{i=1}^{s_{T_k}} \sum_{\ell=1}^{L_k} 
 \mathds{1}_{\{(\boldsymbol{C}_\ell\cap B(x))\neq \varnothing\} } \mathds{1}_{[0,1]^2}(x)
 \mathds{1}_{\boldsymbol{C}_\ell}(x_{T_k}^i), \quad \textrm{ with } g_{k}(x,s_{T_k}) {\stackrel{{{\scriptscriptstyle s_{T_k}\to+\infty}}}{\longrightarrow}} 1,
  \label{hypcub}
\end{align*}
where $\{\boldsymbol{C}_\ell\}_{\ell=1,\ldots, L_k}$ is a family of disjoint boxes with side $a_\mathsf{f}'(s_{T_k})/\sqrt{2}$ that covers $[0,1]^2$.
\item  The dispersion kernel $k^{\mathsf{af}}(\cdot):\mathbb{R}^2\mapsto\mathbb{R_+}$ is  an integrable function that does not change sign in $[0,1]^2$.
\end{enumerate}
 \end{hyp}

 The assumption  $H_2$ enables us to control the accuracy of the approximation of the number of points in the neighborhood of each $x\in[0,1]^2$. Note 
 that $L_k$ is clearly finite (for all $k\geq 1$) together with for each $x,y\in \boldsymbol{C}_\ell$, $|x-y|\leq a_\mathsf{f}'$.
 Note also that the function $g_{k}$ assesses the proportion of points inside the ball $B(x)$ from the number of points inside the boxes that intersect 
 with $B(x)$.  Such an assumption is realistic because, when $s_{T_k}$ is large enough, the side of each  box $a_\mathsf{f}'(s_{T_k})/\sqrt{2}$ is small 
 enough which guarantee a fine mesh and hence $g_{k}\sim 1$ for all $x\in[0,1]^2$. More rigorously, we use the fact that all open of 
 $\mathbb{R}^d$ is a countable reunion of open pavers.  Furthermore, the assumption is realistic  since 
 asymptotically  all the boxes contains points (proved in Lemma \ref{regul}) and is needed in practice to discard  regions with strong variations. Concerning 
 the point $(ii)$ in assumption $H_2$ it's mainly needed for computational issues. 
 
 In order to use $H_2$ for estimating the probability of isolated members, we need to make sure that the global affinity rate stays under control.
 
 \begin{lem}
 Admit the hypothesis  $H_2$. For any $x\in[0,1]^2$ and any bounded region $B(x)$, let consider the (finite) integral of $\lambda_k^\star $ 
 over region $B(x)$  
  \begin{align}
  \Lambda_k^\star (x)= \int_{B(x)}\lambda_k^\star(z)dz =\int_{B(x)} \Big(v_r k(y,z) + d_r + \sum_{y\in r_{T_k}}\mathsf{aff}(z,y)k^{\mathsf{af}}(z)\Big)dz,
  \end{align}
  with $v_r,d_r,A_\mathsf{f}>0,\sigma,a_\mathsf{f}>0$. Then, there exists two constants $\tilde{c}_0$ and $\tilde{c}_1$ in  $B(x)$ such that for $s_{T_k}$ large enough
  \begin{align}
   \Lambda_k^\star (x) =  \Big(v_r k(y,\tilde{c}_0) + d_r + A_\mathsf{f} L_k'\frac{\gamma_k^2\log(s_{T_k})k^{\mathsf{af}}(\tilde{c}_1)}{2} \Big) \emph{vol}(x),
  \end{align}
where the integer $L_k'\in \mathcal{S}_0:=\Big\{\lceil\frac{ \pi a_\mathsf{f}^2}{2 (a_\mathsf{f}')^2}\rceil,\ldots,
\lceil\frac{2\pi a_\mathsf{f}^2}{(a_\mathsf{f}')^2}\rceil\Big\}$ and $\gamma_k>\sqrt{6}/\delta_k$ with $\delta_k\in[0,\delta_k']$ and $\delta_k'>0$.
 \end{lem}
\proof 
We use the assumption $H_2$ to compute the integral of $\lambda_k^\star $ over region $B(x)$. First, for any $z\in B(x)$ we easily remark that the number of 
boxes that intersect with the bounded region $B(z)$ is  between $\lceil\frac{ \pi a_\mathsf{f}^2/4}{(a_\mathsf{f}'/\sqrt{2})^2}\rceil$ (when $z$ is located at 
the corners of $[0,1]^2$) and $\lceil\frac{\pi a_\mathsf{f}^2}{(a_\mathsf{f}'/\sqrt{2})^2}\rceil$ (when $B(z)$ is fully contained in $[0,1]^2$). For all 
$k\geq 1$, let $L_k'$ 
denote the number of boxes whose intersects with $B(z)$. If $s_{T_k}$ is sufficiently large, every square that intersect with $B(z)$ is fully contained in 
$B(z)$ and the large number of boxes allows us to take advantage of the approximation $g_{k}(z,s_{T_k})\sim 1$ for all $z\in[0,1]^2$.  Now, since every square 
is $\nu$-regular  for every $z\in[0,1]^2$, it follows that
\begin{align*}
\int_{B(x)} A_\mathsf{f}  \sum_{i=1}^{s_{T_k}} \mathds{1}_{\{|z-x_{T_k}^i|\leq a_\mathsf{f}\}} dz &= A_\mathsf{f}    \int_{B(x)} \sum_{i=1}^{s_{T_k}}
\sum_{\ell=1}^{L_k} \mathds{1}_{\{(\boldsymbol{C}_\ell\cap B(z))\neq \varnothing\} }
\mathds{1}_{\boldsymbol{C}_\ell}(x_{T_k}^i)\mathds{1}_{B(x)}(z) dz\\
&= A_\mathsf{f}    \int_{B(x)} \sum_{i=1}^{s_{T_k}}
\sum_{\ell=1}^{L_k'} \mathds{1}_{\{(\boldsymbol{C}_\ell\cap B(z))\neq \varnothing\} }
\mathds{1}_{\boldsymbol{C}_\ell}(x_{T_k}^i)\mathds{1}_{B(x)}(z) dz\\
 &= \Big(A_\mathsf{f} L_k'  \frac{s_{T_k}\times (a_\mathsf{f}')^2}{2} \Big) \int_{B(x)} \mathds{1}_{B(x)}(z)dz \\
 &= \Big(A_\mathsf{f} L_k'  \frac{s_{T_k}\times \gamma_k^2\log(s_{T_k})/s_{T_k}}{2} \Big) \textrm{vol}(x),
\end{align*}
and the claim follows easily by application of the mean-value theorem. \carre 
 
 We are now in position to assess the probability of an isolated member inside the network. 
 
 \begin{prop}
 Under the assumption $H_2$ and  for any  $x_{T_k}^i \in [0,1]^2$ with $v_r,d_r,A_\mathsf{f},\sigma,a_\mathsf{f}>0$ and $0<p\leq 1$, we have for $s_{T_k}$  
 sufficiently large and $i\in\{1,\ldots,s_{T_k}\}$ with $k\in\{1,\ldots,N\}$
  \begin{align*}
   \mathbb{P}(E_{ik}) = & \frac{\big(v_r k(y,\tilde{c}_0) + d_r + A_\mathsf{f} L_k'\frac{\gamma_k^2\log(s_{T_k})k^{\mathsf{af}}(\tilde{c}_1)}{2} 
   \big)^{s_{T_k}-1}}{(s_{T_k}-1)!} \emph{e}^{-\big(v_rk(y,\tilde{c}_0) + d_r + A_\mathsf{f} 
 L_k'\frac{\gamma_k^2\log(s_{T_k})k^{\mathsf{af}}(\tilde{c}_1)}{2} \big) } \\
 &\qquad \qquad \times \int_{[0,1]^2} \Big(1-p\emph{vol}(x_{T_k}^i)\Big)^{s_{T_k}-1}dx_{T_k}^i,
  \end{align*}
  where  $L_k'\in \mathcal{S}_0$, $k(y,\tilde{c}_0)> 0$, $k^{\mathsf{af}}(\tilde{c}_1)> 0$ and $\gamma_k>\sqrt{6}/\delta_k$ with $\delta_k\in[0,\delta_k']$ and $\delta_k'>0$.
  \label{propisol}
 \end{prop}
\proof At first sight and for simplicity, let $\Omega^{ik}$ be the event that, for any $x_{T_k}^i \in [0,1]^2$ and for any $j\neq i$, 
$x_{T_k}^j\notin B(x_{T_k}^i)$ or $x_{T_k}^j\in B(x_{T_k}^i)$ but inactive.  By  direct calculation we have
\begin{align*}
 \mathbb{P}(E_{ik}) &=  \int_{[0,1]^2}  \mathbb{P}\Big(E_{ik}\Big|x_{T_k}^i\in [0,1]^2\Big) dx_{T_k}^i =  \int_{[0,1]^2}  \mathbb{P}(\Omega^{ik}) dx_{T_k}^i \\
 &=\int_{[0,1]^2}  \sum_{\ell=0}^{s_{T_k}-1} \Big((1-p)^\ell \frac{(\int_{B(x_{T_k}^i)} \lambda_k^\star(x) dx)^\ell}{\ell!}\textrm{e}^{-\int_{B(x_{T_k}^i)} 
 \lambda_k^\star(x) dx}
 \Big) \\
 &\qquad \qquad  \times \Big( \frac{(\int_{[0,1]^2\setminus B(x_{T_k}^i)} \lambda_k^\star(x) dx)^{s_{T_k}-1-\ell}}{(s_{T_k}-1-\ell)!}
 \textrm{e}^{-\int_{[0,1]^2\setminus B(x_{T_k}^i)} \lambda_k^\star(x) dx}\Big) dx_{T_k}^i \\
 &= \int_{[0,1]^2}  \sum_{\ell=0}^{s_{T_k}-1}\Big((1-p)^\ell \frac{\Big(\Big(v_r k(y,\tilde{c}_0) + d_r + A_\mathsf{f} L_k'\frac{\gamma_k^2\log(s_{T_k})k^{\mathsf{af}}(\tilde{c}_1)}{2} \Big) 
 \textrm{vol}(x_{T_k}^i)\Big)^\ell}{\ell!} \\
& \qquad\times  \textrm{e}^{-\big(v_r k(y,\tilde{c}_0) + d_r + A_\mathsf{f} L_k'\frac{\gamma_k^2\log(s_{T_k})k^{\mathsf{af}}(\tilde{c}_1)}{2} \big) \textrm{vol}(x_{T_k}^i)}
 \Big) \\
 & \qquad  \times 
 \Big( \frac{\Big(\Big(v_rk(y,\tilde{c}_0) + d_r + A_\mathsf{f} L_k'\frac{\gamma_k^2\log(s_{T_k})k^{\mathsf{af}}(\tilde{c}_1)}{2} \Big) 
 \big(1-\textrm{vol}(x_{T_k}^i)\big)\Big)^{(s_{T_k}-1-\ell)}}{(s_{T_k}-1-\ell)!} \\
&\qquad \times  \textrm{e}^{-\big(v_rk(y,\tilde{c}_0) + d_r + A_\mathsf{f} L_k'\frac{\gamma_k^2
 \log(s_{T_k})k^{\mathsf{af}}(\tilde{c}_1)}{2} \big) 
 \big(1-\textrm{vol}(x_{T_k}^i)\big)}\Big) dx_{T_k}^i\\
 &= \frac{\textrm{e}^{-\big(v_r k(y,\tilde{c}_0) + d_r + A_\mathsf{f} L_k'\frac{\gamma_k^2\log(s_{T_k})k^{\mathsf{af}}(\tilde{c}_1)}{2} \big) }}{(s_{T_k}-1)!} 
 \int_{[0,1]^2}  \sum_{\ell=0}^{s_{T_k}-1}  \frac{(s_{T_k}-1)!}{\ell!(s_{T_k}-1-\ell)!} \\
&\qquad \times \Big((1-p)\big(v_rk(y,\tilde{c}_0) + d_r + A_\mathsf{f} L_k'\frac{\gamma_k^2\log(s_{T_k})k^{\mathsf{af}}(\tilde{c}_1)}{2} \big) 
 \textrm{vol}(x_{T_k}^i)\Big)^\ell \\
 &\qquad \times \Big(\big(v_rk(y,\tilde{c}_0) + d_r + A_\mathsf{f} L_k'\frac{\gamma_k^2\log(s_{T_k})k^{\mathsf{af}}(\tilde{c}_1)}{2} \big) 
 \big(1-\textrm{vol}(x_{T_k}^i)\big)\Big)^{s_{T_k}-1-\ell} dx_{T_k}^i
 \end{align*}
 \begin{align*}
 = \frac{\big(v_rk(y,\tilde{c}_0) + d_r + A_\mathsf{f} L_k'\frac{\gamma_k^2\log(s_{T_k})k^{\mathsf{af}}(\tilde{c}_1)}{2} \big)^{s_{T_k}-1}}{(s_{T_k}-1)!} \textrm{e}^{-\big(v_rk(y,\tilde{c}_0) + d_r + A_\mathsf{f} 
 L_k'\frac{\gamma_k^2\log(s_{T_k})k^{\mathsf{af}}(\tilde{c}_1)}{2} \big) } \\
  \qquad \times \int_{[0,1]^2} \Big(1-p\textrm{vol}(x_{T_k}^i)\Big)^{s_{T_k}-1}dx_{T_k}^i,
\end{align*}
where the last line is obtained from the binomial theorem.  This completes the proof. \carre

It is  interesting to remark that the result cited in Proposition \ref{propisol} suggests that the probability of being a member isolated in the network is inversely 
proportional to the size of the network and at the same time to the volume of the ball portion (around the member) that intersect with $[0,1]^2$. Thus, this 
result seems intuitive. Now, we shall assess the probability that more than one member (say $\kappa\geq 2$ members) are isolated inside the network. For any 
$(x^1,\ldots,x^\kappa)\in[0,1]^{2\kappa} $ and to shorten notation, we use from now on  $B(x^1,\ldots,x^\kappa)=B(x^1)\cup \cdots \cup B(x^\kappa)$ and 
$\textrm{vol}(x^1,\ldots,x^\kappa)$ for the volume of $B(x^1,\ldots,x^\kappa)$.

 \begin{prop}
 Under the assumption $H_2$ and  for any $ \kappa\geq 2 $ and $(x_{T_k}^1,\ldots,x_{T_k}^\kappa)\in[0,1]^{2\kappa}$ with $v_r,d_r,A_\mathsf{f},\sigma,a_\mathsf{f}>0$ 
 and $0<p\leq 1$, we have for $s_{T_k}$  sufficiently large and $k\in\{1,\ldots,N\}$
  \begin{align*}
   \mathbb{P}(E_{1k}\cap\cdots\cap E_{\kappa k}) &\leq  \frac{\big(v_rk(y,\tilde{c}_0) + d_r + A_\mathsf{f} L_k'\frac{\gamma_k^2\log(s_{T_k})k^{\mathsf{af}}(\tilde{c}_1)}{2}
   \big)^{s_{T_k}-\kappa}}{(s_{T_k}-\kappa)!} \emph{e}^{-\big(v_r k(y,\tilde{c}_0) + d_r + A_\mathsf{f}  L_k'\frac{\gamma_k^2\log(s_{T_k})k^{\mathsf{af}}(\tilde{c}_1)}{2} \big) } \\
 &\qquad \qquad  \times \int_{[0,1]^{2\kappa}} \Big(1-p\emph{vol}(x_{T_k}^1,\ldots,x_{T_k}^\kappa)\Big)^{s_{T_k}-\kappa}dx_{T_k}^1\cdots dx_{T_k}^\kappa,
  \end{align*}
 where $L_k'\in \mathcal{S}_0$, $k(y,\tilde{c}_0)> 0$, $k^{\mathsf{af}}(\tilde{c}_1)> 0$ and $\gamma_k>\sqrt{6}/\delta_k$ with $\delta_k\in[0,\delta_k']$ and $\delta_k'>0$.
  \label{propisolk}
 \end{prop}

\proof For any $ \kappa\geq 2 $ and $(x_{T_k}^1,\ldots,x_{T_k}^\kappa)\in[0,1]^{2\kappa}$, let $\Omega_0^{\kappa k}$ be the event that 
$B(x_{T_k}^1,\ldots,x_{T_k}^\kappa)$ contains no active members in $x_{T_k}^{\kappa+1},\ldots,x_{T_k}^{s_{T_k}}$. Then we have
\begin{align*}
 \mathbb{P}(E_{1k}\cap\cdots\cap E_{\kappa k}) &=  \int_{[0,1]^{2\kappa}}  \mathbb{P}\Big(E_{1k}\cap\cdots\cap E_{\kappa k}\Big|  
  (x_{T_k}^1,\ldots,x_{T_k}^\kappa)\in[0,1]^{2\kappa} \Big) dx_{T_k}^1\cdots dx_{T_k}^\kappa  \\
  & \leq \int_{[0,1]^{2\kappa}}  \mathbb{P}(\Omega_0^{\kappa k}) dx_{T_k}^1\cdots dx_{T_k}^\kappa \\
 &=\int_{[0,1]^{2\kappa}}  \sum_{\ell=0}^{s_{T_k}-\kappa} \Big((1-p)^\ell \frac{(\int_{B(x_{T_k}^1,\ldots,x_{T_k}^\kappa)} \lambda_k^\star(x)
 dx)^\ell}{\ell!}\textrm{e}^{-\int_{B(x_{T_k}^1,\ldots,x_{T_k}^\kappa)} \lambda_k^\star(x) dx}
 \Big) \\
 &\;  \times \Big( \frac{(\int_{[0,1]^2\setminus B(x_{T_k}^1,\ldots,x_{T_k}^\kappa)} \lambda_k^\star(x) dx)^{s_{T_k}-\kappa-\ell}}{(s_{T_k}-\kappa-\ell)!}
 \textrm{e}^{-\int_{[0,1]^2\setminus B(x_{T_k}^1,\ldots,x_{T_k}^\kappa)} \lambda_k^\star(x) dx}\Big)  dx_{T_k}^1\cdots dx_{T_k}^\kappa,
\end{align*}
and the claim follows easily by similar arguments as the end of proof of Proposition \ref{propisol}. \carre

The elicitation of the probability of several isolated members is somewhat more difficult than the case of one isolated member. Then, attention shows that the 
elicitation of the event $ \Omega_0^{\kappa k}$ used in the proof of Proposition \ref{propisolk} gives only an upper bound. 
Unfortunately, this seems to be difficult to prove in a general setting. To establish an asymptotic  
expression for this probability we need a deeper development and more notations. Let $\mathbb{S}_{a_\mathsf{f}}(x^1,\ldots,x^\kappa)$ denotes the sub-network 
over $(x^1,\ldots,x^\kappa)$ in which two members are connected (by affinity) if and only if their distance is at most $a_\mathsf{f}$. For any integer 
$\mathbf{n}$ satisfying $1\leq  \mathbf{n} \leq \kappa$, we denote by $\mathcal{C}_{\kappa\mathbf{n}}$ the set of $\kappa$-tuples 
$(x^1,\ldots,x^\kappa)\in[0,1]^{2\kappa}$ 
satisfying that $\mathbb{S}_{2a_\mathsf{f}}(x^1,\ldots,x^\kappa)$ has exactly $\mathbf{n}$ connected components. Note that the set $\mathcal{C}_{\kappa\kappa}$
consists of those tuples of $\kappa$ points which satisfies, for $i=1,\ldots,\kappa$, $B(x^i)$ contains none of the other points of the tuple. 
In the following result, we derive an interesting formula for the computation of the probability of several isolated members.
 \begin{prop}
 Under the assumption $H_2$ and  for any $ \kappa\geq 2 $ and $(x_{T_k}^1,\ldots,x_{T_k}^\kappa)\in\mathcal{C}_{\kappa\kappa}$ with 
 $v_r,d_r,A_\mathsf{f},\sigma,a_\mathsf{f}>0$ 
 and $0<p\leq 1$, we have for $s_{T_k}$  sufficiently large and $k\in\{1,\ldots,N\}$
  \begin{align*}
   \mathbb{P}(E_{1k}\cap\cdots\cap E_{\kappa k}) &=  \frac{\big(v_r k(y,\tilde{c}_0) + d_r + A_\mathsf{f} L_k'\frac{\gamma_k^2\log(s_{T_k})k^{\mathsf{af}}(\tilde{c}_1)}{2} \big)^{s_{T_k}-\kappa}}{
   (s_{T_k}-\kappa)!} \emph{e}^{-\big(v_r k(y,\tilde{c}_0) + d_r + A_\mathsf{f} 
 L_k'\frac{\gamma_k^2\log(s_{T_k})k^{\mathsf{af}}(\tilde{c}_1)}{2} \big) } \\
 &\qquad \qquad  \times \int_{\mathcal{C}_{\kappa\kappa}} \Big(1-p\emph{vol}(x_{T_k}^1,\ldots,x_{T_k}^\kappa)\Big)^{s_{T_k}-\kappa}dx_{T_k}^1\cdots dx_{T_k}^\kappa,
  \end{align*}
 where $L_k'\in \mathcal{S}_0$, $k(y,\tilde{c}_0)> 0$, $k^{\mathsf{af}}(\tilde{c}_1)> 0$ and $\gamma_k>\sqrt{6}/\delta_k$ with $\delta_k\in[0,\delta_k']$ and $\delta_k'>0$.
  \label{propisolkeq}
 \end{prop}

 \proof For any $ \kappa\geq 2 $ and $(x_{T_k}^1,\ldots,x_{T_k}^\kappa)\in\mathcal{C}_{\kappa\kappa}$, let $\Omega_1^{\kappa k}$ be the event that, 
 for all $1\leq i\leq \kappa$, $B(x_{T_k}^i)$ contains no active members in $x_{T_k}^{\kappa+1},\ldots,x_{T_k}^{s_{T_k}}$. In addition, 
 let $\Omega_2^{\kappa k}$ be the event that, for all $1\leq i\leq \kappa$, $B(x_{T_k}^i)$ contains $\mathbf{n}_i$ inactive members and no active members 
 in $x_{T_k}^{\kappa+1},\ldots,x_{T_k}^{s_{T_k}}$. Thanks to the previous introducing two events, it follows that
\begin{align*}
 \mathbb{P}(E_{1k}\cap\cdots\cap E_{\kappa k}) &=  \int_{\mathcal{C}_{\kappa\kappa}}  \mathbb{P}\Big(E_{1k}\cap\cdots\cap E_{\kappa k}\Big|  
  (x_{T_k}^1,\ldots,x_{T_k}^\kappa)\in\mathcal{C}_{\kappa\kappa} \Big) dx_{T_k}^1\cdots dx_{T_k}^\kappa  \\
  & = \int_{\mathcal{C}_{\kappa\kappa}}  \mathbb{P}(\Omega_1^{\kappa k}) dx_{T_k}^1\cdots dx_{T_k}^\kappa 
  =\int_{\mathcal{C}_{\kappa\kappa}}  \sum_{\mathbf{n}_1+\cdots+\mathbf{n}_\kappa}^{s_{T_k}-\kappa}   \mathbb{P}(\Omega_2^{\kappa k}) dx_{T_k}^1\cdots dx_{T_k}^\kappa\\
   &=\int_{\mathcal{C}_{\kappa\kappa}}  \sum_{\mathbf{n}_1+\cdots+\mathbf{n}_\kappa=0}^{s_{T_k}-\kappa}  \Big\{ \prod_{i=1}^\kappa 
 \Big((1-p)^{\mathbf{n}_i} \frac{(\int_{B(x_{T_k}^i)} \lambda_k^\star(x)
 dx)^{\mathbf{n}_i}}{\mathbf{n}_i!}\textrm{e}^{-\int_{B(x_{T_k}^i)} \lambda_k^\star(x) dx}
 \Big) \Big\}\\
 &\;  \times \Big( \frac{(\int_{[0,1]^2\setminus B(x_{T_k}^1,\ldots,x_{T_k}^\kappa)} \lambda_k^\star(x) dx)^{s_{T_k}-\kappa-\sum_{i=1}^\kappa \mathbf{n}_i}}{
 (s_{T_k}-\kappa- \sum_{i=1}^\kappa \mathbf{n}_i)!}
 \textrm{e}^{-\int_{[0,1]^2\setminus B(x_{T_k}^1,\ldots,x_{T_k}^\kappa)} \lambda_k^\star(x) dx}\Big)  dx_{T_k}^1\cdots dx_{T_k}^\kappa 
 \end{align*}
 \begin{align*}
 &=\int_{\mathcal{C}_{\kappa\kappa}} \frac{\textrm{e}^{-\big(v_r k(y,\tilde{c}_0) + d_r + A_\mathsf{f} L_k'\frac{\gamma_k^2\log(s_{T_k})k^{\mathsf{af}}(\tilde{c}_1)}{2} \big)
 \big(\sum_{i=1}^\kappa\textrm{vol}(x_{T_k}^i) +1-\textrm{vol}(x_{T_k}^1,\ldots,x_{T_k}^\kappa) \big) }}{(s_{T_k}-\kappa)!} \\
 & \times \sum_{\mathbf{n}_1+\cdots+\mathbf{n}_\kappa=0}^{s_{T_k}-\kappa}  
 \binom {s_{T_k}-\kappa} {\mathbf{n}_1,\ldots,\mathbf{n}_\kappa,s_{T_k}-\kappa- \sum_{i=1}^\kappa \mathbf{n}_i} \\
 &\times \Big\{ \prod_{i=1}^\kappa \Big( (1-p) \big(v_r k(y,\tilde{c}_0) + d_r + A_\mathsf{f} L_k'\frac{\gamma_k^2\log(s_{T_k})k^{\mathsf{af}}(\tilde{c}_1)}{2} \big) 
 \big(\textrm{vol}(x_{T_k}^i)\big)\Big)^{\mathbf{n}_i}\Big\}  \\
 &\times \Big( \big(v_rk(y,\tilde{c}_0)  + d_r + A_\mathsf{f} L_k'\frac{\gamma_k^2\log(s_{T_k})k^{\mathsf{af}}(\tilde{c}_1)}{2} \big) \\
& \times  \big(1-\textrm{vol}(x_{T_k}^1,\ldots,x_{T_k}^\kappa)\big)\Big)^{s_{T_k}-\kappa-\sum_{i=1}^\kappa \mathbf{n}_i} dx_{T_k}^1\cdots dx_{T_k}^\kappa \\
 &= \frac{\big(v_r k(y,\tilde{c}_0) + d_r + A_\mathsf{f} L_k'\frac{\gamma_k^2\log(s_{T_k})k^{\mathsf{af}}(\tilde{c}_1)}{2} \big)^{s_{T_k}-\kappa}}{(s_{T_k}-\kappa)!}
 \textrm{e}^{-\big(v_r k(y,\tilde{c}_0) + d_r + A_\mathsf{f} L_k'\frac{\gamma_k^2\log(s_{T_k})k^{\mathsf{af}}(\tilde{c}_1)}{2} \big) } \\
 &\qquad \qquad \times  \int_{\mathcal{C}_{\kappa\kappa}} \Big(1-p\textrm{vol}(x_{T_k}^1,\ldots,x_{T_k}^\kappa)\Big)^{s_{T_k}-\kappa}
 dx_{T_k}^1 \cdots dx_{T_k}^\kappa,
\end{align*}
where the last line is obtained from the multinomial theorem and by remarking that, for any set 
$(x_{T_k}^1,\ldots,x_{T_k}^\kappa)\in\mathcal{C}_{\kappa\kappa}$, we have $\sum_{i=1}^\kappa \textrm{vol}(x_{T_k}^i)=\textrm{vol}(x_{T_k}^1,\ldots,x_{T_k}^\kappa)$. 
This completes the proof. \carre

 Our previous result (Proposition \ref{propisolkeq}) gives an asymptotic equivalence for the probability of several isolated members not in the whole 
 domain of the torus but in a more restricted domain given by $(x_{T_k}^1,\ldots,x_{T_k}^\kappa)\in\mathcal{C}_{\kappa\kappa}$ where $\kappa\geq 2$. Even 
 if this result appears somewhat restrictive, we will show in the sequel that the probability in the domain 
 $[0,1]^{2\kappa} \setminus \mathcal{C}_{\kappa\kappa}$ (which is not covered by our result) converges asymptotically to nothing.  
 
 Let us highlight the major formulas of the previous results. Let denote 
 \begin{align}
  \alpha_k^\star(s_{T_k}) = \big(v_rk(y,\tilde{c}_0) + d_r + A_\mathsf{f} L_k'\frac{\gamma_k^2\log(s_{T_k})k^{\mathsf{af}}(\tilde{c}_1)}{2} \big),
  \label{alphastar}
 \end{align}
which may be infinite. From Proposition \ref{propisol}, we have asymptotically
\begin{align}
   \mathbb{P}(E_{ik}) \propto  \frac{\big(\alpha_k^\star(s_{T_k}) \big)^{s_{T_k}-1}}{(s_{T_k}-1)!} \textrm{e}^{-\big(\alpha_k^\star(s_{T_k}) \big) } ,
  \end{align}
 and we know from Proposition \ref{propisolkeq} that
 \begin{align}
   \mathbb{P}(E_{1k}\cap\cdots\cap E_{\kappa k}) \propto 
   \frac{\big(\alpha_k^\star(s_{T_k}) \big)^{s_{T_k}-\kappa}}{(s_{T_k}-\kappa)!} \textrm{e}^{-\big(\alpha_k^\star(s_{T_k}) \big) } .
  \end{align}

 The numerical problem that we aim to state now (to establish the distribution of isolated members) is the following; We have to constrain $\alpha^\star$ to fulfill the 
 following equation (needed for the proof
of the law of $N_0^k$  based on a version of Brun's sieve theorem and Bonferroni inequalities)
 \begin{align}
 \frac{\big(\alpha_k^\star(s_{T_k}) \big)^{s_{T_k}-\kappa}}{(s_{T_k}-\kappa)!} \textrm{e}^{-\big(\alpha_k^\star(s_{T_k}) \big) } = 
 \Big(\frac{\big(\alpha_k^\star(s_{T_k}) \big)^{s_{T_k}-1}}{(s_{T_k}-1)!} \textrm{e}^{-\big(\alpha_k^\star(s_{T_k}) \big) }  \Big)^\kappa,
 \label{constraint}
 \end{align}
 which is straightforwardly  solved by 
 \begin{align}
  \alpha_k^\star(s_{T_k}) = - s_{T_k} W\Big(-\frac{\big( (s_{T_k}-1)! ((s_{T_k}-\kappa)!)^{\frac{-1}{\kappa}}\big)^{\frac{\kappa}{(\kappa-1)s_{T_k}}}}{ s_{T_k}} \Big),
 \label{solconstraint}
 \end{align}
where $W(\cdot)$ is the Lambert function. From now on we impose an additional constraint on the constants $v_r,d_r,A_\mathsf{f},\sigma,\gamma_k$ in order that 
$\alpha_k^\star$ fulfills (\ref{constraint}) together with of course the fact that $v_r,d_r,A_\mathsf{f},\sigma>0$ and $\gamma_k>\sqrt{6}/\delta_k$ with $\delta_k\in[0,\delta_k']$ 
and $\delta_k'>0$. One may easily verify the constraint (\ref{constraint}) if we fix the intrinsic network parameters $v_r,d_r,A_\mathsf{f},\sigma$ and we 
compute the needed value of $\gamma_k$ from formula (\ref{solconstraint}).  We state this more precisely in the next Theorem \ref{theopois} below, which requires also the assumption $H_1$. Our 
main result (its  proof is delayed until Section \ref{proofs}) is as follows.

\begin{theo}
 Assume that hypothesis $H_1$ and $H_2$ are satisfied. Let 
 \begin{align}
  a_\mathsf{f}^\star(s_{T_k})= \Big(\frac{ \log\Big(\frac{s_{T_k} \pi}{4(s_{T_k}-1)!}\Big)+\mathfrak{f}_k({s_{T_k}})+ (s_{T_k}-1)\log(\alpha_k^\star(s_{T_k}))-
 \alpha_k^\star(s_{T_k}) }{\pi p s_{T_k}}\Big)^{1/2},
 \label{threshold}
 \end{align} 
 for any $p\in(0,1]$ and where $\alpha_k^\star$ is given by (\ref{alphastar}) and satisfying (\ref{constraint}). Then, we have  
 \begin{align}
  s_{T_k}  \Big(\frac{\big(\alpha_k^\star(s_{T_k}) \big)^{s_{T_k}-1}}{(s_{T_k}-1)!} \emph{e}^{-\big(\alpha_k^\star(s_{T_k}) \big) }  \Big)  
  \int_{[0,1]^2} \emph{e}^{-s_{T_k} p \emph{vol}(x)}dx  \longrightarrow   \beta_k,
 \label{beta}
 \end{align}
 as $s_{T_k}\to \infty$ and where 
 \begin{align*}
  \beta_k =   \emph{e} ^{-\displaystyle\lim_{s_{T_k}\to \infty} \mathfrak{f}_k({s_{T_k}})}.
 \end{align*}
If $\beta_k\in(0,\infty)$, then as $s_{T_k}\to \infty$, we have for $m\in\mathbb{N}$ and $k\geq 1$ that
\begin{align}
 \mathbb{P}\Big(N_0^k=m\Big) \longrightarrow
 \emph{e}^{-\beta_k} \beta_k^m/m!.
 \label{pois}
\end{align}
If $\beta_k=0$, then $\mathbb{P}\Big(N_0^k=0\Big) \to 1$, and if $\beta_k=\infty$, then $\mathbb{P}\Big(N_0^k=m\Big) \to 0$ 
for all $m\in\mathbb{N}$.
 \label{theopois}
\end{theo}

The  affinity threshold (\ref{threshold}) established for the  Poisson point process (resultant from the three processes) that generates the data 
$\mathbf{r}_{T_k}=\{x_{T_k}^1,\ldots,x_{T_k}^{s_{T_k}}\}$, for $1\leq k\leq N$, looks somewhat complicated than the threshold obtained with uniform random points. The 
particularity of the threshold (\ref{threshold}) is that, in addition to the size of the data, it depends on the parameters of the Poisson process 
$(v_r,d_r,A_\mathsf{f},\sigma)$ and the parameters of the discretization $(\gamma_k,L_k')$.

We now discuss related work and open problems.  Note that the Poisson distribution (\ref{pois}) [but not with the same threshold (\ref{threshold})] was already 
proved by \cite{Penr16}  and \cite{Yi+al06} in the special case  of  points uniformly distributed (respectively in the unit square and in a disk 
of unit area). Here we are considering a much more general 
class of random point processes. The results of this  paper  goes  a  step  beyond  the  literature  in  that  it  considers a Poisson point 
process with connection (affinity) function that is zero beyond the optimal threshold (\ref{threshold})  in  the  same  model. To our best knowledge, this is 
the first work where  results about the optimal threshold and the distribution of isolated members are shown for geometric networks with points generated 
from  Poisson point measures (rather than the uniform random points considered usually in literature). Other ideas for the choose 
of connection functions are proposed in the literature that deals with the subject of the connectivity of random geometric graphs initiated by 
\cite{Det+Hen89}.  For instance, using a step connection function ($\mathds{1}_{[0,a_\mathsf{f}]}(\|x-y\|)$), \cite{Det+Hen89} showed that 
for any constant $C>0$, the disk graph on $s_{T_k}$ uniform random points  with connectivity threshold $\sqrt{(\ln(s_{T_k})+C)/(\pi s_{T_k}) }$ has no isolated nodes with probability 
$\exp(-e^{-C})$ when $s_{T_k}$ tends to infinity. The theory has been generalized after in the unit square $[0,1]^d$ (with $d\geq 2$) by \cite{Penr16} for a 
class of connection functions that decay exponentially in some fixed positive power of distance. In some applications, it is desirable to use the Rayleigh 
fading connection function given by $\exp(-\bar{\beta}(\|x-y\|/\bar{\rho})^{\bar{\gamma}})$ for some fixed positive $\bar{\beta},\bar{\rho},\bar{\gamma}>0$ 
(typically $\bar{\gamma}=2$). It would be interesting to try to extend our results to these connection functions but this would be a nontrivial task because 
the discretization method developed here can be quite hard to adaptation.

Another related problem is the connectivity of the network. It is known that the main obstacle to connectivity is the existence of isolated members. More 
clearly, for the geometric (Gilbert) graph  $G(\mathcal{X}_{s_{T_k}}, a_\mathsf{f}({s_{T_k}}))$ with vertex set $\mathcal{X}_{s_{T_k}}$ given by a set of 
$s_{T_k}$ independently uniformly distributed points in $[0,1]^d$ with $d\geq 2$, and with an edge included between each pair of vertices at distance at most 
$a_\mathsf{f}({s_{T_k}})$, \cite{Penr97} showed that the probability that the graph is disconnected but free of isolated vertices tends to zero as 
$s_{T_k}$ tends to infinity, for any choice of  $(a_\mathsf{f}({s_{T_k}}))_{s_{T_k}\in\mathbb{N}}$.  The same result happens with the Erd\"os-Rényi graph 
but the proof for the geometric graph is much harder as pointed out by \cite{Bol01}. More formally, for random graphs $G$, the number of isolated vertices 
(denoted $N_0(G)$) has  (asymptotically) a Poisson distribution, hence with $\mathcal{K}$ denoting the class of connected graphs we have 
$\mathbb{P}(G\in \mathcal{K})\sim \mathbb{P}(N_0(G)=0)\sim \exp(-\mathbb{E}[N_0])$ as $s_{T_k}\to\infty$. We would expect something similar to hold for our 
social network $(r_t)_{t\geq 0}$. Under additional assumption on $\mathfrak{f}_k$ the  connectivity result of graphs $G$ presented here might naturally 
be conjectured for $(r_t)_{t\geq 0}$ as follows
\begin{align*}
 \forall k\geq 1, \qquad \quad \mathbb{P}\Big(r_{T_k}\in\mathcal{K}\Big) \to \emph{e}^{-\beta_k} \qquad\qquad \textrm{as } s_{T_k}\to\infty,
\end{align*}
with $\emph{e}^{-\beta_k}$ interpreted as 0 for $\beta_k=\infty$. More generally, it would be of interest (in its own right) to extend this to the case 
of $(r_t)_{t\geq 0}$.

We conclude our results by  checking  a statement shown essentially that the active isolated members of the social network  enjoys also a Poisson 
approximation at each time but with a slightly different mean. The following theorem states this.
\begin{theo}
 Admit assumptions of Theorem \ref{theopois} and consider the affinity threshold (\ref{threshold}). If $\beta_k\in(0,\infty)$, then as 
 $s_{T_k}\to \infty$, we have for $m\in\mathbb{N}$ and $k\geq 1$ that
\begin{align}
 \mathbb{P}\Big(N_a^k=m\Big) \longrightarrow
 \emph{e}^{-p\beta_k} (p\beta_k)^m/m!,
 \label{pois}
\end{align}
where $p\in(0,1]$ and $N_a^k$ denoting the number of active isolated members at time $T_k$. If $\beta_k=0$, then $\mathbb{P}\Big(N_a^k=0\Big) \to 1$, 
and if $\beta_k=\infty$, then 
$\mathbb{P}\Big(N_a^k=m\Big) \to 0$ for all $m\in\mathbb{N}$.
 \label{theopoisac}
\end{theo}

In section \ref{proofs}, we prove Theorems \ref{theopois} and \ref{theopoisac}.

We finish this section with a short discussion about the limit of $\mathfrak{f}_k(\cdot)$. As we have seen this limit plays a crucial role in determining the Poisson 
distribution of isolated members in the network. We fix $k\geq 1$. From the formula (\ref{threshold}) of optimal threshold, we deduce the following limits:
\begin{align*}
 \lim_{s_{T_k}\to \infty} \log\Big(\frac{s_{T_k} \pi}{4(s_{T_k}-1)!}\Big) = -\infty ,\;  \lim_{s_{T_k}\to \infty}(s_{T_k}-1)\log(\alpha_k^\star(s_{T_k}))=+\infty 
 \textmd{ and } \lim_{s_{T_k}\to \infty} -\alpha_k^\star(s_{T_k}) = - \infty.
\end{align*}
Furthermore, by series expansion at $s_{T_k}=+\infty$, we find 
\begin{align*}
 \lim_{s_{T_k}\to \infty} \log\Big(\frac{s_{T_k} \pi}{4(s_{T_k}-1)!}\Big) + (s_{T_k}-1)\log(\alpha_k^\star(s_{T_k})) - \alpha_k^\star(s_{T_k}) = - \infty,
\end{align*}
and 
\begin{align*}
 \lim_{s_{T_k}\to \infty} \frac{\log\Big(\frac{s_{T_k} \pi}{4(s_{T_k}-1)!}\Big) + (s_{T_k}-1)\log(\alpha_k^\star(s_{T_k})) - \alpha_k^\star(s_{T_k})}{\pi p s_{T_k}}
 = -\infty.
\end{align*}
Then, in order that the  function $\psi_k$ satisfy the three conditions of assumption $H_1$, there exists only  function  
$\mathfrak{f}_k({s_{T_k}})$ such that its limit is $ +\infty $. The two others cases ($\lim_{s_{T_k}\to \infty} \mathfrak{f}_k({s_{T_k}}) = \{-\infty,\in(0,\infty)\}$) discussed in 
Theorems \ref{theopois} and \ref{theopoisac} are in practice impossible since there exists no functions with these limits that satisfy at the same time $H_1$. Hence, 
$\beta_k$ takes only the convenient value 0 and we conclude that there is no isolated members in the network with probability one as $s_{T_k} \to \infty$.

\section{Dynamic clustering with $a_\mathsf{f}^\star$ }
\label{sec:3}

We apply in this section our main result, Theorem \ref{theopois} of Section \ref{sec:2}, to the problem of members clustering is that of grouping similar communities (components) of 
the social network, and estimating these groups from the random points  $\mathbf{r}_{T_k}$ at each time $T_k$. When cluster
similarity is defined via latent models, in which groups are relative to a partition of the index set $\{1,\ldots , s_{T_k}\}$, the most natural 
clustering strategy is K-means. We explain why this strategy cannot lead to perfect cluster in our context (especially towards the detection of isolated 
members) and offer another strategy, based on the optimal threshold $a_\mathsf{f}^\star$, that can be viewed as a density-based spatial clustering of 
applications with noise \citep{Ester+al96}. We introduce a cluster separation method tailored to our random network. The clusters estimated by this method 
are shown to be adaptively from the data. We compare this method with appropriate K-means-type procedure, and show that the former outperforms the latter 
for  cluster with detection of isolated members.

The solutions to the problem of clustering are typically algorithmic and entirely data based. They include applications of K-means, spectral clustering, 
density-based spatial clustering or versions of them. The statistical properties of these procedures have received a very limited amount of investigation. It 
is not currently known what statistical cluster method can be estimated by these popular techniques, or by their modifications. We try here to offer an 
answer to this question for the case of random points data issued from Poisson point process.

To describe our procedure, we begin by defining the function $\mathfrak{f}_k(\cdot)$  for all $k\geq 1$, for instance, by
\begin{align}
  \mathfrak{f}_k({s_{T_k}}) := (s_{T_k})^{l} -\log\Big(\frac{s_{T_k} \pi}{4(s_{T_k}-1)!}\Big) - (s_{T_k}-1)\log(\alpha_k^\star(s_{T_k})) +
 \alpha_k^\star(s_{T_k}) ,
 \label{choixf}
\end{align}
where $l\in(0,1)$. This gives $a_\mathsf{f}^\star(s_{T_k})= \sqrt{\frac{(s_{T_k})^{l} }{\pi p s_{T_k}}}$ and, with this particular choice, the threshold 
$a_\mathsf{f}^\star$ is not affected by the intrinsic parameters of the network ($v_r,d_r,A_\mathsf{f},\sigma$) and the parameters of the discretization 
$(\gamma_k,L_k')$. The clustering algorithm has three steps, and the main step 2 produces an estimator of one cluster from which we derive the estimated 
members (active and inactive) of this cluster. The three steps of the procedure are:
\begin{enumerate}
 \item Start with an arbitrary starting random point $x_{T_k}^i$ where $i \sim U_{\{1,\ldots , s_{T_k}\}}$.
 \item Compute $\rho(x_{T_k}^i)$ and this point's $a_\mathsf{f}^\star$-neighborhood is retrieved, and if it contains at least one active point, 
 except $x_{T_k}^i$ itself, a cluster is started. Otherwise, the point is labeled as isolated. Note 
 that we test if a point is active or not thanks to one realization of Bernoulli($p$). Note also that  this point might later be found in a sufficiently 
 sized $a_\mathsf{f}^\star$-neighborhood of a different point and hence be made part of a cluster. If a point is found to be a dense part of a cluster, its 
 $a_\mathsf{f}^\star$-neighborhood is also part of that cluster. Hence, all points that are found within the  $a_\mathsf{f}^\star$-neighborhood are added, as 
 is their own $a_\mathsf{f}^\star$-neighborhood when they are also dense. This process continues until the cluster is completely found.
 \item A new unvisited point is retrieved and processed, leading to the discovery of a further cluster or isolated point.
\end{enumerate}

The construction of an accurate function $\mathfrak{f}_k(\cdot)$ is a crucial step for guaranteeing the statistical optimality of the clustering procedure and 
for which the results of Theorems \ref{theopois} and  \ref{theopoisac} hold. Estimating $p$ (if it is unknown) before estimating the partition itself is a non-trivial task, and needs to be done with 
extreme care. The required inputs for Step 2 of the procedure are: $(i)$ $s_{T_k}$, the current size of the network; $(ii)$ $p$, the probability to be active 
in the network. Hence, this step is done at no additional accuracy cost. Remark  that, unlikely a K-means procedure, the $a_\mathsf{f}^\star$-neighborhood 
procedure does not require the  number of groups, which need an approach for selecting it in a data adaptive fashion.

\begin{figure}
\center
\subfigure[\it State of the network at $10^5$ iterations.]
    {\includegraphics[width=5.5cm,height=5cm]
          {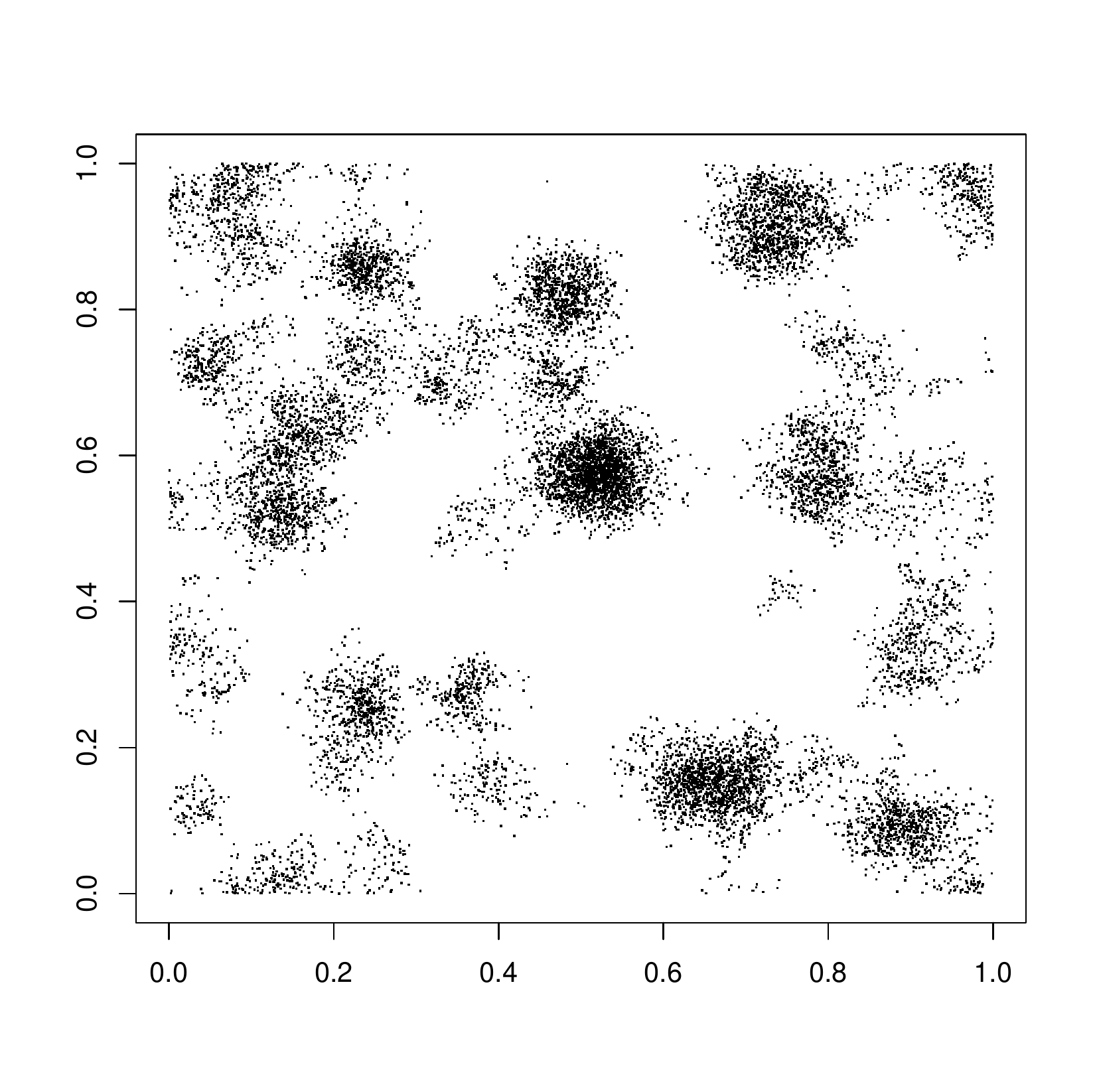}}        
\subfigure[\it Clustering with $a_\mathsf{f}^\star$-neighborhood procedure.]
    {\includegraphics[width=5.5cm,height=5cm]
                 {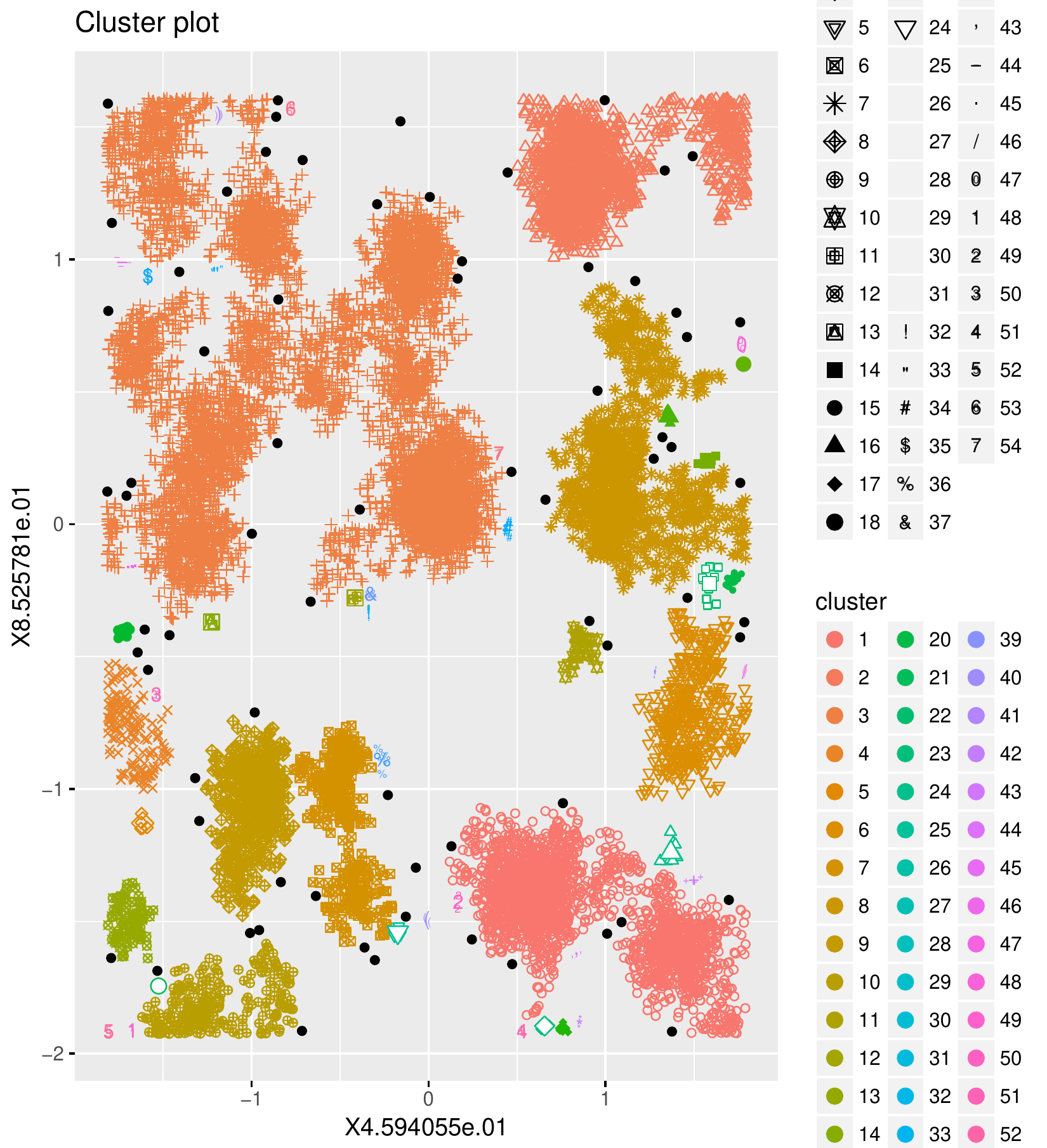}}                
\subfigure[\it Clustering with K-means procedure.]
    {\includegraphics[width=5.5cm,height=5cm]
                 {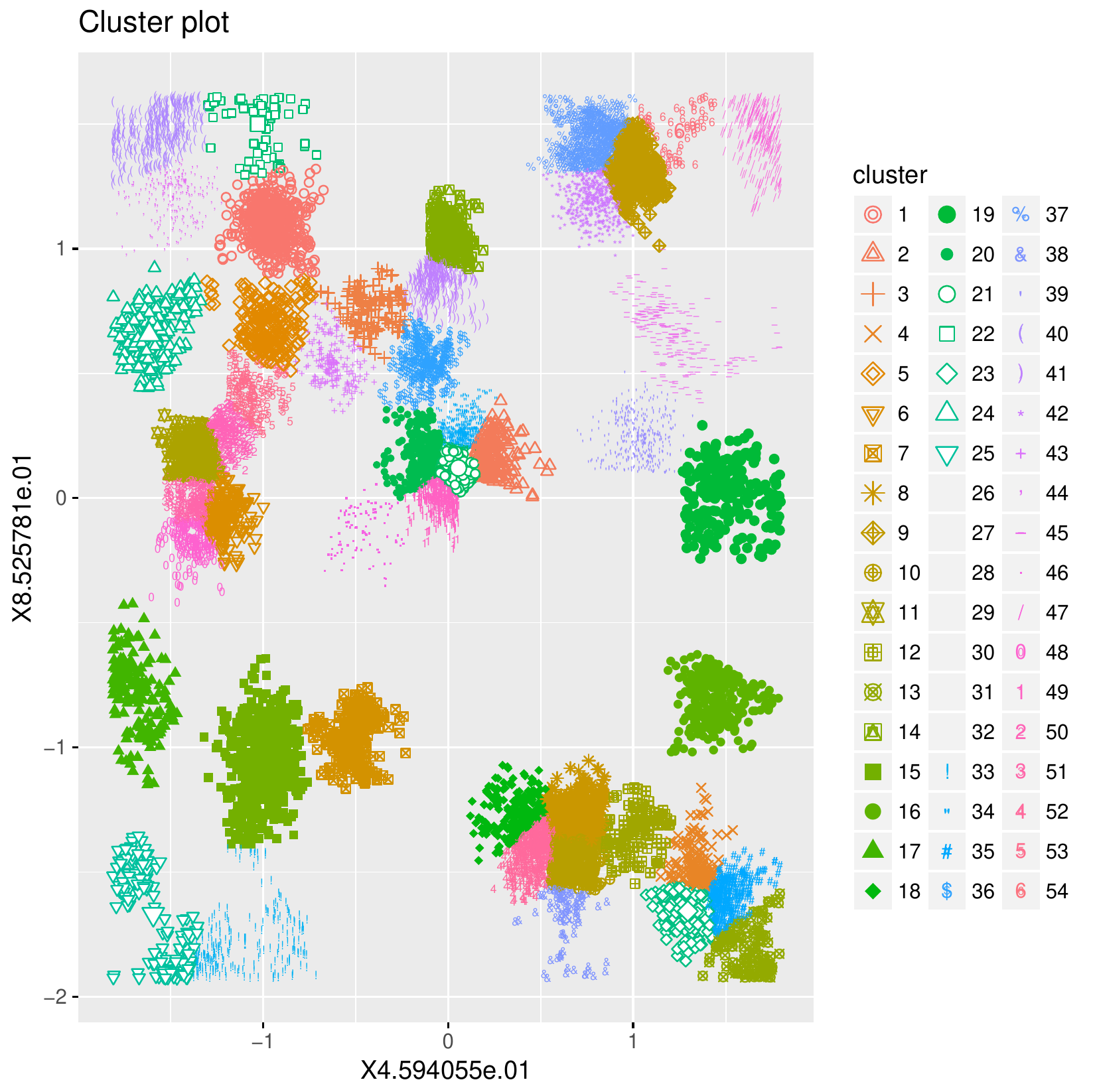}}             
\caption{{\small Figures shown a state of the network generated by Algorithm \ref{algo} with parameters $v_r=4$, $d_r=2$, $A_{\mathsf{f}}=2$, $a_\mathsf{f}=0.1$ and 
$\sigma=0.01$ and it's cluster analysis by two procedures.}}
\label{Figclus}
\end{figure}

Figure \ref{Figclus} shows a cluster analysis of a network at some time. The size of this network at $10^5$ updates is 14102 members. For the cluster analysis 
shown in Figure \ref{Figclus}(b), the function (\ref{choixf}) is used with $l=\frac{\log((a_\mathsf{f})^2 p\pi s_{T_k})}{\log(s_{T_k})}$  and the threshold 
$a_\mathsf{f}^\star$ is computed with $p=1$ (all  members are active) and $a_\mathsf{f}=0.1$. These gived $l = 0.6378047$. The number of isolated members finded by the $a_\mathsf{f}^\star$-neighborhood clustering procedure is 71 and the number of groups is 
 54. To compare this clustering result, we plot in Figure  \ref{Figclus}(c) a K-means procedure applied to the same network state by calibrating the 
 number of groups to 54 (to assess the difference between the two procedures at least by visual inspection). As seen, the clustering obtained now does not 
 detect isolated members. As known, the strategy of K-means is not based on the point's neighborhood and this has important repercussions on the analysis 
 and detection of isolated points in cluster estimates. The analysis of these geometric networks is non-standard, and needs to be done with care, as 
 illustrated by the proof of our Theorem \ref{theopois}. Moreover, in contrast to the $a_\mathsf{f}^\star$-neighborhood procedure tailored to our 
 spatial-temporal underlying model, K-means and spectral methods for this kind of models need to be corrected in a non-trivial 
fashion. 

\begin{figure}
\center
\subfigure[\it State of the network at $10^5$ iterations.]
    {\includegraphics[width=5.5cm,height=5cm]
          {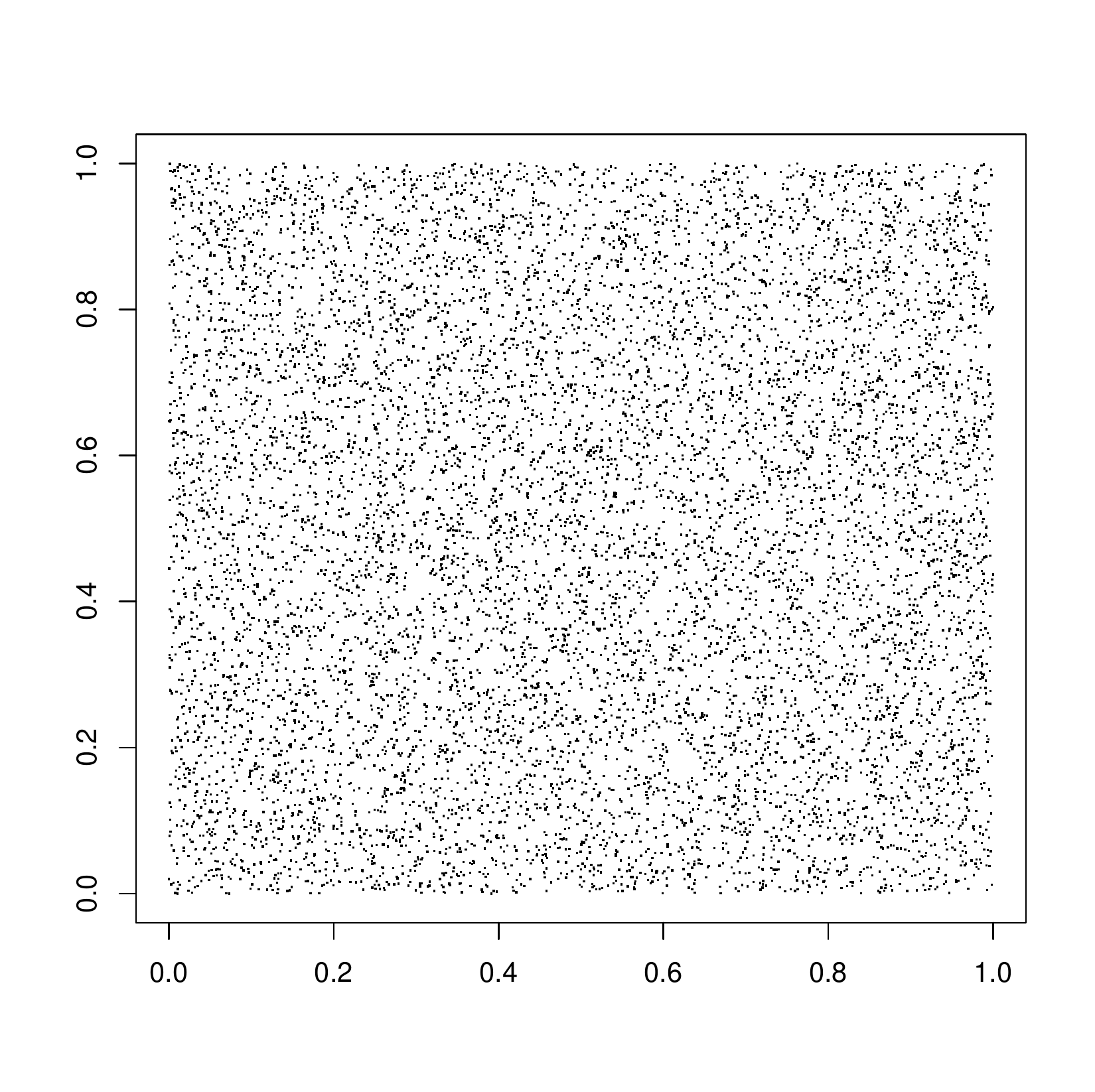}}        
\subfigure[\it Clustering with $a_\mathsf{f}^\star$-neighborhood procedure.]
    {\includegraphics[width=5.5cm,height=5cm]
                 {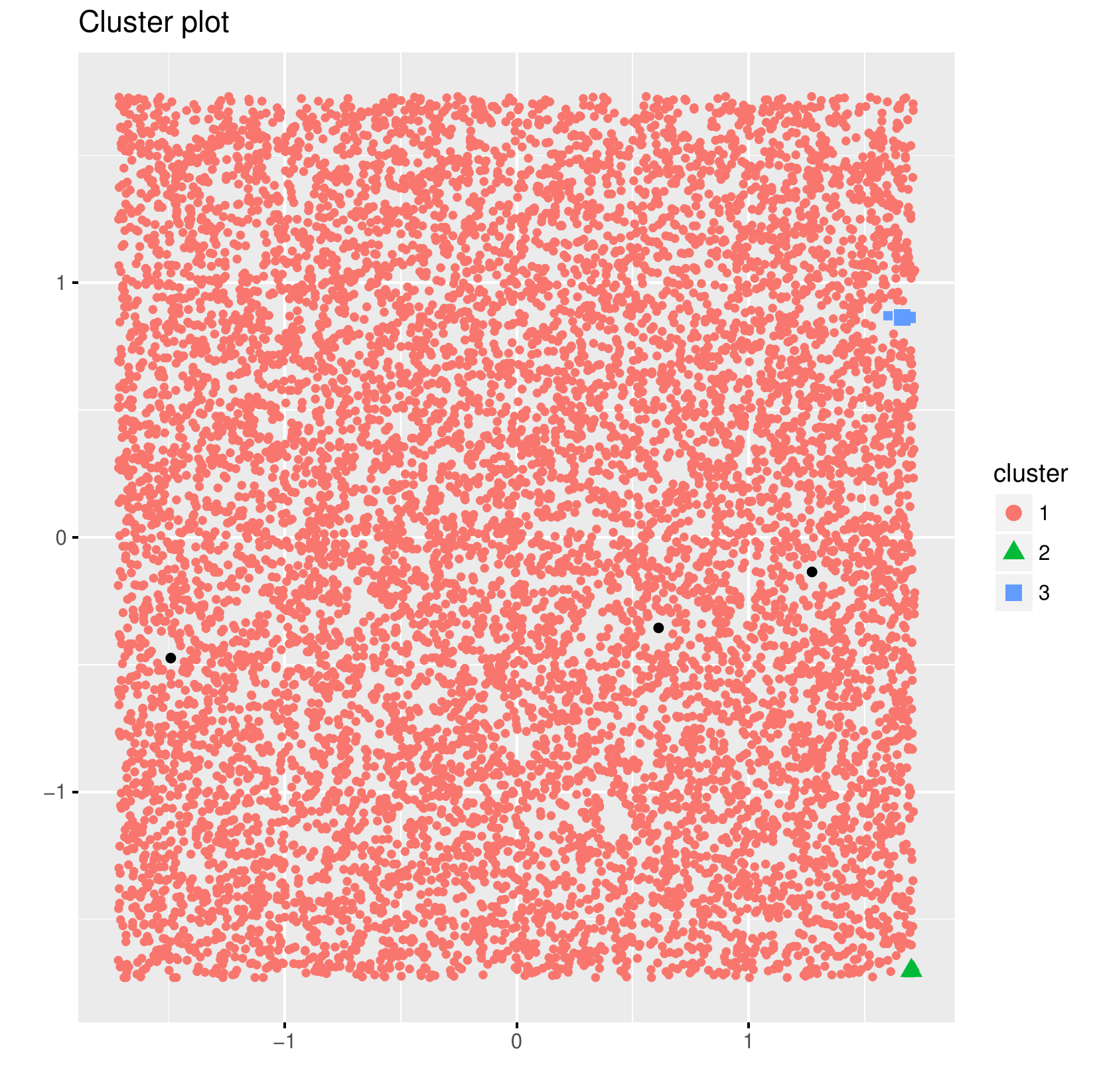}}                
\subfigure[\it Clustering with K-means procedure.]
    {\includegraphics[width=5.5cm,height=5cm]
                 {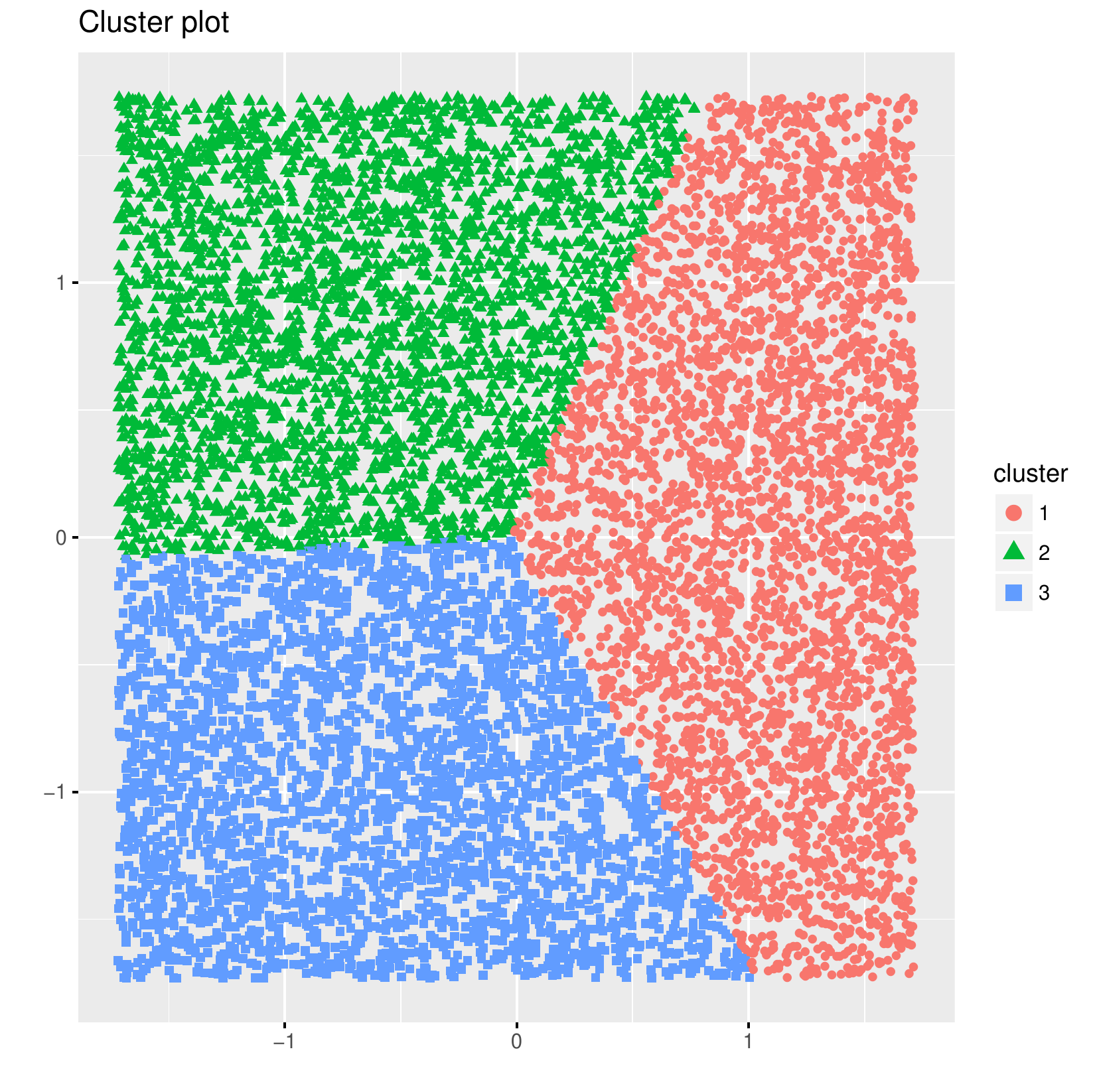}}             
\caption{{\small Figures shown a state of the network generated by Algorithm \ref{algo} with parameters $v_r=3$, $d_r=2$, $A_{\mathsf{f}}=2$, $a_\mathsf{f}=0.1$ and 
$\sigma=0.15$ and it's cluster analysis by two procedures.}}
\label{Figclus2}
\end{figure}

\begin{figure}
\center
\subfigure[\it State of the network at $10^5$ iterations.]
    {\includegraphics[width=5.5cm,height=5cm]
          {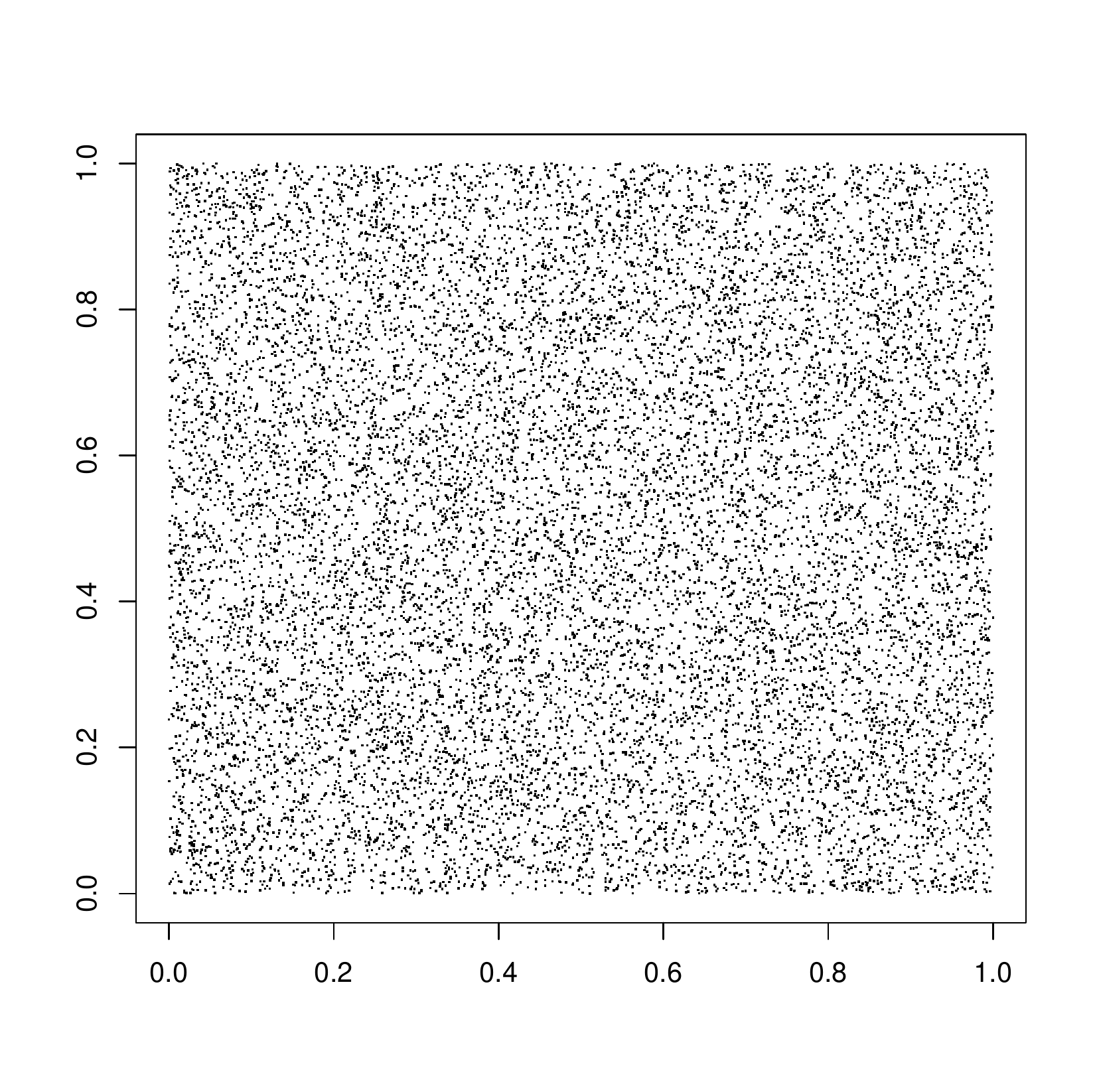}}        
\subfigure[\it Clustering with $a_\mathsf{f}^\star$-neighborhood procedure.]
    {\includegraphics[width=5.5cm,height=5cm]
                 {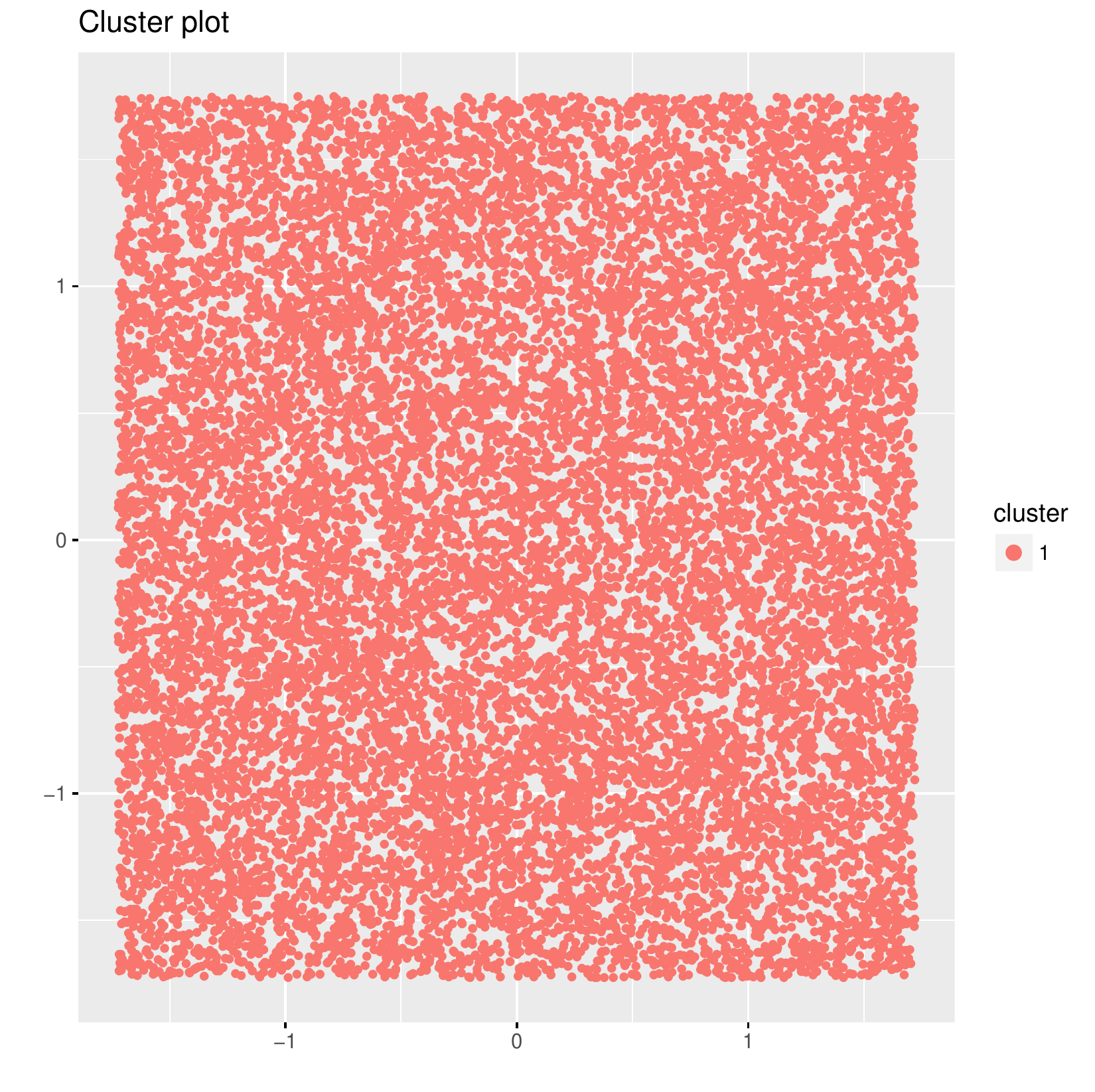}}                             
\caption{{\small Figures shown a state of the  network generated by Algorithm \ref{algo} with parameters $v_r=3.5$, $d_r=2$, $A_{\mathsf{f}}=2$, $a_\mathsf{f}=0.1$ and 
$\sigma=0.21$ and it's cluster analysis by the $a_\mathsf{f}^\star$-neighborhood procedure.}}
\label{Figclus3}
\end{figure}

Now, we consider another configuration when we increase the dispersion ($\uparrow\sigma$) in order that the members of the network tend to occupy almost all 
the space $[0,1]^2$. Figure \ref{Figclus2} shows a cluster analysis of a network with high dispersion ($\sigma=0.15$). A cluster analysis using the 
$a_\mathsf{f}^\star$-neighborhood  procedure (with again $p=1$) has found 3 isolated members and 3 groups (one very large and 
two very small), here $s_{T_k}=10125$ members. We report the clustering result obtained from the K-means procedure to highlight again that this type of 
clustering is not perfect for our kind of network modeling. Finally, we increase a little more  the dispersion to $\sigma=0.21$ and the invited rate to 
$v_r=4$ in order to strengthen the occupation of the space by the members of the network. We conclude  by observing that in Figure 
\ref{Figclus3} the isolated members have disappeared and the network forms one connected group which confirms, even numerically, our previous conjecture on 
the connectivity of the network ($\mathbb{P}(r_{T_k}\in\mathcal{K}) \to \emph{e}^{-\beta_k}=1$ as $s_{T_k}\to\infty$). 

\section{Discussion}
\label{sec:4}

In this framework, a geometric and dynamic social network constructed from Poisson point measures is investigated  by the distribution of isolated members 
in the network. We presume a social network is composed of particles in interaction in continuous time represented by a set of points over the unit square. We 
assume a Bernoulli probability on the activity state of each point in the network. We prove  that if all points have the same radius of affinity 
$a_\mathsf{f}^\star=\sqrt{\frac{(s_t)^l}{p\pi s_t}}$ for $l\in(0,1)$, the number of isolated points is asymptotically Poisson with mean $\beta_k$ 
where $\beta_k\sim 0$ (at each time $T_k$ for $k\geq 1$) as the current size $s_{T_k}$ of the network tends to infinity. This offers a natural threshold 
for the construction of a  $a_\mathsf{f}^\star$-neighborhood  procedure tailored to the dynamic clustering of the network adaptively from the data. The question 
whether vanishment of isolated members almost surely ensures connectivity of the network or not remains open.

\section{Proofs}
\label{proofs}

\begin{figure}[h]
\begin{center}
\begin{tikzpicture}
\draw (0, 0)--(5, 0)--(5, 5)--(0, 5)--(0, 0);
\draw (2.5, 2.5) circle (1.5cm);
\draw[black, ultra thick, dotted] (2.5, 2.5) circle (2.3cm);
\draw[red, thick, dashed] (2.5, 2.5) circle (2.5cm);
\draw[ultra thick, <-] (1.5, 4.8) -- (1.5, 5.5) node[above]{$\mathcal{C}_3$};
\draw[thick, <-] (2.8, 4.8) -- (2.8, 5.5) node[above]{$\mathcal{C}_2$};
\draw[thick, <-] (2.3, 4.0) -- (2.3, 5.5) node[above]{$\mathcal{C}_1$};
\draw[thick, <-] (4.5, 5.01) -- (4.5, 5.5) node[above]{$\bar{\mathcal{D}}=[0,1]^2$};
\node[black,very thick](A) at (2.5, 2.2) {$\mathcal{D}_1$};
\node[black,very thick](A) at (1.2, 1.2) {$\mathcal{D}_2$};
\node[black,very thick](A) at (0.4, 0.4) {$\mathcal{D}_4$};
\node[black,very thick](A) at (4.6, 0.4) {$\mathcal{D}_4$};
\node[black,very thick](A) at (4.6, 4.6) {$\mathcal{D}_4$};
\node[black,very thin](A) at (0.4, 4.6) {$\mathcal{D}_4$};
\node[black,very thick](A) at (-0.1, -0.1) {$0$};
\node[black,very thick](A) at (-0.1, 5.0) {$1$};
\node[black,very thick](A) at (5.1, -0.1) {$1$};
\draw[thick, <-] (4.35, 1.0) -- (5.2, 1.0) node[right]{$\mathcal{D}_3$};
\draw[.-.] (2.5, 2.5) -- (4.8, 2.5);
\draw[.-.] (4.8, 3.5) -- (4.8, 2.5);
\draw[.-.] (4.8, 1.5) -- (4.8, 2.5);
\draw[.-.] (0.0, 2.5) -- (1.0, 2.5);
\draw[.-.] (2.5, 2.5) -- (4.8, 3.5);
\draw[color=blue,decorate,decoration={brace,raise=0.05cm}](2.5,2.5) -- (4.8,3.5) node[above=0.05cm,pos=0.37] {$1/2$};
\draw[color=blue,decorate,decoration={brace,raise=0.05cm}](0,2.5) -- (1,2.5) node[above,pos=0.5] {$a_\mathsf{f}$};
\draw[color=blue,decorate,decoration={brace,raise=0.05cm}](4.8,3.5) -- (4.8,2.5) node[right=2,pos=0.5] {$a_\mathsf{f}$};
\draw[color=blue,decorate,decoration={brace,raise=0.05cm}](4.8,2.5) -- (4.8,1.5) node[right=2,pos=0.5] {$a_\mathsf{f}$};
\end{tikzpicture}
\caption{The different levels of partitioning of the torus $[0,1]^2$ are shown here with three circles $\mathcal{C}_1$ (--), $\mathcal{C}_2$ ($\cdots$) and 
$\mathcal{C}_3$ (- -). The torus is subdivided into 4 regions; $\mathcal{D}_1$ (given by the disk of radius $1/2-a_\mathsf{f}$ centred at $ \mathbb{O}=(1/2,1/2)$), 
$\mathcal{D}_2$ (given by the annulus of radii $1/2-a_\mathsf{f}$ and $\sqrt{1/4-a_\mathsf{f}^2}$ centered at $\mathbb{O}$), 
 $\mathcal{D}_3$ (given by the annulus of radii $\sqrt{1/4-a_\mathsf{f}^2}$ and $1/2$ centered at $\mathbb{O}$) and $\mathcal{D}_4$ (given by the rest of the torus).  }
\label{partition}
 \end{center}
\end{figure}
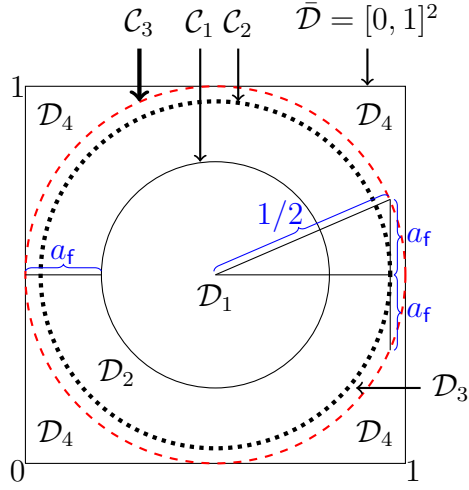

We split the proof of Theorem \ref{theopois}  into several lemmas. Before presenting the following statements, we partition the unit square $[0,1]^2$ into 
four regions $\mathcal{D}_1,\ldots,\mathcal{D}_4$ as explained in Figure \ref{partition}. From this partition, it is easy to see that for any $x$ in the torus we have 
$\pi a_\mathsf{f}^2/4 \leq \textrm{vol}(x)\leq \pi a_\mathsf{f}^2$ (in particular $\textrm{vol}(x)= \pi a_\mathsf{f}^2/4$ when $x$ is at the corners and 
$\textrm{vol}(x)= \pi a_\mathsf{f}^2$ for all $x\in\mathcal{D}_1$). 

In the following lemma, we shall study  a lower bound on $\textrm{vol}(x)$ for all $x\in\mathcal{D}_2$.
\begin{lem}
 Under the previous partition of the torus $[0,1]^2$ given in Figure \ref{partition}, then, for all $x\in\mathcal{D}_2$, we have
 \begin{align}
  \emph{vol}(x) \geq \frac{\pi a_\mathsf{f}^2}{2} + a_\mathsf{f}\Big( \frac{1}{2}- \eta \|x\| \Big), \qquad \textrm{ with } \qquad
  \eta\geq \frac{1}{\frac{\sqrt{1/2-2a_\mathsf{f}^2}}{1/2-2a_\mathsf{f}^2}-1},
 \label{lem1}
 \end{align}
 where $0< a_\mathsf{f} \leq 1/2$.
 \label{lemvol}
\end{lem}

\proof
Without loss of generality, let consider two points $x, y\in[0,1]^2$ such that $y$ be a point in $\partial\bar{\mathcal{D}}$ (the boundary of the torus) with coordinates 
$(0.5,0)$ or $(1,0.5)$ or $(0.5,1)$ or $(0,0.5)$. Straightforwardly, at these fourth coordinates the distance  between the circle $\mathcal{C}_2$ and 
the boundary of the torus achieves its minimum. It is sufficient to prove the lemma for one of these coordinates. The point $x$ be in $\partial \mathcal{C}_2$ 
together with $\|y-x\|=1/2-\|x-\mathbb{O}\|$ and $ab$ be the diameter of $B(x,a_\mathsf{f})$ 
perpendicular to $xy$ (see Figure \ref{coor}). 

\begin{figure}[h]
\begin{center}
\begin{tikzpicture}
\draw (5, 0)--(5, 5);
\node[black,very thick](A) at (5, -0.2) {$\vdots$};
\node[black,very thick](A) at (5, 5.4) {$\vdots$};
\draw (4, 2.5) circle (1.5cm);
\draw (4, 1)--(4, 4);
\draw (4, 2.5)--(5, 2.5);
\draw (4, 4)--(5, 2.5);
\draw (4, 1)--(5, 2.5);
\node[black,very thick](A) at (3.8, 2.5) {$x$};
\node[black,very thick](A) at (5.2, 2.5) {$y$};
\node[black,very thick](A) at (1.2, 2.2) {$\mathbb{O}$};
\node[black,very thick](A) at (1.2, 2.5) {$\bullet$};
\node[black,very thick](A) at (4, 4.2) {$a$};
\node[black,very thick](A) at (4, 0.8) {$b$};
\draw[thick, <-] (5, 4.5) -- (5.5, 4.5) node[right]{$\partial \bar{\mathcal{D}}$};
\draw[thick, <-] (5.4, 3.0) -- (5.9, 3.0) node[right]{$B(x,a_\mathsf{f})$};
\end{tikzpicture}
\caption{A partial configuration of the torus boundary and the ball $B(x,a_\mathsf{f})$ that contains the point $y=(1,0.5)$ and centered at 
$x=(0.5+\sqrt{1/4-a_\mathsf{f}^2},0.5)$. The volume $\textrm{vol}(x)$ 
is higher than the volume of the half of $B(x,a_\mathsf{f})$ plus the area of the triangle $aby$.}
\label{coor}
 \end{center}
\end{figure}

Thus, for all $x\in\mathcal{D}_2$, we can write
\begin{align}
 \textrm{vol}(x) \geq \frac{\pi a_\mathsf{f}^2}{2} + a_\mathsf{f}  \|y-x\|  = \frac{\pi a_\mathsf{f}^2}{2} + a_\mathsf{f}\Big( 1/2-\|x-\mathbb{O}\| \Big).
 \label{Circle2}
\end{align}
Now, for all $x\in \partial\mathcal{C}_2$ and from Figure \ref{partition}, we have 
\begin{align*}
 \inf_{x\in \partial\mathcal{C}_2} \|x\|=\frac{\sqrt{2}}{2} - \sqrt{1/4-a_\mathsf{f}^2}=\frac{\sqrt{2}}{2}-\|x-\mathbb{O}\|  \leq \|x\|.
\end{align*}
One may easily find an upper bound for $\|x-\mathbb{O}\|$ in function of $\|x\|$ by assuming that exist $\eta\in(0,\infty)$ such that 
$\|x-\mathbb{O}\|\leq \eta \inf_{x\in \partial\mathcal{C}_2}\|x\|$ and by easy  calculation we find
\begin{align*}
 \|x-\mathbb{O}\|  \leq \eta \Big( \frac{\sqrt{2}}{2}-\|x-\mathbb{O}\| \Big) \qquad & \Leftrightarrow   \qquad  \sqrt{1/4-a_\mathsf{f}^2} \leq \eta 
 \Big( \frac{\sqrt{2}}{2} - \sqrt{1/4-a_\mathsf{f}^2}  \Big) \\
 & \Leftrightarrow   \eta\geq \frac{1}{\frac{\sqrt{1/2-2a_\mathsf{f}^2}}{1/2-2a_\mathsf{f}^2}-1}.
\end{align*}
Hence, we complete the proof by  using  $\|x-\mathbb{O}\| \leq \eta \inf_{x\in \partial\mathcal{C}_2}\|x\| \leq \eta \|x\|$ in  (\ref{Circle2}) which leads to (\ref{lem1}).

\carre 

We may now formulate a lower bound on the volume of more than one member inside the network, for which the proof follows the same geometric spirit than 
Lemma \ref{lemvol}

\begin{lem}
 Admit the partition of the torus $[0,1]^2$ given in Figure \ref{partition}. For any sequence $(x^1,\ldots,x^\kappa)\in[0,1]^{2\kappa}$ where $\kappa\geq 2$ and such 
 that $x^1$ has the largest norm with $\|x^i - x^j\|\leq 2a_\mathsf{f}$ if and only if $|i-j|\leq 1$, we have
\begin{align}
 \emph{vol}(x^1,\ldots,x^\kappa) \geq \emph{vol}(x^1) + \frac{\pi a_\mathsf{f}}{16} \sum_{i=1}^{\kappa-1} \Delta x^i,
 \label{volkappa}
\end{align}
where $\Delta x^i=\|x^{i+1} - x^i\|$.
\label{lemsev}
\end{lem}

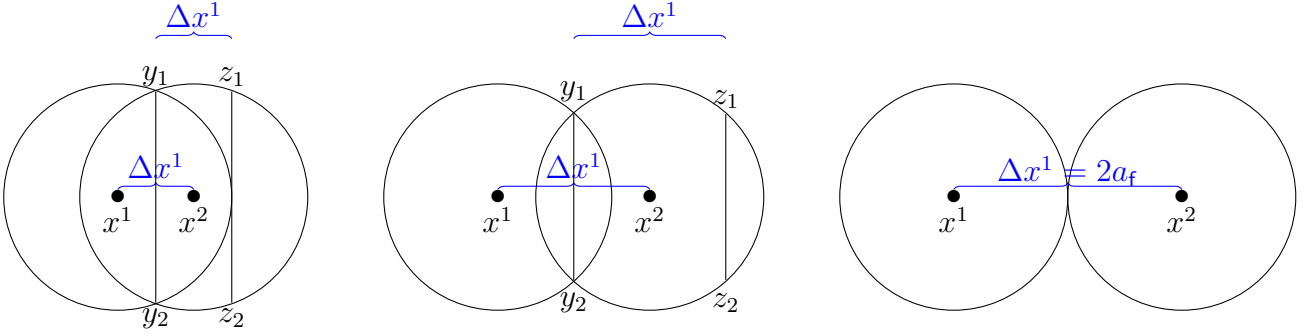
\begin{figure}[h]
\begin{center}
\begin{tikzpicture}

\draw (-05, 2.5) circle (1.5cm);
\draw (-4, 2.5) circle (1.5cm);
\draw (0, 2.5) circle (1.5cm);
\draw (2, 2.5) circle (1.5cm);
\draw (6, 2.5) circle (1.5cm);
\draw (9, 2.5) circle (1.5cm);
\node[black,very thick](A) at (-5, 2.5) {$\bullet$};
\node[black,very thick](A) at (-5, 2.2) {$x^1$};
\node[black,very thick](A) at (-4, 2.5) {$\bullet$};
\node[black,very thick](A) at (-4, 2.2) {$x^2$};
\node[black,very thick](A) at (0, 2.5) {$\bullet$};
\node[black,very thick](A) at (0, 2.2) {$x^1$};
\node[black,very thick](A) at (2, 2.5) {$\bullet$};
\node[black,very thick](A) at (2, 2.2) {$x^2$};
\node[black,very thick](A) at (6, 2.5) {$\bullet$};
\node[black,very thick](A) at (6, 2.2) {$x^1$};
\node[black,very thick](A) at (9, 2.5) {$\bullet$};
\node[black,very thick](A) at (9, 2.2) {$x^2$};
\draw[color=blue,decorate,decoration={brace,raise=0.1cm}](6,2.5) -- (9,2.5) node[above=0.1,pos=0.5] {$ \Delta x^1=2a_\mathsf{f}$};
\draw[color=blue,decorate,decoration={brace,raise=0.1cm}](0,2.5) -- (2,2.5) node[above=3,pos=0.5] {$ \Delta x^1$};
\draw[color=blue,decorate,decoration={brace,raise=0.1cm}](-5,2.5) -- (-4,2.5) node[above=3,pos=0.5] {$ \Delta x^1$};
\draw (-4.5, 3.9)--(-4.5, 1.1);
\node[black,very thick](A) at (-4.5, 4.1) {$y_1$};
\node[black,very thick](A) at (-4.5, 0.9) {$y_2$};
\draw (-3.5, 3.9)--(-3.5, 1.1);
\node[black,very thick](A) at (-3.5, 4.1) {$z_1$};
\node[black,very thick](A) at (-3.5, 0.9) {$z_2$};
\draw[color=blue,decorate,decoration={brace,raise=0.1cm}](-4.5,4.5) -- (-3.5,4.5) node[above=3,pos=0.5] {$ \Delta x^1$};
\draw (1, 3.6)--(1, 1.4);
\node[black,very thick](A) at (1, 3.9) {$y_1$};
\node[black,very thick](A) at (1, 1.1) {$y_2$};
\draw (3, 3.6)--(3, 1.4);
\node[black,very thick](A) at (3, 3.8) {$z_1$};
\node[black,very thick](A) at (3, 1.1) {$z_2$};
\draw[color=blue,decorate,decoration={brace,raise=0.1cm}](1,4.5) -- (3,4.5) node[above=3,pos=0.5] {$ \Delta x^1$};
\end{tikzpicture}
\caption{Three examples of configuration for the chords $y_1y_2$ and $z_1z_2$ (of two intersecting balls $ B(x^1,a_\mathsf{f}) $ and 
$ B(x^2,a_\mathsf{f}) $), where we see that $\|y_1-y_2\|$ is decreasing over $[0,2a_\mathsf{f}] $.}
\label{chords}
 \end{center}
\end{figure}

\proof We show the result  (\ref{volkappa}) for $\kappa=2$ at first. Then, we prove by induction that the result  holds for any $\kappa\geq 2$. For simplicity, 
let fix $\kappa=2$  and consider the function $\varphi(\Delta x^1)=\textrm{vol}\big(B(x^2,a_\mathsf{f}) \setminus B(x^1,a_\mathsf{f})\big)$. Our aim is to prove that 
$\varphi(\Delta x^1)\geq \frac{\pi a_\mathsf{f}}{2}\Delta x^1$ for any $x^1,x^2\in[0,1]^2$ with $\|x^1\| \geq \|x^2\|$. To do this, let $y_1y_2$ be the common chord of $\partial B(x^1,a_\mathsf{f}) $ and 
$\partial B(x^2,a_\mathsf{f}) $. Let also $z_1z_2$ be another chord of $\partial B(x^2,a_\mathsf{f}) $ that is parallel to $y_1y_2$ and has the same length as $y_1y_2$ as 
explained in Figure \ref{chords}.

\begin{figure}[h]
\begin{center}
\begin{tikzpicture}

\draw (0, 2.5) circle (1.5cm);
\draw (2, 2.5) circle (1.5cm);

\node[black,very thick](A) at (0, 2.55) {$\bullet$};
\node[black,very thick](A) at (0, 2.8) {$x^1$};
\node[black,very thick](A) at (2, 2.55) {$\bullet$};
\node[black,very thick](A) at (2, 2.8) {$x^2$};
\draw[color=blue,decorate,decoration={brace,raise=0.1cm}](0,2.5) -- (2,2.5) node[above=3,pos=0.5] {$ \Delta x^1$};
\draw[thick, <-] (2, 0) -- (2.5, 0) node[right]{$\partial \bar{\mathcal{D}}$};
\draw[thick, <-] (-3, 2.55) -- (-3, 3) node[above]{$\partial \bar{\mathcal{D}}$};
\draw (1, 3.6)--(1, 1.4);
\draw (-4, 2.55)--(2, 2.55);
\draw (2, 2.55)--(2, -1);
\node[black,very thick](A) at (2, -1.3) {$\vdots$};
\node[black,very thick](A) at (-4.4, 2.55) {$\cdots$};
\node[black,very thick](A) at (1, 3.9) {$y_1$};
\node[black,very thick](A) at (1, 1.1) {$y_2$};
\draw (3, 3.6)--(3, 1.4);
\node[black,very thick](A) at (3, 3.8) {$z_1$};
\node[black,very thick](A) at (3, 1.1) {$z_2$};
\draw[color=blue,decorate,decoration={brace,raise=0.1cm}](1,4.5) -- (3,4.5) node[above=3,pos=0.5] {$ \Delta x^1$};
\end{tikzpicture}
\caption{An example of configuration for the chords $y_1y_2$ and $z_1z_2$ (of two intersecting balls $ B(x^1,a_\mathsf{f}) $ and 
$ B(x^2,a_\mathsf{f}) $), where $x^2$ is at the corner of $[0,1]^2$.}
\label{corner}
 \end{center}
\end{figure}
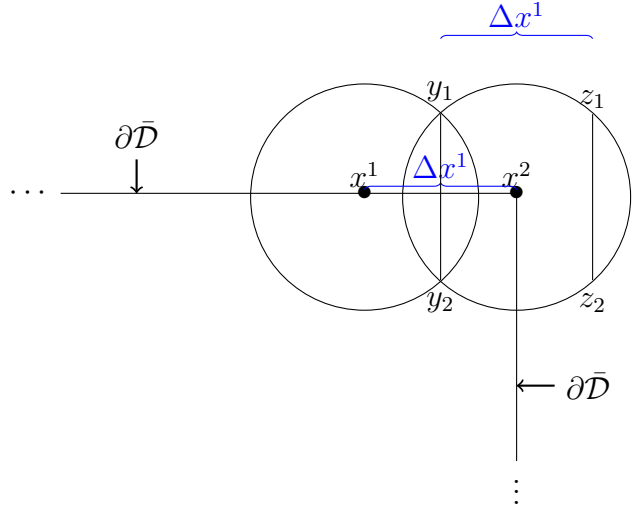

 From Figure \ref{chords}, its clear that $\varphi(\Delta x^1)$ is equal to the volume of the portion of $ B(x^2,a_\mathsf{f}) $ between the chords $y_1y_2$ and $z_1z_2$ and 
 immediately we deduce that the second derivative $\varphi''(\Delta x^1)\leq 0$ (since $\varphi'(\Delta x^1)=\|y_1-y_2\|$ which is decreasing). Hence, the function  $\varphi$ is concave 
 with  $\varphi(0)=0$ and $\varphi(2a_\mathsf{f})=\pi a_\mathsf{f}^2$ which enables us to write
 \begin{align*}
 \forall \Delta x, \Delta x'\in [0,2a_\mathsf{f}] \quad \textrm{ and }  \quad \forall \ell^0\in[0,1], \qquad         \varphi\big(\ell^0 \Delta x  +(1-\ell^0) \Delta x'\big)  \geq 
 \ell^0 \varphi(\Delta x ) +(1-\ell^0)\varphi(\Delta x'),
 \end{align*}
and by taking $\Delta x=0$ and $\Delta x' =2a_\mathsf{f}$,
\begin{align}
 \varphi\big((1-\ell^0)2a_\mathsf{f}\big) \geq  (1-\ell^0) \pi a_\mathsf{f}^2,
 \label{concave}
\end{align}
where $\varphi(\Delta x^1)\geq \frac{\pi a_\mathsf{f}}{2}\Delta x^1$ holds by choosing $\ell^0$ in (\ref{concave}) such that $\Delta x^1=(1-\ell^0)2a_\mathsf{f}$. We discuss now the lower 
bound of $\textrm{vol}(x^1,x^2) - \textrm{vol}(x^1)$ following the 
position of $x^1$ in the partition of $[0,1]^2$. Indeed, if $x^1\in \mathcal{D}_1$, then the two balls $ B(x^1,a_\mathsf{f}) $ and 
$ B(x^2,a_\mathsf{f}) $ are completely contained in $[0,1]^2$ which enables us to obtain
\begin{align*}
 \textrm{vol}(x^1,x^2) - \textrm{vol}(x^1) = \textrm{vol}\big(B(x^2,a_\mathsf{f}) \setminus B(x^1,a_\mathsf{f})\big) =\varphi(\Delta x^1) \geq 
 \frac{\pi a_\mathsf{f}}{2}\Delta x^1,
\end{align*}
and the lemma is proved for $\kappa=2$ and if $x^1\in \mathcal{D}_1$. If $x^1\notin \mathcal{D}_1$, we remark that for the same distance $\Delta x^1$, the value of 
$\textrm{vol}(x^1,x^2) - \textrm{vol}(x^1)$ achieves its minimum if $x^1,x^2\in \partial\bar{\mathcal{D}}$ together with  $x^2$ is at the corner.  It suffice 
to show the lemma in this case. For the easy of exposition, let consider again the previous chords such that $y_2\in [0,1]^2$  as shown in Figure \ref{corner}.

 From Figure \ref{corner}, it appeared that $ \textrm{vol}(x^1,x^2) - \textrm{vol}(x^1)= \frac{\varphi(\Delta x^1)}{4} \geq \frac{\pi a_\mathsf{f}}{8}\Delta x^1$ and the 
 lemma follows also in this case. Now, let consider $\kappa=3$ and write
 \begin{align*}
  \textrm{vol}(x^1,x^2,x^3) &\geq \textrm{vol}(x^1) + \textrm{vol}(x^3)  \geq \textrm{vol}(x^1) + \frac{\pi a_\mathsf{f}^2}{4}\\
  &= \textrm{vol}(x^1) + \frac{\pi a_\mathsf{f}}{16} 4a_\mathsf{f} \geq \textrm{vol}(x^1) + \frac{\pi a_\mathsf{f}}{16} \sum_{i=1}^2 \Delta x^i,
 \end{align*}
since $\Delta x^i\leq 2a_\mathsf{f}$. Finally, we generalize for $\kappa>3$ by induction as follows
\begin{align*}
  \textrm{vol}(x^1,\ldots,x^\kappa) &\geq \textrm{vol}(x^1,\ldots,x^{\kappa-2}) +\textrm{vol}(x^\kappa)  \\
  &  \geq  \textrm{vol}(x^1)  + \Big(\frac{\pi a_\mathsf{f}}{16} \sum_{i=1}^{\kappa-3} \Delta x^i \Big)  + \frac{\pi a_\mathsf{f}^2}{4}\\
  & \geq \textrm{vol}(x^1)  + \frac{\pi a_\mathsf{f}}{16} \sum_{i=1}^{\kappa-1} \Delta x^i,
\end{align*}
which completes the proof. \carre

We prove now a lower bound on the volume of more than one member inside the sub-network $\mathbb{S}_{2a_\mathsf{f}}(x^1,\ldots,x^\kappa)$  
over $(x^1,\ldots,x^\kappa)\in \mathcal{C}_{\kappa 1}$ in which two members are connected (by affinity) if and only if their distance is at most 
$2a_\mathsf{f}$. 
\begin{lem}
 For any sequence $(x^1,\ldots,x^\kappa)\in\mathcal{C}_{\kappa 1}$ where $\kappa\geq 2$ and such 
 that $x^1$ being the one of the largest norm among $x^1,\ldots,x^\kappa$, we have
\begin{align}
 \emph{vol}(x^1,\ldots,x^\kappa) \geq \emph{vol}(x^1) + \frac{\pi a_\mathsf{f}}{16} \sup_{2\leq i\leq \kappa} \|x^i - x^1 \|.
 \label{volkappa}
\end{align}
\label{lemCk1}
\end{lem}

\proof  For the easy of the proof, we assume that $\|x^\kappa - x^1 \|=\sup_{2\leq i\leq \kappa} \|x^i - x^1 \|$ and let $\mathcal{P}$ be a minimum-hop path between 
$x^1$ and $x^\kappa$ in $\mathbb{S}_{2a_\mathsf{f}}(x^1,\ldots,x^\kappa)$ with $\Delta^\star x$ be the total length of $\mathcal{P}$. Thus, every pair of members in 
$\mathcal{P}$ that are not adjacent members in $\mathcal{P}$ are distant by more than $2a_\mathsf{f}$ and by application of Lemma \ref{lemsev} to the members in the 
path $\mathcal{P}$ we find that 
\begin{align*}
 \textrm{vol}\big(\{x^i:x^i\in \mathcal{P}\}\big) \geq \textrm{vol}(x^1) + \frac{\pi a_\mathsf{f}}{16} \Delta^\star x.
\end{align*}
We conclude by remarking that $\textrm{vol}(x^1,\ldots,x^\kappa) \geq \textrm{vol}(\{x^i:x^i\in \mathcal{P}\})$ and $\Delta^\star x \geq \|x^\kappa - x^1 \|$. This 
completes the proof. \carre 

In proving (\ref{beta}), we shall use the following lemma. 
\begin{lem}
 Assume that hypothesis $H_1$ and $H_2$ are satisfied. Let 
 \begin{align*}
  a_\mathsf{f}^\star(s_{T_k})= \Big(\frac{ \log\Big(\frac{s_{T_k} \pi}{4(s_{T_k}-1)!}\Big)+\mathfrak{f}_k({s_{T_k}})+ (s_{T_k}-1)\log(\alpha_k^\star(s_{T_k}))-
 \alpha_k^\star(s_{T_k}) }{\pi p s_{T_k}}\Big)^{1/2},
 \end{align*} 
 for any $p\in(0,1)$ and where $\alpha_k^\star$ is given by (\ref{alphastar}). Then, we have  
 \begin{align*}
  s_{T_k}  \Big(\frac{\big(\alpha_k^\star(s_{T_k}) \big)^{s_{T_k}-1}}{(s_{T_k}-1)!} \emph{e}^{-\big(\alpha_k^\star(s_{T_k}) \big) }  \Big)  
  \int_{[0,1]^2} \emph{e}^{-s_{T_k} p \emph{vol}(x)}dx  \longrightarrow     \emph{e} ^{-\displaystyle\lim_{s_{T_k}\to \infty}\mathfrak{f}_k({s_{T_k}})}, \qquad \textrm{ as } 
   s_{T_k}\to \infty.
 \end{align*}
 \label{lempois}
\end{lem}

\proof We shall proceed by approximating the integral in the left-hand side using the four regions $\mathcal{D}_1,\ldots,\mathcal{D}_4$ of the unit square $[0,1]^2$. If 
$x\in\mathcal{D}_1$, we know that $\textrm{vol}(x)=\pi (a_\mathsf{f}^\star)^2$ and it's straightforward that 
\begin{align*}
  s_{T_k}  \Big(\frac{\big(\alpha_k^\star(s_{T_k}) \big)^{s_{T_k}-1}}{(s_{T_k}-1)!} \textrm{e}^{-\big(\alpha_k^\star(s_{T_k}) \big) }  \Big)  
  \int_{\mathcal{D}_1} \textrm{e}^{-s_{T_k} p \textrm{vol}(x)}dx  &=
  s_{T_k}  \textrm{e}^{-\log\Big(\frac{s_{T_k} \pi}{4}\Big)-\mathfrak{f}_k({s_{T_k}})}  \int_{\mathcal{D}_1}dx \\
  &= \frac{4}{\pi}\textrm{e}^{-\mathfrak{f}_k({s_{T_k}})} \pi\Big(\frac{1}{2} - a_\mathsf{f}^\star \Big)^2 \\
  &= \textrm{e}^{-\mathfrak{f}_k({s_{T_k}})} \Big(1 - 2 a_\mathsf{f}^\star \Big)^2 \\
  &  {\stackrel{{{\scriptscriptstyle s_{T_k}\to+\infty}}}{\longrightarrow}}      \textrm{e} ^{-\displaystyle\lim_{s_{T_k}\to \infty}\mathfrak{f}_k({s_{T_k}})},
 \end{align*}
where we note that for sufficiently large  $s_{T_k}$ we have  $a_\mathsf{f}^\star \to 0$ (by hypothesis $H_1$). Therefore, if $x\in\mathcal{D}_3$, we know that 
$\textrm{vol}(x) \geq \frac{1}{2}\pi (a_\mathsf{f}^\star)^2$ and by using the upper  bound 
\begin{align*}
 \forall k\geq 1, \qquad \Big(\frac{\big(\alpha_k^\star(s_{T_k}) \big)^{s_{T_k}-1}}{(s_{T_k}-1)!} \textrm{e}^{-\big(\alpha_k^\star(s_{T_k}) \big) }  \Big) \leq 1,
\end{align*}
we find
\begin{align*}
 s_{T_k}  \Big(\frac{\big(\alpha_k^\star(s_{T_k}) \big)^{s_{T_k}-1}}{(s_{T_k}-1)!} \textrm{e}^{-\big(\alpha_k^\star(s_{T_k}) \big) }  \Big)  
  \int_{\mathcal{D}_3} \textrm{e}^{-s_{T_k} p \textrm{vol}(x)}dx  & \leq 
 s_{T_k}  \textrm{e}^{-\frac{1}{2}  s_{T_k} p \pi (a_\mathsf{f}^\star)^2}  \int_{\mathcal{D}_3}dx \\
  &= s_{T_k}  \pi (a_\mathsf{f}^\star)^2 \textrm{e}^{-\frac{1}{2}  s_{T_k} p \pi (a_\mathsf{f}^\star)^2}
  {\stackrel{{{\scriptscriptstyle s_{T_k}\to+\infty}}}{\longrightarrow}} 0.
\end{align*}
The same think happens if $x\in\mathcal{D}_4$ where $\textrm{vol}(x) \geq \frac{1}{4}\pi (a_\mathsf{f}^\star)^2$ but we need more notation. Let consider 
the points $A=(0,0), B=(\frac{1}{2},0),C=(0,\frac{1}{2})$ and the triangle $\mathbb{T}$ formed by $ABC$. From Figure \ref{partition}, it's clear that 
$\textrm{vol}(\mathcal{D}_4)\leq 4 \textrm{vol}(\mathbb{T})=\frac{4}{8}$ and we find 
\begin{align*}
 s_{T_k}  \Big(\frac{\big(\alpha_k^\star(s_{T_k}) \big)^{s_{T_k}-1}}{(s_{T_k}-1)!} \textrm{e}^{-\big(\alpha_k^\star(s_{T_k}) \big) }  \Big)  
  \int_{\mathcal{D}_4} \textrm{e}^{-s_{T_k} p \textrm{vol}(x)}dx  & \leq 
  s_{T_k}  \textrm{e}^{-\frac{1}{4}  s_{T_k} p \pi (a_\mathsf{f}^\star)^2} \int_{\mathcal{D}_4}dx \\
  & \leq  \frac{s_{T_k}}{2} \textrm{e}^{-\frac{1}{4}  s_{T_k} p \pi (a_\mathsf{f}^\star)^2} 
   {\stackrel{{{\scriptscriptstyle s_{T_k}\to+\infty}}}{\longrightarrow}} 0.
\end{align*}
Let us now turn out to the region  $x\in\mathcal{D}_2$. Hence, by application of Lemma \ref{lemvol}  and the polar coordinate system we have
\begin{align*}
 s_{T_k} \Big(\frac{\big(\alpha_k^\star(s_{T_k})  \big)^{s_{T_k}-1}}{(s_{T_k}-1)!}& \textrm{e}^{-\big(\alpha_k^\star(s_{T_k}) \big) }  \Big)  
  \int_{\mathcal{D}_2} \textrm{e}^{-s_{T_k} p \textrm{vol}(x)}dx \\  
  & \leq s_{T_k}   \textrm{e}^{-\frac{1}{2}s_{T_k} p \pi (a_\mathsf{f}^\star)^2} \int_{\mathcal{D}_2} \textrm{e}^{ - s_{T_k} p  a_\mathsf{f}^\star \Big( \frac{1}{2}-\eta \|x\| \Big)}dx  \\
  &= 2\pi s_{T_k}   \textrm{e}^{-\frac{1}{2}s_{T_k} p \pi (a_\mathsf{f}^\star)^2} \int_{\frac{1}{2}-a_\mathsf{f}^\star}^{\sqrt{\frac{1}{4}-(a_\mathsf{f}^\star)^2}}
  \rho_x  \textrm{e}^{ - s_{T_k} p  a_\mathsf{f}^\star \Big( \frac{1}{2}-\eta \rho_x \Big)} d\rho_x \\
  & \leq 2\sqrt{2}\pi s_{T_k}   \textrm{e}^{-\frac{1}{2}s_{T_k} p \pi (a_\mathsf{f}^\star)^2} \int_{\frac{1}{2}-a_\mathsf{f}^\star}^{\frac{1}{2}}  
  \textrm{e}^{ - s_{T_k} p  a_\mathsf{f}^\star \Big( \frac{1}{2}-\eta \rho_x \Big)} d\rho_x \\
  &= \frac{2\sqrt{2}\pi}{\eta } s_{T_k} \textrm{e}^{-\frac{1}{2}s_{T_k} p \pi (a_\mathsf{f}^\star)^2}
  \int_{\frac{1}{2}-\frac{1}{2}\eta}^{\frac{1}{2}-\eta(\frac{1}{2}-a_\mathsf{f}^\star)}  
  \textrm{e}^{ - s_{T_k} p  a_\mathsf{f}^\star z} dz \\
  & = \frac{2\sqrt{2}\pi}{\eta p a_\mathsf{f}^\star} \textrm{e}^{-\frac{1}{2}s_{T_k} p \pi (a_\mathsf{f}^\star)^2} 
  \Big(1-\textrm{e}^{-s_{T_k} p \eta (a_\mathsf{f}^\star)^2} \Big) \\
  & \leq  \frac{2\sqrt{2}\pi}{ p }   \frac{1}{a_\mathsf{f}^\star} \Big( \frac{\sqrt{\frac{1}{2}-2(a_\mathsf{f}^\star)^2}}{\frac{1}{2}-2(a_\mathsf{f}^\star)^2} -1\Big)
  \textrm{e}^{-\frac{1}{2}s_{T_k} p \pi (a_\mathsf{f}^\star)^2} =\frac{O(1)}{\sqrt{\psi_k(s_{T_k})}}= o(1),
\end{align*}
by our choice of $\mathfrak{f}_k$, this shows $\psi_k$ tends to infinity, completing the proof. \carre

\begin{lem}
 Assume that hypothesis $H_1$ and $H_2$ are satisfied. Let 
 \begin{align*}
  a_\mathsf{f}^\star(s_{T_k})= \Big(\frac{ \log\Big(\frac{s_{T_k} \pi}{4(s_{T_k}-1)!}\Big)+\mathfrak{f}_k({s_{T_k}})+ (s_{T_k}-1)\log(\alpha_k^\star(s_{T_k}))-
 \alpha_k^\star(s_{T_k}) }{\pi p s_{T_k}}\Big)^{1/2},
 \end{align*} 
 for any $p\in(0,1)$ and where $\alpha_k^\star$ is given by (\ref{alphastar})  and satisfying (\ref{constraint}). Then, for any $\kappa\geq 2$ and 
 $(x^1,\ldots,x^\kappa)\in\mathcal{C}_{\kappa \mathbf{n}}$ with $1\leq \mathbf{n}<\kappa$, we have  
 \begin{align*}
  (s_{T_k})^\kappa  \Big(\frac{\big(\alpha_k^\star(s_{T_k}) \big)^{s_{T_k}-\kappa}}{(s_{T_k}-\kappa)!} \emph{e}^{-\big(\alpha_k^\star(s_{T_k}) \big) }  \Big)  
  \int_{\mathcal{C}_{\kappa \mathbf{n}}} \emph{e}^{-s_{T_k} p \emph{vol}(x^1,\ldots,x^\kappa)}dx^1\cdots dx^\kappa  \longrightarrow   0, \qquad \textrm{ as } 
   s_{T_k}\to \infty.
 \label{beta0}
 \end{align*}
 \label{sevto0}
\end{lem}

\proof We divide the proof into two steps. \\ 

\paragraph{Step 1.} Let us first show the lemma for $\mathbf{n}=1$.   For simplicity and without loss of generality, let consider the subset $\mathcal{C}_{\kappa 1}^0$ 
denoting the set of 
$(x^1,\ldots,x^\kappa)\in\mathcal{C}_{\kappa 1}$ satisfying that $x^1$ being the one of the largest norm among $x^1,\ldots,x^\kappa$ and $x^\kappa$ 
being the one with longest distance from $x^1$  among $x^1,\ldots,x^\kappa$ which enables us to write
 \begin{align*}
  (s_{T_k})^\kappa & \Big(\frac{\big(\alpha_k^\star(s_{T_k}) \big)^{s_{T_k}-\kappa}}{(s_{T_k}-\kappa)!} \textrm{e}^{-\big(\alpha_k^\star(s_{T_k}) \big) }  \Big)  
  \int_{\mathcal{C}_{\kappa 1}} \textrm{e}^{-s_{T_k} p \textrm{vol}(x^1,\ldots,x^\kappa)}dx^1\cdots dx^\kappa  \\
  &\leq \kappa (\kappa-1)(s_{T_k})^\kappa  \Big(\frac{\big(\alpha_k^\star(s_{T_k}) \big)^{s_{T_k}-\kappa}}{(s_{T_k}-\kappa)!} \textrm{e}^{-\big(\alpha_k^\star(s_{T_k}) \big) }  \Big)  
  \int_{\mathcal{C}_{\kappa 1}^0} \textrm{e}^{-s_{T_k} p \textrm{vol}(x^1,\ldots,x^\kappa)}dx^1\cdots dx^\kappa.
  \end{align*}
We claim next that (\ref{volkappa}) (in Lemma \ref{lemCk1}) holds also for $(x^1,\ldots,x^\kappa)\in\mathcal{C}_{\kappa 1}^0$ with $\frac{\pi}{16}$ replaced by some 
constant $\tilde{C}$, that is,
\begin{align}
 \textrm{vol}(x^1,\ldots,x^\kappa) \geq \textrm{vol}(x^1) + \tilde{C} a_\mathsf{f}^\star  \|x^\kappa - x^1 \|,
\end{align}
with $x^\kappa\in B(x^1,2(\kappa-1)a_\mathsf{f}^\star)$ and, for $i\in\{2,\ldots,\kappa-1\}$, $x^i\in B(x^1,\|x^\kappa - x^1 \|)$. Indeed, by the constraint 
(\ref{constraint}) imposed on $\alpha^\star$ and the fact that 
\begin{align*}
\forall k\geq 1,\qquad \Big(\frac{\big(\alpha_k^\star(s_{T_k}) \big)^{s_{T_k}-1}}{(s_{T_k}-1)!} \emph{e}^{-\big(\alpha_k^\star(s_{T_k}) \big) } 
\Big)^{(\kappa-1)} \leq 1,
\end{align*}
we obtain 
\begin{align*}
 (s_{T_k})^\kappa & \Big(\frac{\big(\alpha_k^\star(s_{T_k}) \big)^{s_{T_k}-\kappa}}{(s_{T_k}-\kappa)!} \textrm{e}^{-\big(\alpha_k^\star(s_{T_k}) \big) }  \Big)  
  \int_{\mathcal{C}_{\kappa 1}^0} \textrm{e}^{-s_{T_k} p \textrm{vol}(x^1,\ldots,x^\kappa)}dx^1\cdots dx^\kappa \\
  & \leq (s_{T_k})^\kappa  \Big(\frac{\big(\alpha_k^\star(s_{T_k}) \big)^{s_{T_k}-1}}{(s_{T_k}-1)!} \textrm{e}^{-\big(\alpha_k^\star(s_{T_k}) \big) }  \Big)  
  \int_{\mathcal{C}_{\kappa 1}^0} \textrm{e}^{-s_{T_k} p \Big(\textrm{vol}(x^1) + \tilde{C} a_\mathsf{f}^\star  \|x^\kappa - x^1 \| \Big)}dx^1\cdots dx^\kappa \\
  & \leq (s_{T_k})^\kappa  \Big(\frac{\big(\alpha_k^\star(s_{T_k}) \big)^{s_{T_k}-1}}{(s_{T_k}-1)!} \textrm{e}^{-\big(\alpha_k^\star(s_{T_k}) \big) }  \Big) 
  \int_{[0,1]^2} \textrm{e}^{-s_{T_k} p \textrm{vol}(x^1)} dx^1 \\
  &\qquad \times \int_{ B(x^1,2(\kappa-1)a_\mathsf{f}^\star)} \textrm{e}^{-s_{T_k} p \tilde{C} a_\mathsf{f}^\star  \|x^\kappa - x^1 \|} dx^\kappa
  \prod_{i=2}^{\kappa-1} \int_{B(x^1,\|x^\kappa - x^1 \|) } dx^i \\
 & = 2(\pi s_{T_k})^{\kappa-1} \Big(s_{T_k}\Big(\frac{\big(\alpha_k^\star(s_{T_k}) \big)^{s_{T_k}-1}}{(s_{T_k}-1)!} \textrm{e}^{-\big(\alpha_k^\star(s_{T_k}) \big) }  \Big) 
  \int_{[0,1]^2} \textrm{e}^{-s_{T_k} p \textrm{vol}(x^1)} dx^1 \Big) \\
  &\qquad \times \Big( \int_{ 0}^{2(\kappa-1)a_\mathsf{f}^\star} (\rho_x)^{2\kappa-3} \textrm{e}^{-s_{T_k} p \tilde{C} a_\mathsf{f}^\star 
  \rho_x} d\rho_x \Big ) 
  \end{align*}
  \begin{align*}
  & < 2(\pi s_{T_k})^{\kappa-1} \Big(s_{T_k}\Big(\frac{\big(\alpha_k^\star(s_{T_k}) \big)^{s_{T_k}-1}}{(s_{T_k}-1)!} \textrm{e}^{-\big(\alpha_k^\star(s_{T_k}) \big) }  \Big) 
  \int_{[0,1]^2} \textrm{e}^{-s_{T_k} p \textrm{vol}(x^1)} dx^1 \Big) \\
  &\qquad \times \Big( \int_{ 0}^{\infty} (\rho_x)^{2\kappa-3} \textrm{e}^{-s_{T_k} p \tilde{C} a_\mathsf{f}^\star 
  \rho_x} d\rho_x \Big ) \\
  &= \frac{\Gamma(2\kappa-2)2(\pi s_{T_k})^{\kappa-1}}{(s_{T_k} p \tilde{C} a_\mathsf{f}^\star)^{2\kappa-2}}
  \Big(s_{T_k}\Big(\frac{\big(\alpha_k^\star(s_{T_k}) \big)^{s_{T_k}-1}}{(s_{T_k}-1)!} \textrm{e}^{-\big(\alpha_k^\star(s_{T_k}) \big) }  \Big) 
  \int_{[0,1]^2} \textrm{e}^{-s_{T_k} p \textrm{vol}(x^1)} dx^1 \Big) \\
  & = O(1) \frac{s_{T_k}\Big(\frac{\big(\alpha_k^\star(s_{T_k}) \big)^{s_{T_k}-1}}{(s_{T_k}-1)!} \textrm{e}^{-\big(\alpha_k^\star(s_{T_k}) \big) }  \Big) 
  \int_{[0,1]^2} \textrm{e}^{-s_{T_k} p \textrm{vol}(x^1)} dx^1}{(\psi_k(s_{T_k}))^{\kappa-1}}=o(1),
\end{align*}
by application of Lemma \ref{lempois} and where $\Gamma(\kappa)=\int_0^\infty x^{\kappa-1} \textrm{e}^{-x}dx$ denoting the gamma function. 

\paragraph{Step 2.} We show now that the same result holds for $(x^1,\ldots,x^\kappa)\in\mathcal{C}_{\kappa \mathbf{n}}$ for any $2\leq \mathbf{n}<\kappa$. For any 
random $\mathbf{n}-$partition 
\begin{align*}
 \Pi_\kappa=\{P_1,\ldots,P_\mathbf{n}\} \qquad \textrm{ of the subset } \; [\kappa]:=\{1,\ldots,\kappa\},
\end{align*}
where each component $P_j$, $j=1,\ldots,\mathbf{n}$, is of cardinal $|P_j|$ and let denote by $\bar{\mathcal{D}}^\kappa(\Pi_\kappa)$ the set of 
$(x^1,\ldots,x^\kappa)\in\bar{\mathcal{D}}^\kappa$ such that the points $\{x^i: i \in P_j\}$ formed a connected component of 
$\mathbb{S}_{2a_\mathsf{f}^\star}(x^1,\ldots,x^\kappa)$. Hence, 
\begin{align*}
 \mathcal{C}_{\kappa \mathbf{n}} =\displaystyle \bigcup_{\small{\textrm{all }\mathbf{n}-\textrm{partitions }\Pi_\kappa}}\bar{\mathcal{D}}^\kappa(\Pi_\kappa),
\end{align*}
and it suffice to prove that for any $\mathbf{n}-$partition $\Pi_\kappa$
\begin{align*}
 (s_{T_k})^\kappa  \Big(\frac{\big(\alpha_k^\star(s_{T_k}) \big)^{s_{T_k}-\kappa}}{(s_{T_k}-\kappa)!} \textrm{e}^{-\big(\alpha_k^\star(s_{T_k}) \big) }  \Big)  
  \int_{\bar{\mathcal{D}}^\kappa(\Pi_\kappa)} \textrm{e}^{-s_{T_k} p \textrm{vol}(x^1,\ldots,x^\kappa)}dx^1\cdots dx^\kappa  \longrightarrow   0, \qquad \textrm{ as } 
   s_{T_k}\to \infty.
\end{align*}
Second, without loss of generality, let now fix one arbitrary  $\mathbf{n}-$partition $\Pi_\kappa$ and observe that for any 
$(x^1,\ldots,x^\kappa)\in\bar{\mathcal{D}}^\kappa(\Pi_\kappa)$, we have
\begin{align}
 \bar{\mathcal{D}}^\kappa(\Pi_\kappa) \subseteq \prod_{j=1}^\mathbf{n}  \mathcal{C}_{|P_j| 1} \qquad \textrm{ and } \qquad 
 \textrm{vol}(x^1,\ldots,x^\kappa) = \sum_{j=1}^{\mathbf{n}} \textrm{vol}(\{x^i: i\in P_j\}).
 \label{volclass}
\end{align}
It follows from (\ref{volclass}) that
\begin{align*}
 (s_{T_k})^\kappa & \Big(\frac{\big(\alpha_k^\star(s_{T_k}) \big)^{s_{T_k}-\kappa}}{(s_{T_k}-\kappa)!} \textrm{e}^{-\big(\alpha_k^\star(s_{T_k}) \big) }  \Big)  
  \int_{\bar{\mathcal{D}}^\kappa(\Pi_\kappa)} \textrm{e}^{-s_{T_k} p \textrm{vol}(x^1,\ldots,x^\kappa)}dx^1\cdots dx^\kappa \\
  & =  \Big(s_{T_k}\frac{\big(\alpha_k^\star(s_{T_k}) \big)^{s_{T_k}-1}}{(s_{T_k}-1)!} \textrm{e}^{-\big(\alpha_k^\star(s_{T_k}) \big) }  \Big)^{\kappa}  
  \int_{\bar{\mathcal{D}}^\kappa(\Pi_\kappa)} \prod_{j=1}^{\mathbf{n}} \textrm{e}^{-s_{T_k} p  \textrm{vol}(\{x^i: i\in P_j\})}dx^1\cdots dx^\kappa \\
  & \leq  \Big(s_{T_k}\frac{\big(\alpha_k^\star(s_{T_k}) \big)^{s_{T_k}-1}}{(s_{T_k}-1)!} \textrm{e}^{-\big(\alpha_k^\star(s_{T_k}) \big) }  \Big)^{\kappa}
  \prod_{j=1}^{\mathbf{n}}  \int_{\mathcal{C}_{|P_j| 1} }\textrm{e}^{-s_{T_k} p  \textrm{vol}(\{x^i: i\in P_j\})} \prod_{i\in P_j} dx^i\\
  & = \prod_{j=1}^{\mathbf{n}} \Big(s_{T_k}\frac{\big(\alpha_k^\star(s_{T_k}) \big)^{s_{T_k}-1}}{(s_{T_k}-1)!} \textrm{e}^{-\big(\alpha_k^\star(s_{T_k}) \big) }  \Big)^{|P_j|}
    \int_{\mathcal{C}_{|P_j| 1} }\textrm{e}^{-s_{T_k} p  \textrm{vol}(\{x^i: i\in P_j\})} \prod_{i\in P_j} dx^i,
\end{align*}
which tends to zero as shown in Step 1, completing the proof. \carre 

Next we study  the limit in $\mathcal{C}_{\kappa \kappa}$.   

\begin{lem}
 Assume that hypothesis $H_1$ and $H_2$ are satisfied. Let 
 \begin{align*}
  a_\mathsf{f}^\star(s_{T_k})= \Big(\frac{ \log\Big(\frac{s_{T_k} \pi}{4(s_{T_k}-1)!}\Big)+\mathfrak{f}_k({s_{T_k}})+ (s_{T_k}-1)\log(\alpha_k^\star(s_{T_k}))-
 \alpha_k^\star(s_{T_k}) }{\pi p s_{T_k}}\Big)^{1/2},
 \end{align*} 
 for any $p\in(0,1)$ and where $\alpha_k^\star$ is given by (\ref{alphastar})  and satisfying (\ref{constraint}). Then, for any $\kappa\geq 2$ and 
 $(x^1,\ldots,x^\kappa)\in\mathcal{C}_{\kappa \kappa}$, we have  
 \begin{align*}
  (s_{T_k})^\kappa  \Big(\frac{\big(\alpha_k^\star(s_{T_k}) \big)^{s_{T_k}-\kappa}}{(s_{T_k}-\kappa)!} \emph{e}^{-\big(\alpha_k^\star(s_{T_k}) \big) }  \Big)  
  \int_{\mathcal{C}_{\kappa \kappa}} \emph{e}^{-s_{T_k} p \emph{vol}(x^1,\ldots,x^\kappa)}dx^1\cdots dx^\kappa  
  {\stackrel{{{\scriptscriptstyle s_{T_k}\to+\infty}}}{\longrightarrow}}   
  \emph{e} ^{-\displaystyle \kappa \lim_{s_{T_k}\to \infty}\mathfrak{f}_k({s_{T_k}})}.
 \end{align*}
 \label{sevtono0}
\end{lem}

\proof Recall that for any  $(x^1,\ldots,x^\kappa)\in\mathcal{C}_{\kappa \kappa}$ we observe $\textrm{vol}(x^1,\ldots,x^\kappa) = \sum_{i=1}^{\kappa}\textrm{vol}(x^i)$. 
This enables to write
\begin{align}
\nonumber
 (s_{T_k})^\kappa  &\Big(\frac{\big(\alpha_k^\star(s_{T_k}) \big)^{s_{T_k}-\kappa}}{(s_{T_k}-\kappa)!} \textrm{e}^{-\big(\alpha_k^\star(s_{T_k}) \big) }  \Big)  
  \int_{\mathcal{C}_{\kappa \kappa}} \textrm{e}^{-s_{T_k} p \textrm{vol}(x^1,\ldots,x^\kappa)}dx^1\cdots dx^\kappa \\
  \nonumber
  & = (s_{T_k})^\kappa  \Big(\frac{\big(\alpha_k^\star(s_{T_k}) \big)^{s_{T_k}-\kappa}}{(s_{T_k}-\kappa)!} \textrm{e}^{-\big(\alpha_k^\star(s_{T_k}) \big) }  \Big)  
  \Big\{  \int_{\bar{\mathcal{D}}^\kappa} \textrm{e}^{-s_{T_k} p \sum_{i=1}^{\kappa}\textrm{vol}(x^i)}dx^1\cdots dx^\kappa \\ 
  \nonumber
  & \qquad - \int_{\bar{\mathcal{D}}^\kappa \setminus \mathcal{C}_{\kappa \kappa}} \textrm{e}^{-s_{T_k} p 
  \sum_{i=1}^{\kappa}\textrm{vol}(x^i)}dx^1\cdots dx^\kappa \Big\}  \\
  \nonumber
  & = \Big\{ \prod_{i=1}^\kappa  
    \Big(s_{T_k}\frac{\big(\alpha_k^\star(s_{T_k}) \big)^{s_{T_k}-1}}{(s_{T_k}-1)!} \textrm{e}^{-\big(\alpha_k^\star(s_{T_k}) \big) }  \Big)  
  \int_{[0,1]^2} \textrm{e}^{-s_{T_k} p \textrm{vol}(x^i)}dx^i  \Big\} \\
  \label{minus}
  & \qquad-(s_{T_k})^\kappa  \Big(\frac{\big(\alpha_k^\star(s_{T_k}) \big)^{s_{T_k}-\kappa}}{(s_{T_k}-\kappa)!} \textrm{e}^{-\big(\alpha_k^\star(s_{T_k}) \big) }  \Big) 
  \int_{\bar{\mathcal{D}}^\kappa \setminus \mathcal{C}_{\kappa \kappa}} \textrm{e}^{-s_{T_k} p \sum_{i=1}^{\kappa}\textrm{vol}(x^i)}dx^1\cdots dx^\kappa \\
  \nonumber
  & {\stackrel{{{\scriptscriptstyle s_{T_k}\to+\infty}}}{\longrightarrow}}   
  \Big(\textrm{e} ^{-\displaystyle\lim_{s_{T_k}\to \infty}\mathfrak{f}_k({s_{T_k}})} \Big)^\kappa -0.
\end{align}
There may be some doubt as to why the  term (\ref{minus}) tends to zero. Let us verify this by observing that for any  
$(x^1,\ldots,x^\kappa)\in \bar{\mathcal{D}}^\kappa \setminus\mathcal{C}_{\kappa \kappa}$ we have 
$\textrm{vol}(x^1,\ldots,x^\kappa) \leq  \sum_{i=1}^{\kappa}\textrm{vol}(x^i)$ which enables us to find 
\begin{align*}
 (s_{T_k})^\kappa & \Big(\frac{\big(\alpha_k^\star(s_{T_k}) \big)^{s_{T_k}-\kappa}}{(s_{T_k}-\kappa)!} \textrm{e}^{-\big(\alpha_k^\star(s_{T_k}) \big) }  \Big) 
  \int_{\bar{\mathcal{D}}^\kappa \setminus \mathcal{C}_{\kappa \kappa}} \textrm{e}^{-s_{T_k} p \sum_{i=1}^{\kappa}\textrm{vol}(x^i)}dx^1\cdots dx^\kappa \\
  & \leq (s_{T_k})^\kappa  \Big(\frac{\big(\alpha_k^\star(s_{T_k}) \big)^{s_{T_k}-\kappa}}{(s_{T_k}-\kappa)!} \textrm{e}^{-\big(\alpha_k^\star(s_{T_k}) \big) }  \Big)
  \int_{\bar{\mathcal{D}}^\kappa \setminus \mathcal{C}_{\kappa \kappa}} \textrm{e}^{-s_{T_k} p \textrm{vol}(x^1,\ldots,x^\kappa)}dx^1\cdots dx^\kappa \\
  &= \sum_{\mathbf{n}=1}^{\kappa-1} \Big\{ (s_{T_k})^\kappa  \Big(\frac{\big(\alpha_k^\star(s_{T_k}) \big)^{s_{T_k}-\kappa}}{(s_{T_k}-\kappa)!} \textrm{e}^{-\big(\alpha_k^\star(s_{T_k}) \big) }  \Big)
  \int_{\mathcal{C}_{\kappa \mathbf{n}}} \textrm{e}^{-s_{T_k} p \textrm{vol}(x^1,\ldots,x^\kappa)}dx^1\cdots dx^\kappa\Big\},
\end{align*}
where it's straightforward that the sum tends to zero thanks to Lemma \ref{sevto0}. \carre 

A key result  needed in the sequel for the proofs of  Theorems \ref{theopois} and \ref{theopoisac} is as follows:
\begin{lem}
 Given a sequence $E_{1k},E_{2k},\ldots,E_{s_{T_k} k}$ of events such that $E_{ik}$ be the event that the point $x_{T_k}^i$ is isolated, define $U_{s_{T_k}}$ to be 
 the random number of $E_{ik}$ that hold. If for any set $\{i_1,\ldots,i_\kappa\}$ it is true that 
 \begin{align}
  \mathbb{P}\Big(\bigcap_{j=1}^\kappa E_{jk}\Big) = \mathbb{P}\Big(\bigcap_{j=1}^\kappa E_{i_j k}\Big),
 \label{equivset}
 \end{align}
and there is a constant $u$ such that for any fixed $\kappa$
\begin{align}
 (s_{T_k})^\kappa \times  \mathbb{P}\Big(\bigcap_{j=1}^\kappa E_{j k}\Big) {\stackrel{{{\scriptscriptstyle s_{T_k}\to+\infty}}}{\longrightarrow}} u^\kappa,
\label{meanpois}
 \end{align}
hence the sequence $U_{s_{T_k}}$ converges in distribution to a Poisson random variable with mean $u$.
\label{bonfero}
\end{lem}
We do not claim originality of the Lemma \ref{bonfero} and, in fact, similar result have been proved in \cite{Alon+Spen00} using a probabilistic version of Brun's 
sieve theorem and the Bonferroni inequalities. Since we have found this particular result in the literature, we not provide a detailed proof.

We now proceed to conclude the proof of Theorem \ref{theopois}. 

\proof For $(x_{T_k}^1,\ldots,x_{T_k}^\kappa)\in\mathcal{C}_{\kappa\kappa}$, we proved in Proposition \ref{propisolkeq} that 
  \begin{align*}
   \mathbb{P}(E_{1k}\cap\cdots\cap E_{\kappa k}) &=  
   \Big(\frac{\big(\alpha_k^\star(s_{T_k}) \big)^{s_{T_k}-\kappa}}{(s_{T_k}-\kappa)!} \emph{e}^{-\big(\alpha_k^\star(s_{T_k}) \big) }  \Big)   
  \\
 &\qquad \qquad  \times \int_{\mathcal{C}_{\kappa\kappa}} \Big(1-p\textrm{vol}(x_{T_k}^1,\ldots,x_{T_k}^\kappa)\Big)^{s_{T_k}-\kappa}dx_{T_k}^1\cdots dx_{T_k}^\kappa.
  \end{align*}
  Or, for $s_{T_k}$ sufficiently large we find
  \begin{align*}
   \Big(1-p\textrm{vol}(x_{T_k}^1,\ldots,x_{T_k}^\kappa)\Big)^{s_{T_k}-\kappa} &=
   \frac{\Big(1-p\textrm{vol}(x_{T_k}^1,\ldots,x_{T_k}^\kappa)\Big)^{s_{T_k}}}{\Big(1-p \sum_{i=1}^\kappa \textrm{vol}(x_{T_k}^i)\Big)^{\kappa}} \\
   &= \frac{\Big(1-p\textrm{vol}(x_{T_k}^1,\ldots,x_{T_k}^\kappa)\Big)^{s_{T_k}}}{\Big(1-p \kappa \pi (a_\mathsf{f}^\star)^2)\Big)^{\kappa}} \\
   & \backsimeq \frac{\textrm{e}^{-s_{T_k} p \textrm{vol}(x_{T_k}^1,\ldots,x_{T_k}^\kappa)}}{1}.
  \end{align*}
Thus, by application of Lemma \ref{sevtono0}, 
\begin{align*}
 (s_{T_k})^\kappa \times  \mathbb{P}(E_{1k}\cap\cdots\cap E_{\kappa k}) {\stackrel{{{\scriptscriptstyle s_{T_k}\to+\infty}}}{\longrightarrow}}   
  \Big(\textrm{e} ^{-\displaystyle\lim_{s_{T_k}\to \infty}\mathfrak{f}_k({s_{T_k}})} \Big)^\kappa = (\beta_k)^\kappa <\infty,
\end{align*}
and immediately the condition (\ref{meanpois}) in Lemma \ref{bonfero} is verified. It is easily seen that the condition (\ref{equivset}) is also verified 
(its proof is left to the reader) and hence Theorem \ref{theopois} follows for $(x_{T_k}^1,\ldots,x_{T_k}^\kappa)\in\mathcal{C}_{\kappa\kappa}$ and generally for 
$(x_{T_k}^1,\ldots,x_{T_k}^\kappa)\in\bar{\mathcal{D}}^{\kappa}$ since for 
$(x_{T_k}^1,\ldots,x_{T_k}^\kappa)\in\bar{\mathcal{D}}^{\kappa} \setminus \mathcal{C}_{\kappa\kappa}$ and for $s_{T_k}$ sufficiently large 
  \begin{align*}
   \Big(1-p\textrm{vol}(x_{T_k}^1,\ldots,x_{T_k}^\kappa)\Big)^{s_{T_k}-\kappa} \leq 
    \frac{\Big(1-p\textrm{vol}(x_{T_k}^1,\ldots,x_{T_k}^\kappa)\Big)^{s_{T_k}}}{\Big(1-p \kappa \pi (a_\mathsf{f}^\star)^2)\Big)^{\kappa}} 
   \backsimeq \textrm{e}^{-s_{T_k} p \textrm{vol}(x_{T_k}^1,\ldots,x_{T_k}^\kappa)}.
  \end{align*}
Then, the probability of the event 
$(E_{1k}\cap E_{2k}\cap \cdots \cap E_{\kappa k})$ times $(s_{T_k})^\kappa$ tends to zero as $s_{T_k} \to \infty$ by Lemma \ref{sevto0}, this completes the proof.  \carre 

Finally, we conclude the proof of Theorem \ref{theopoisac}. 
\proof Note that 
\begin{align*}
 \mathbb{P}(F_{1k}\cap\cdots\cap F_{\kappa k}) = p^\kappa  \mathbb{P}(E_{1k}\cap\cdots\cap E_{\kappa k}), 
\end{align*}
and Theorem \ref{theopoisac} follows using the same arguments as  for Theorem \ref{theopois}. \carre


\bibliographystyle{livre} 
\bibliography{biblio}

\end{document}